\documentclass[preprint,3p,times]{elsarticle}

\usepackage{array}
\usepackage{color}
\usepackage{tabularx}
\usepackage{graphicx}
\usepackage{amsmath}
\usepackage{amssymb}
\usepackage{amsfonts}	
\usepackage{moreverb}
\usepackage{dsfont}
\usepackage{bm}
\usepackage{multirow}
\usepackage{soul}

\newcommand\BibTeX{{\rmfamily B\kern-.05em \textsc{i\kern-.025em b}\kern-.08em
T\kern-.1667em\lower.7ex\hbox{E}\kern-.125emX}}

\newcommand{\x}{\mbf{x}}

\newcommand{\mbf}[1]{\mathbf{#1}}			%

\renewcommand{\u}{\mathbf{u}}

\newcommand{\q}{\mathbf{q}}
\newcommand{\F}{\mathbf{F}}
\newcommand{\f}{\mathbf{f}}
\newcommand{\g}{\mathbf{g}}
\newcommand{\h}{\mathbf{h}}
\renewcommand{\v}{\mathbf{v}}
\newcommand{\B}{\mathbf{B}}

\newcommand{\Ell}{\mathcal{L}}
\newcommand{\emm}{m}

\newcommand{\PMM}{\mathbb{P}_N\mathbb{P}_N}

\newcommand{\halb}{\frac{1}{2}}

\newcommand{\be}{\begin{equation}}
\newcommand{\ee}{\end{equation}}
\newcommand{\bdm}{\begin{displaymath}}
\newcommand{\edm}{\end{displaymath}}

\newcommand{\apriori}{\textit{a priori} }
\newcommand{\aposteriori}{\textit{a posteriori} }

\newfont{\numerikEleven}{ecrm1000}
\newfont{\numerikTen}{cmss10}
\newfont{\numerikNine}{cmss9}
\newfont{\numerikEight}{cmss8}

\journal{Computers and Fluids}

\begin{document} 
\begin{frontmatter}
\title{Space-time adaptive ADER discontinuous Galerkin finite element schemes with a posteriori sub-cell finite volume limiting} 
\author[UniTN]{Olindo Zanotti}
\ead{olindo.zanotti@unitn.it}

\author[UniTN]{Francesco Fambri}
\ead{francesco.fambri@unitn.it}

\author[UniTN]{Michael Dumbser$^{*}$}
\ead{michael.dumbser@unitn.it}
\cortext[cor1]{Corresponding author}

\author[UPM]{Arturo Hidalgo}
\ead{arturo.hidalgo@upm.es}

\address[UniTN]{Department of Civil, Environmental and Mechanical Engineering, University of Trento, Via Mesiano, 77 - 38123 Trento, Italy.}
\address[UPM]{Departamento de Matem\'atica Aplicada y M\'etodos Inform\'aticos, Universidad Polit\'ecnica de Madrid, Madrid, Spain}

\begin{abstract}
In this paper we present a novel arbitrary high order accurate discontinuous Galerkin (DG) finite element method on space-time adaptive Cartesian meshes (AMR) 
for hyperbolic conservation laws in multiple space dimensions, using a high order \aposteriori sub-cell ADER-WENO finite volume \emph{limiter}.  
Notoriously, the original DG method produces strong oscillations in the presence of discontinuous solutions and several types of limiters have been 
introduced over the years to cope with this problem.
Following the innovative idea recently proposed in \cite{Dumbser2014}, the discrete solution within the troubled cells is \textit{recomputed} by 
scattering the DG polynomial at the previous time step onto a suitable number of sub-cells along each direction. Relying on the robustness of  
classical finite volume WENO schemes, the sub-cell averages are recomputed and then gathered back into the DG polynomials over the main grid. 
In this paper this approach is implemented for the first time within a space-time adaptive AMR framework in two and three space dimensions, after 
assuring the proper averaging and projection between sub-cells that belong to different levels of refinement.  The combination of the sub-cell 
resolution with the advantages of AMR allows for an unprecedented ability in resolving even the finest details in the dynamics of the fluid. 
The spectacular resolution properties of the new scheme have been shown through a wide number of test cases performed in two and in three space 
dimensions, both for the Euler equations of compressible gas dynamics and for the magnetohydrodynamics (MHD) equations.  
\end{abstract}

\end{frontmatter}

\section{Introduction} \label{sec:introduction}
The numerical solution of hyperbolic problems has attracted a lot of attention in recent years, as they arise in many physical and technological applications. Many of them are in the field of  computational fluid dynamics, such as compressible gas dynamics, multiphase flows, air flow around aircraft or cars, astrophysical flows, free surface flows, environmental and geophysical flows like avalanches, dam break problems and water flow in channels, rivers and oceans, to mention but a few. 
Among the numerical methods specifically developed to solve hyperbolic problems, there are finite volume (FV)  methods and discontinuous Galerkin (DG) methods.  
While until a few years ago FV methods were comparatively more popular, the situation is now rapidly changing and DG schemes, 
first introduced by Reed and Hill in \cite{reed} to solve a first order neutron transport equation,
are now widely applied in several  different fields, in particular those related to fluid dynamics. 
In a series of masterpiece works \cite{cbs0,cbs1,cbs2,cbs3,cbs4}, Cockburn and Shu provided a rigorous formal framework of these methods, contributing significantly to their widespread use. 
DG methods are very robust and, among high order numerical methods, they show high flexibility and adaptivity strategies in handling complex geometries \cite{Biswas_94}.
Moreover, Jiang and Shu proved in \cite{jiangshu} that DG methods verify an entropy condition which confers them nonlinear $L_2$ stability. Despite this interesting property, explicit DG methods have a  strong stability limitation, since usually the CFL restriction for these schemes is very severe and the time step in $d$ space dimensions is constrained as $\Delta t\leq h/[d(2N+1)|\lambda_{max}|]$, 
where $d$ is the number of space dimensions, $h$ is a characteristic mesh size, $\lambda_{max}$ is the maximum signal velocity and $N$ is the degree of the 
basis polynomial. 

In DG schemes a high order time integration is typically performed by means of TVD Runge-Kutta schemes \cite{shu2}, leading to the family of so-called RKDG schemes. 
These methods are certainly efficient, but they have a maximum reachable order of accuracy in time, which is four. However, due to the high complexity of the fourth order TVD Runge-Kutta scheme, only
up to third order TVD Runge-Kutta methods are used in practice. In the presence of stiff source terms, usually the so-called IMEX Runge-Kutta schemes are employed, see 
\cite{pareschi_2005_ier,Pareschi2005}. 
To overcome these limitations, in our approach we follow the so-called ADER strategy for time integration, which was first introduced by Toro and Titarev 
in the finite volume context \cite{toro1,toro3,toro4,titarevtoro,Toro:2006a}, and it is a very attractive tool allowing to achieve arbitrary order of accuracy in \textit{space and time} 
in one single step by incorporating the approximate solution of a Generalised Riemann Problem (GRP) at the element interfaces. There are essentially two different families of approximate 
GRP solvers: those who first interact the spatial derivatives and subsequently compute a temporal expansion at the interface \cite{Artzi,Raviart.GRP.1,Raviart.GRP.2,toro1,toro3,toro4,titarevtoro,Toro:2006a,dumbser_jsc}, and those who first evolve the data locally \textit{in the small} inside each element and then interact the evolved data at the element interfaces via a classical Riemann solver, 
see e.g. \cite{eno,stedg1,stedg2,DumbserEnauxToro,DumbserZanotti,Dumbser2008}. For a more detailed discussion on the approximate solution of the GRP, see \cite{CastroToro,Montecinos2012,GoetzIske}. 
Nevertheless, the  original ADER approach has two main drawbacks: first, it makes use of the rather cumbersome and problem-dependent Cauchy-Kowalewski procedure and, second, it fails in the presence of stiff source terms.  
A subsequent version of the ADER approach that solves both these difficulties was developed in \cite{DumbserEnauxToro}, where 
the Cauchy-Kowalewski procedure was replaced with a local space-time DG predictor approach based on a weak formulation of the problem in space-time. 
This formulation is usually referred to as the \emph{local space-time discontinuous Galerkin} (LSTDG) \emph{predictor} and it has been successfully adopted in a variety of mathematical and physical problems \cite{DumbserZanotti,HidalgoDumbser,DumbserNSE,AMR3DCL,AMR3DNC,BalsaraMultiDRS}. We remark that, although this LSTDG approach is locally implicit, the full formulation remains explicit and, therefore, the above mentioned CFL restriction still holds. The ADER time stepping method has been also applied successfully to the discontinuous Galerkin finite element framework, see e.g. 
\cite{dumbser_jsc,QiuDumbserShu,LTS}.  

DG schemes are very efficient in smooth regions, but in the presence of sharp gradients and/or shock waves, they cannot escape from the Gibbs phenomenon and, as a consequence, they give rise to undesirable oscillations in the solution, since they are linear in the sense of Godunov. 
In fact, according to Godunov's theorem \citep{godunov} there are no \textit{linear} and monotone schemes of order higher than the first. In the finite volume framework Godunov's theorem is circumvented by carrying out a \textit{nonlinear} reconstruction within each cell. Here, TVD slope limiters \cite{toro-book} and ENO/WENO reconstructions \cite{eno,shu_efficient_weno,balsarashu,FVWENO3D} 
are among the most popular.    
In the discontinuous Galerkin approach, on the other hand, even if in principle no spatial reconstruction is needed, in practice 
it is necessary to introduce some sort of \textit{limiters} to avoid oscillations in the presence of discontinuities. 
Among the most relevant limiters proposed so far 
we mention the use of artificial viscosity \cite{Hartman_02,Persson_06,Feistauer4,Feistauer5,Feistauer6,Feistauer7}, of spectral 
filtering \citep{Radice2011}, of (H)WENO limiting procedures \cite{QiuShu1,QiuShu3,Qiu_2004,balsara2007,Zhu_2008,Zhu_13,Luo_2007}, 
and of slope and moment limiting \cite{cbs1,Krivodonova_2007,Biswas_94,Burbeau_2001,Yang_parameterfree_09,Kuzmin2014}. 
In \cite{Dumbser2014} we have recently proposed a totally different and alternative solution to this longstanding problem, which relies on a new 
\textit{a posteriori} sub-cell finite volume limiting approach. In practice, we first 
compute the solution by means of an \textit{unlimited} DG scheme, and subsequently we find the \textit{troubled cells} by using some very simple but 
effective \aposteriori detection criteria, namely the positivity of the solution and a relaxed discrete maximum principle in the sense of polynomials. 
Once the troubled cells have been identified, a sub-grid of size $(2N+1)^d$ is created within these cells and a more robust ADER-WENO finite volume 
approach is used to recompute the solution on the subgrid. A peculiar aspect of this new paradigm is that the size of the sub-grid is chosen as to make sure 
that the maximum admissible time step of the finite volume scheme on the sub-cells matches the time step of the DG scheme on the main grid. 
The idea of introducing an \textit{a posteriori} approach to the problem of limiting has been recently established by Clain, Diot and Loub{\`e}re 
in the finite volume context, by means of the so-called Multi-dimensional Optimal Order Detection (MOOD) method \cite{CDL1,CDL2,CDL3,ADER_MOOD_14}. 
The MOOD paradigm may in fact be considered as the progenitor of our \textit{a posteriori} limiting procedure for DG schemes \cite{Dumbser2014}. 
 
In the present work we combine the new ADER-DG paradigm with \aposteriori subcell limiters of \cite{Dumbser2014} with 
Adaptive Mesh Refinement (AMR) techniques, thus significantly enhancing the resolution capabilities compared to simple uniform grids. AMR  was first proposed by M.~J. Berger and collaborators in a series of well-known papers \cite{Berger-Oliger1984,Berger-Colella1989,BergerLeveque2011}. They introduced a patch-boxed block-structured AMR approach developed within the finite volume framework and later used  extensively for astrophysical applications. Other interesting developments have been presented in \cite{Mulet1}, where the first higher order AMR algorithms based on WENO finite volume schemes 
have been introduced. 
Applications of AMR techniques in the field of shallow water equations have been reported,  for instance, by \cite{Donat2014}. Other interesting implementations of adaptive mesh refinement are based on the so-called quadtree/octree refinement; for an overview of these techniques in different contexts see 
\cite{agbaghla2011} and \cite{tsai2013}. 
AMR techniques have also been successfully implemented with Central WENO (CWENO) schemes, such as in \cite{IvanGroth2009,Ivan2014}. 
Following the "cell-by-cell"' refinement approach of \cite{Khokhlov98}, 
the first implementation of a high order ADER-WENO finite volume scheme with AMR was proposed in \cite{AMR3DCL,AMR3DNC} for conservative
and nonconservative hyperbolic PDE in two and three space dimensions. It was subsequently extended to special relativistic hydrodynamics 
(RHD) and magnetohydrodynamics (RMHD) in \cite{Zanotti2013}. 

The combination of DG schemes with AMR has been considered in a significant number of papers, although in this context
the concept of adaptive mesh refinement is commonly absorbed into that of \emph{hp-adaptivity}. 
Two well-known early series of papers on hp-adaptive DG schemes are due to Baumann and Oden \cite{Baumann1999311,Baumann199979} and 
Houston, S\"uli and Schwab, see \cite{Houston20022133,Houston20001618,Houston20021226}. 
Furthermore, in \cite{Leicht2008} a DG scheme was proposed with anisotropic AMR for the compressible Navier--Stokes equations, while in  
\cite{LuoBaumLoehner,Yu2011} the Euler equations have been solved on adaptive unstructured meshes.  
In the context of atmospheric simulations, on the other hand, 
\cite{Kopera2014} implemented a numerical scheme which includes implicit-explicit RKDG, artificial viscosity and adaptive mesh refinement on two dimensional 
non-conforming elements. Other relevant results have been obtained in \cite{Georgoulis2009}, \cite{Lu2014}. 
Our goal is to improve with respect to these approaches by proposing a space-time adaptive ADER-DG scheme with time-accurate local time stepping 
that can be arbitrarily high order accurate both in space and time, that avoids Runge--Kutta sub-steps as well as artificial viscosity of any kind, 
and that incorporates a proper \aposteriori subcell limiter within the full advantages of AMR. 

The plan of the paper is the following: in Section~\ref{sec:basic} 
we present the basic mathematical framework, while in Section~\ref{sec:ADER-DG} 
we explain the ADER discontinuous Galerkin method.  
Section~\ref{sec:A-posteriori-sub-cell-limiter} is devoted to the description of the \aposteriori sub-cell limiter, 
whereas the incorporation within the AMR framework is deferred to Section~\ref{sec:AMR}. The numerical results are discussed in Section \ref{sec:numerics} and 
the conclusions are given in Section \ref{sec:conclusion}. 

\section{Mathematical framework: an overview}
\label{sec:basic}

We consider nonlinear systems of hyperbolic equations written in conservative form as
\begin{equation}
\label{eqn.pde.nc}
\frac{\partial \u}{\partial t} + \nabla \cdot \F\left(\u\right) = 0,  
\qquad \x \in \Omega \subset \mathds{R}^d, \quad t \in \mathds{R}_0^+,
\end{equation}
where $\u$ is the vector of so-called conserved quantities, while
$\F(\u) = (\f,\g, \h)$ is a non-linear flux tensor that depends on the state $\u$. 
The computational domain $\Omega$ is discretized by a Cartesian grid composed by  
elements $T_i$, namely
\begin{equation}
\label{eqn.tetdef}
 \Omega = \bigcup \limits_{i=1}^{N_E} T_i\,,
\end{equation}
where the index $i$ ranges from 1 to the total number of elements $N_E$, which, in
our adaptively mesh refinement framework, is of course a time-dependent quantity.
In the following, we denote the cell volume by $|T_i| = \int_{T_i} d\x$.
At the beginning of each time-step, the numerical solution of Eq.~ \eqref{eqn.pde.nc} 
is represented within each cell $T_i$ 
by piecewise polynomials of maximum degree $N \geq 0$ as
\begin{equation}
\label{eqn.ansatz.uh}
  \u_h(\x,t^n) = \sum_{l=0}^{N}\Phi_l(\x) \hat{\u}^n_l= \Phi_l(\x) \hat{\u}^n_l \quad \x \in T_i\,,
\end{equation}
where $\u_h$ is referred to as the \emph{discrete representation} of the solution, 
while the coefficients
$\hat{\u}^n_l$ are usually called the \emph{degrees of freedom}.\footnote{Throughout this paper we use the Einstein summation convention, 
implying summation over indices appearing twice, although there is no need to distinguish among covariant and contra-variant indices.}
In the expansion expressed by Eq.~(\ref{eqn.ansatz.uh}), 
the basis functions $\Phi_l(\x)$ are chosen as tensor-products of 
Lagrange interpolation polynomials of maximum degree $N$ 
which pass through the tensor-product of $(N+1)$ Gauss-Legendre 
quadrature points \citep{stroud,Kopriva,GassnerKopriva,KoprivaGassner}. \\
The numerical method used in this paper is the combination of several crucial steps, which will be described below and that can be listed schematically  as
\begin{itemize}
\item \emph{a predictor step}, 
in which Eq.~\eqref{eqn.pde.nc} is solved within each element \textit{in the small} \citep{eno} 
by means of a locally implicit space-time discontinuous Galerkin scheme (see Sect.~\ref{sec.galerkin});
\item \emph{a pure discontinuous Galerkin (DG) scheme}, i.e. a $\PMM$ scheme according to \cite{Dumbser2008}, which,
by exploiting the information obtained by the predictor, allows to compute the solution at the next time level through a single one-step corrector 
(see Sect.~\ref{sec.ADERNC});
\item \emph{an a posteriori sub-cell limiter}, which recomputes the solution of the troubled cells needing a limiter through an ADER-WENO finite volume scheme acting at the sub-cell level (see Sect.~\ref{sec:A-posteriori-sub-cell-limiter});
\item \emph{an adaptive mesh refinement} (AMR) approach, which is implemented according to a cell-by-cell strategy  and must be properly nested within the sub-cell philosophy (see Sect.~\ref{sec:AMR}).
\end{itemize}
We emphasize that the adaptivity of the main grid provided by the AMR approach has nothing to do with the subcell limiter. The AMR technology is used to refine and recoarsen the 
computational grid according to physical features that one wants to follow, while the subgrid limiter is only used to cope with shock waves or other discontinuities that require 
limiting of the DG scheme. In the following we provide the necessary minimum details for each of the above items, while addressing the reader 
to 
\cite{Dumbser2008,DumbserZanotti,HidalgoDumbser,GassnerDumbserMunz,Balsara2013934,AMR3DCL,Dumbser2014} for an exhaustive discussion of the subtleties that may be involved. 

\section{The ADER-DG scheme}
\label{sec:ADER-DG}
\subsection{The local space-time predictor}
\label{sec.galerkin}

At the heart of the ADER approach, either in the original version proposed in 
\cite{toro3} and \cite{titarevtoro} or in the later version proposed in \cite{DumbserEnauxToro,Dumbser2008,Balsara2013934}, that we also follow in this paper,  
there is the solution of the generalised or derivative Riemann problem. This requires a time evolution of known spatial derivatives of the polynomials 
approximating the solution at time $t^n$ and is in our case performed locally for each cell and independently from the neighbor cells. 
In the FV framework, such polynomials are obtained via reconstruction from the known cell averages of the conserved quantities. 
In the DG framework, on the contrary, no reconstruction is needed and the time evolution acts directly on  
the representation polynomials $\u_h(\x,t^n)$ of Eq.~(\ref{eqn.ansatz.uh}).
To show how the predictor works, 
we first transform the PDE system of  Eq.~\eqref{eqn.pde.nc} into a space-time reference coordinate system 
$(\xi,\eta,\zeta, \tau)$. In practice, the 
space-time control volume 
$\mathcal{C}_{ijkn}=[x_{i-\halb};x_{i+\halb}] \times [y_{j-\halb};y_{j+\halb}] \times [z_{k-\halb};z_{k+\halb}] \times [t^n;t^{n+1}]$ is mapped into the 
space-time reference element $[0;1]^{4}$ through the definitions
\begin{equation}
\label{eq:xi}
x = x_{i-\halb} + \xi   \Delta x_i, \quad 
y = y_{j-\halb} + \eta  \Delta y_j, \quad 
z = z_{k-\halb} + \zeta \Delta z_k, \quad
t = t^n + \tau \Delta t\,.
\end{equation} 
In general, we will use $T_E=[0;1]^d$ to denote the spatial reference 
elements in $d$ spatial dimensions.
As a result, Eq.~\eqref{eqn.pde.nc} will be rewritten as
\begin{equation}
\label{eqn.pde.nc.2d}
 \frac{\partial \u}{\partial \tau}
    + \nabla_\xi \cdot \F^* \left( \u \right) = 0,
\end{equation}
where
\begin{equation}
  \F^* := \Delta t \left( \partial \boldsymbol{\xi} / \partial \x  \right)^{T} \cdot \F(\u)
\end{equation}
with $(\boldsymbol{\xi}=\xi,\eta,\zeta)$ and
$\nabla_\xi = \partial \boldsymbol{\xi} / \partial \x \cdot \nabla$.
Multiplication of \eqref{eqn.pde.nc.2d} with a space-time test function 
$\theta_k=\theta_k(\boldsymbol{\xi},\tau)$ and integration over 
the space-time reference control volume $T_E \times [0;1]$ yields
\begin{equation}
\label{eqn.pde.nc.weak1}
 \int \limits_{0}^{1} \int \limits_{T_E} \theta_k \frac{\partial \u}{\partial \tau} \, d \boldsymbol{\xi} \, d \tau\,
    + \int \limits_{0}^{1} \int \limits_{T_E}\theta_k \nabla_\xi \cdot \F_h^* \left(\u\right)\, d \boldsymbol{\xi} \, d \tau\,
    = 0\,.
\end{equation}
In analogy to Eq.~(\ref{eqn.ansatz.uh}), we now introduce the 
\emph{discrete spacetime solution} of equation \eqref{eqn.pde.nc.weak1}, denoted by $\q_h$, as well as the corresponding one for the flux, i.e.
\begin{equation}
\label{eqn.st.state}
 \q_h = \q_h(\boldsymbol{\xi},\tau) =
 \theta_l \hat{\q}_l\,.
\end{equation}
\begin{equation}
\label{eqn.st.flux}
 \F^*_h = \F^*_h(\boldsymbol{\xi},\tau) =
 \theta_l \hat{\F}^*_l,
\end{equation}
The space-time test and basis functions $\theta_l$ are chosen again as tensor products of Lagrange interpolation polynomials passing through the Gauss-Legendre quadrature points. 
Due to this choice of a nodal basis, the degrees of freedom for the fluxes are simply the point--wise evaluation of the physical fluxes, namely
\begin{equation}
  \hat{\F}^*_l = \F^*(\hat{\q}_l)\,.
\end{equation}
The next crucial step of the approach consists of
integrating by parts in time the first term in \eqref{eqn.pde.nc.weak1}, while keeping the information local in space. This  yields
\begin{equation}
\label{eqn.pde.nc.dg1}
  \int \limits_{T_E}  
	\theta_k (\boldsymbol{\xi},1) \q_h \, d \boldsymbol{\xi}  - 
	\int \limits_{T_E}  
	\theta_k (\boldsymbol{\xi},0) \u_h \, d \boldsymbol{\xi} 
	   - \int \limits_{0}^{1} \int \limits_{T_E}
		\frac{\partial \theta_k}{\partial \tau}  \q_h  \, d \boldsymbol{\xi} \, d \tau
    + \int \limits_{0}^{1} \int \limits_{T_E} \theta_k \nabla_\xi \cdot \F^*_h \, d \boldsymbol{\xi} \, d \tau
    = 0.
\end{equation}
After substituting (\ref{eqn.st.state}) and (\ref{eqn.st.flux}) into Eq.~(\ref{eqn.pde.nc.dg1}) we obtain \citep{Dumbser2008,HidalgoDumbser,DumbserZanotti}
\begin{equation}
\label{eqn.pde.nc.dg2}
   \left( \int \limits_{T_E}  
	\theta_k (\boldsymbol{\xi},1) \theta_l (\boldsymbol{\xi},1)\, d \boldsymbol{\xi}
	- \int \limits_{0}^{1} \int \limits_{T_E}
		\frac{\partial \theta_k}{\partial \tau}   \theta_l  \, d \boldsymbol{\xi} \, d \tau
    \right)
    \hat{\q}_l 
    =   \left( \int \limits_{T_E}  
	\theta_k (\boldsymbol{\xi},0) \Phi_l \, d \boldsymbol{\xi}  \right)\hat \u_l^n
		 - \left(\int \limits_{0}^{1} \int \limits_{T_E} \theta_k \nabla_\xi \theta_l \, d \boldsymbol{\xi} \, d \tau\right) \F^*(\hat{\q}_l)\,.
\end{equation}
Equations \eqref{eqn.pde.nc.dg2} represents a nonlinear system to be solved in the unknown expansion 
coefficients $\hat{\q}_l$ of the local space-time predictor solution, while the terms $\hat \u_l^n$ are the known degrees of freedom of the DG polynomial at time level $t^n$. We note, incidentally, that, although the solution of the above nonlinear system is certainly demanding in terms of computational costs, it has the significant advantage that it can also cope with stiff source terms, which are absent for the system of equations considered in this paper but are quite common in several fields of applied mathematics and physical sciences \citep{DumbserZanotti,HidalgoDumbser,Zanotti2011,Dumbser-Uuriintsetseg2013}.
\subsection{The fully discrete one-step ADER-DG scheme}
\label{sec.ADERNC}

The spacetime solution $\q_h$ of Eq.~(\ref{eqn.pde.nc.dg2}) cannot of course provide the true solution at time level $t^{n+1}$, since it completely
neglects the contribution of fluxes from neighbouring cells. In our approach, the proper correction is obtained
through a fully discrete one-step ADER-DG scheme, which works as follows (see also \cite{dumbser_jsc,taube_jsc,QiuDumbserShu}). We first 
multiply the governing PDE \eqref{eqn.pde.nc} by a test 
function $\Phi_k$, identical to the spatial basis functions of Eq.~\eqref{eqn.ansatz.uh}. Second, 
we integrate over the space-time control volume $T_i \times [t^n;t^{n+1}]$. 
The flux divergence term is then integrated by parts in space, thus yielding 
\begin{equation}
\label{eqn.pde.nc.gw1}
\int \limits_{t^n}^{t^{n+1}} \int \limits_{T_i} \Phi_k \frac{\partial \u_h}{\partial t} d\x dt
+\int \limits_{t^n}^{t^{n+1}} \int \limits_{\partial T_i} \Phi_k \, \F\left(\u_h \right)\cdot\mathbf{n} \, dS dt 
-\int \limits_{t^n}^{t^{n+1}} \int \limits_{T_i} \nabla \Phi_k \cdot \F\left(\u_h \right) d\x dt 
= 0, 
\end{equation}
where $\mathbf{n}$ is the outward pointing unit normal vector on the surface $\partial T_i$ of element $T_i$. 
Since the discrete solution is allowed to be discontinuous at element boundaries, the surface integration involved in the second term of 
\eqref{eqn.pde.nc.gw1} is done through the solution of a Riemann problem, which is therefore deeply rooted in the DG scheme and guarantees 
the overall upwind character of the method \citep{cbs0,cbs1,cbs2,cbs3,cbs4}. 
Whatever numerical flux function (Riemann solver) is chosen, denoted as $\mathcal{G}$, the time integration of the second and of the third term of Eq.~\eqref{eqn.pde.nc.gw1} must be performed to the desired order of accuracy. 
To this extent, we use the local space-time predictor $\q_h$ from Sect.~\ref{sec.galerkin}, which clearly pays off at this stage, 
and allows to compute the numerical flux function of the second term as
$\mathcal{G}\left(\q_h^-, \q_h^+ \right)$ 
and the physical flux of the third term as $\F\left(\q_h \right)$.
We emphasize that  $\q_h^-$ and $\q_h^+$ are the left and right states of the Riemann problem.
On the other hand, after inserting $\u_h$, as given by \eqref{eqn.ansatz.uh}, in the first term of \eqref{eqn.pde.nc.gw1}
we find the following arbitrary high order accurate one-step discontinuous Galerkin (ADER-DG) scheme: 
\begin{equation}
\label{eqn.pde.nc.gw2}
\left( \int \limits_{T_i} \Phi_k \Phi_l d\x \right) \left( \hat{\u}_l^{n+1} -  \hat{\u}_l^{n} \right) +
\int\limits_{t^n}^{t^{n+1}} \int_{\partial T_i} \Phi_k \, \mathcal{G}\left(\q_h^-, \q_h^+ \right)\cdot\mathbf{n} \, dS dt 
-\int\limits_{t^n}^{t^{n+1}} \int_{T_i} \nabla \Phi_k \cdot \F\left(\q_h \right) d\x dt  = 0\,.
\end{equation}
Concerning the choice of the Riemann solver, in this paper we have used both
the simple Rusanov flux\footnote{This is sometimes referred to as the local Lax Friedrichs flux.} \citep{Rusanov:1961a}, and the more sophisticated
Osher-type flux proposed in \cite{OsherUniversal}, which requires the knowledge of the eigenvectors of the system 
\eqref{eqn.pde.nc}. For a very recent alternative family of genuinely multi-dimensional HLL Riemann solvers, see \cite{balsarahlle2d,balsarahllc2d,BalsaraMUSIC1,BalsaraMultiDRS}. 

The effective order of accuracy of the ADER-DG scheme resulting from \eqref{eqn.pde.nc.gw2} is
$N+1$, both in space and in time, as long as the solution remains smooth. 
In spite of its great ability in achieving sub-cell resolution even on very coarse grids, the ADER-DG scheme, as well any other unlimited DG scheme, will
fail at discontinuities due to the Gibbs phenomenon. For this reason it is necessary to introduce some sort of limiter, 
which should ideally preserve the typical sub-cell resolution properties of the DG method. Such an approach, proposed and discussed with all details 
in \cite{Dumbser2014}, is briefly summarized in the next section. Before proceeding, we also comment on the Courant-Friedrichs-Lewy (CFL) restriction 
imposed by explicit DG schemes. In multiple space dimensions, the time step is usually restricted as 
\citep{krivodonova2013analysis}
\begin{equation}
\Delta t < \frac{1}{d} \frac{1}{(2N+1)} \frac{h}{|\lambda_{\max}|}\,,
\label{eqn.cfl} 
\end{equation}
where $d$ is the number of space dimensions, $h$ and $|\lambda_{\max}|$ are a characteristic mesh size and the maximum signal velocity, respectively.  
The factor $2N+1$ in the denominator of \eqref{eqn.cfl} will motivate the choice for the number of sub-cells required by the limiter, as we explain below. 
Note that for ADER-DG schemes and Lax-Wendroff DG schemes, the CFL condition is even slightly more severe, see \cite{Dumbser2008,QiuDumbserShu}. 

\section{A posteriori sub-cell limiter}
\label{sec:A-posteriori-sub-cell-limiter}
 
Let us assume that we have obtained the discrete representation $\textbf{u}_h(\textbf{x},t^n)$,
within a general cell $T_i$, of the solution. In order to update the solution to the next time level, 
we first calculate a so-called \textit{candidate} solution, denoted so forth as $\textbf{u}_h^*(\textbf{x},t^{n+1})$, which 
results from the \textit{unlimited} scheme of Eq.~\eqref{eqn.pde.nc.gw2}. Due to the appearance of possible oscillations, the candidate solution
may not be acceptable everywhere in the computational domain, and a number of 
detection criteria must be fulfilled in order to promote the candidate solution to the \textit{accepted} discrete solution at the new time level. 
The first criterion is based on physical considerations and it consists of checking whether $\textbf{u}_h^*(\textbf{x},t^{n+1})$ verifies the physical 
positivity constraints.  
This is a necessary condition for a number of variables appearing in many conservation laws 
(mass, density, pressure, internal energy, water depth, etc.). For simplicity, and because this will be the case in the rest of the paper, 
we can refer to the Euler equations for gas dynamics, for which density and pressure must remain positive. 
Since the DMP was already a very useful tool to construct high resolution shock capturing finite volume schemes in the past, our second detection criterion 
is a \emph{relaxed discrete maximum principle} (DMP) in the sense of polynomials. To this purpose, the following condition must be verified  
\begin{equation} 
\label{NAD}
\min \limits_{y\in {\cal{V}}_i} (\textbf{u}_h(\textbf{y},t^n))-\delta \leq \textbf{u}_h^*(\textbf{x},t^{n+1})\leq \max \limits_{\textbf{y}\in {\cal{V}}_i}(\textbf{u}_h(y,t^n))+\delta, 
\qquad \forall \textbf{x}  \in T_i\,,
\end{equation}
where the set ${\cal{V}}_i$ contains the cell $T_i$ and the Voronoi neighbor cells which share a common node with $T_i$. 
In practice, Eq.~\eqref{NAD} says that the polynomial representing the candidate solution must lie between the minimum and the maximum 
of the polynomials representing the solution at the previous time step in the set ${\cal{V}}_i$.
The small quantity $\delta$ in (\ref{NAD}) 
is used to relax the maximum principle thus allowing for small undershoots and overshoots and it can avoid problems with roundoff errors. 
The value used here, as recommended in \cite{Dumbser2014}, is
\begin{equation}
	\delta =\max \left( \delta_0, \epsilon \cdot \left( \max \limits_{y\in {\cal{V}}_i}(\textbf{u}_h(y,t^n))- \min \limits_{y\in {\cal{V}}_i}(\textbf{u}_h(\textbf{y},t^n))\right)\, \right),
\end{equation}
with $\delta_0=10^{-4}$ and $\epsilon=10^{-3}$. It is interesting to remark that the physical and the numerical criteria are totally independent, which implies that the relaxation of the maximum principle does not affect the positivity of the solution. Moreover, this approach takes into account the information from two different time levels, $t^n$ and $t^{n+1}$, whereas classical indicators typically use information from one time level only.

During the detection phase, if either the first physical criterion or the second numerical criterion is not fulfilled, the corresponding cell will acquire
a so-called \emph{limiter status} $\beta=1$, which flags the cell as \emph{troubled}. On the other hand, if both criteria are met, the limiter status is set to $\beta=0$. Troubled cells
immediately generate a sub-grid, for which an \textit{alternative data representation} $\textbf{v}_h(\textbf{x},t^n)$ must be provided. 
This new solution is expressed by a set of piecewise constant sub-cell averages $\textbf{v}_{i,j}^n$. These values are computed via $L_2$ projection on the $(N_s)^d$ sub-cells $S_{i,j},j=1,\cdots,(N_s)^d$ in which $T_i$ is divided, where $N_s= 2N+1$, i.e.
\begin{equation} \label{vh}
	v_{i,j}^n=\frac{1}{|S_{i,j}|}\int_{S_{i,j}}{\textbf{u}_h(\textbf{x},t^n)\,d\textbf{x}}=
	\frac{1}{|S_{i,j}|}\int_{S_{i,j}}{\hat{\textbf{u}}_l^n\phi_l(\textbf{x})\,d\textbf{x}}, \qquad \forall S_{i,j}\in {\cal{S}}_i\,,
\end{equation}
and where ${\cal{S}}_i=\bigcup_j S_{i,j}$ is the set of the sub-grid cells. 
We emphasize that the choice $N_s= 2N+1$ is not a heuristic one, but it is properly motivated by an optimality argument. 
With $N_s= 2N+1$ the maximum timestep of the ADER-DG scheme on the main grid (c.f. Eq.~\eqref{eqn.cfl}) matches the maximum 
possible time step of the ADER finite volume scheme on the sub-grid. This leads to the maximum admissible $\textnormal{CFL}$ 
number for the sub-grid finite volume scheme, thus minimizing its dispersion and dissipation error. 
Note that for ADER finite volume methods applied to the linear advection equation in 1D the error terms scale with  
$(1-\textnormal{CFL})$, see \cite{dumbser_diffapprox}. For an alternative subcell finite volume limiter approach that works 
with \apriori indicator functions and a classical TVD scheme on $N_s=N+1$ subgrid cells, see \cite{Sonntag,Fechter1}.

Using the data representation $\textbf{v}_h(\textbf{x},t^n)$ as initial condition, the next step consists of updating  the discrete solution by means of a robust scheme on the sub-grid. While any TVD scheme could serve to the scope, we have preferred to adopt a third order ADER-WENO finite volume scheme to avoid the clipping of smooth extrema. As a result, both the ADER-DG scheme on the main grid and the ADER-WENO finite volume scheme on the sub-grid are \emph{one-step} schemes. This approach has the net effect of reducing the total amount of MPI communications with respect to traditional Runge-Kutta schemes. 
Once the solution for all the troubled cells has been recomputed on the sub-grid, the solution on the main grid is recovered through the requirement that
\begin{equation}
\label{eqn.sub-cell.r} 
  \int \limits_{S_{i,j}} \u_h(\x,t^{n+1}) d\x = \int \limits_{S_{i,j}} \v_h(\x,t^{n+1}) d\x, \qquad \forall S_{i,j} \in \mathcal{S}_i\,.
\end{equation} 
which is a standard reconstruction problem arising both 
within the finite volume context as well as for spectral finite volume methods \citep{spectralfv2d,spectralfv3d,spectralfv.bnd}. 
It may well be the case that a cell is marked as troubled for a sequence of successive time steps. Under these circumstances, the initial data for 
$\v_h(\x,t^n)$ are directly available on the sub-grid from the 
ADER-WENO finite volume scheme of the previous time step.

\section{AMR with the sub-cell limiter}
\label{sec:AMR}

\subsection{Summary of the cell-by-cell AMR implementation}

There are basically two major strategies for implementing an AMR algorithm.
The first strategy employs a nested structure of independent overlaying sub-grid-patches \citep{Berger-Oliger1984,
berger85,Berger-Colella1989}. The second strategy, on the contrary, refines each cell individually and it is referred to as
a "cell-by-cell" refinement. Due to its simple tree-type data structure \citep{Khokhlov98,AMR3DCL}, and also for its  
slightly more general formulation, in our work we have adopted the latter approach. 

By defining the maximum level of refinement $\ell_\text{max}$, each level of refinement is indicated by $\ell$, such 
that $0\le\ell\le\ell_\text{max}$. This means that we are considering up to $\ell_\text{max}$ overlaying uniform lattices, whose cells 
are activated only where and when necessary. The union of all the cells up to level $\ell$ is denoted by $\Ell_\ell$. Every cell
is labeled by a positive integer number $\emm$ and can be denoted as $\mathcal{C}_\emm$, with $\emm\le N_\text{e}$, where 
$N_\text{e}$ is the (time-dependent) total number of the cells, or elements. In $d$ space dimensions, every cell $\mathcal{C}_\emm$ 
has up to $2d$ \emph{Neumann neighbors},  namely neighbor cells sharing a face with 
$\mathcal{C}_\emm$. Furthermore, one can identify up to $3d-1$ \emph{Voronoi neighbors}, namely neighbor cells that share
at least one lattice-node with $\mathcal{C}_\emm$.

Moreover, each cell of a given $\ell$-th level has one among three possible \emph{status} $\sigma$: \emph{active} cells ($\sigma=0$) are updated according
to the finite element ADER-DG scheme described in Section \ref{sec.ADERNC}; \emph{virtual children} cells ($\sigma=1$), or \emph{virtual children}, 
are updated according to 
standard $L_2$ projection of the high order polynomial of the so-called \emph{mother} cell at  $(\ell-1)$-th level; \emph{virtual mother} cells ($\sigma=-1$), or \emph{virtual mothers}, are updated by averaging recursively the \emph{children cells} of the upper refinement levels, from the $(\ell+1)$-th to the level of the corresponding \emph{active} children cells.
More specifically, whenever $\mathcal{C}_\emm$ is refined, it generates $\mathfrak{r}^d$ children cells, such that
\begin{equation}
	\Delta x_\ell = \mathfrak{r}\:\Delta x_{\ell+1}; \;\;\;\Delta y_\ell = \mathfrak{r}\: \Delta y_{\ell+1}; \;\;\;  \Delta z_\ell = \mathfrak{r}\:\Delta z_{\ell+1}\,,
\end{equation}
where $\mathfrak{r}$ is the \emph{refinement factor}.
We emphasize that the time steps can be chosen \emph{locally} \citep{LTS,toro3,stedg1},  depending on the refinement level, such that 
\begin{equation}
\Delta t_\ell = \mathfrak{r}\: \Delta t_{\ell+1},
\end{equation}
with noticeable increase in performance. 
Consistently with the chosen nomenclature, the tree-structure formed by the union of a fixed mother cell of the coarsest refinement level ($\ell =0$), with all the recursive (contained) children can be referred to as a \emph{family-tree}. 
All over the computation, we need a \emph{criterion} to mark any given active cell $\mathcal{C}_\emm$ as a cell requiring refinement or recoarsening.
We therefore introduce a \emph{refinement-estimator function} $\chi_\emm$, built according to \cite{Loehner1987},
which  involves up to the second order derivative of an \emph{indicator function} $\Phi$, i.e.
\begin{equation} \label{eq:chi_m}
	\chi_\emm (\Phi) =  \sqrt{ \frac{\sum_{k,l}{\left(\left. \partial ^2 \Phi \middle/ \partial x_k \partial x_l \right. \right)^2}}{ \sum_{k,l}{ \left[ \left. \Big( \left| \left. \partial \Phi \middle/ \partial x_k\right.\right|_{i+1} + \left| \left.\partial \Phi \middle/ \partial x_k \right. \right|_i \Big) \middle/ \Delta x_l \right. + \epsilon \left|\frac{\partial^2 }{\partial x_k\partial x_l} \right|\left| \Phi \right|    \right]^2}}  }\,.
\end{equation}
Whenever $\chi_\emm > \chi_\text{ref}$, $\mathcal{C}_\emm$ is marked for refinement, while  it is marked for recoarsening if $\chi_\emm < \chi_\text{rec}$.
The sum $\sum_{k,l}$ is intended to be the double summation over all the spatial indexes, so that cross derivatives contributions are properly taken into account. 
$\Phi=\Phi(\u)$ is a generic function of the conservative variables $\u$, and in all the numerical tests reported in Sect.~\ref{sec:numerics} for the Euler equations we have used $\Phi(\u)=\rho$. The two parameters
$\chi_\text{ref}$ and $\chi_\text{rec}$ are moderately model-dependent and they are typically chosen in the range $\sim[0.2,0.25]$ and $\sim[0.05,0.15]$, respectively. Finally, $\epsilon=0.01$ is a filter-parameter that avoids unnecessary mesh-refinement in regions affected by ripples.
For problems involving the propagation of discontinuities, most notably
shock waves, it is advisable to anticipate the arrival of the discontinuity in such a way that it is always surrounded by a few additional refined cells, thus avoiding any pre- or post-shock oscillation. In practice, this is obtained by forcing the refinement of a suitable number of cells, usually one or two, in a neighborhood of the cell which has been marked for refinement according to the standard criterion.
 
In the numerical implementation of our AMR algorithm we have followed a number of basic rules: 
\begin{itemize}
\item along a family-tree an \emph{active cell} can only have recursive \emph{non-active mothers} and \emph{non-active children};
\item only active cells can be refined; 
\item two Voronoi's neighbors belong either to the same, or to an \emph{adjacent} refinement level, which implies that they must have $\Delta\ell \leq 1$.
\end{itemize}
We emphasize that, according to these conventions, the real (active) grid is a non-overlapping, non-conforming grid.  
Full details about the implementation and the parallelization of the AMR framework through the standard Message Passing Interface (MPI) can be found
 in \cite{AMR3DCL,AMR3DNC}. For details on the high order local time stepping (LTS) procedure see also \cite{AMR3DCL}.

\subsection{Incorporation of the sub-cell limiter into the AMR framework}

What we discussed in Sect.~\ref{sec:A-posteriori-sub-cell-limiter}, namely the \aposteriori sub-cell limiter which is activated in the troubled zones
of the ADER-DG scheme, must be properly nested within the AMR framework. 
In order to understand how the interaction works, let us first list the basic rules that we have followed
\begin{itemize}
\item The virtual children cells inherit the limiter status of their active mother cell.
\item If at least one active child is flagged as troubled, then the (virtual) mother is also flagged as troubled.
\item Cells which have been flagged as troubled cannot be recoarsened.
\end{itemize}
Because of the presence of the limiter, the two typical AMR operations represented by \emph{projection} and \emph{averaging} 
must be also extended to the alternative data representation $\textbf{v}_h(\textbf{x},t^n)$.  
Let us denote the sub-grid of a generic cell $\mathcal{C}_n$ at level $\ell$ as $\mathcal{S}_n^{\ell}$ and the data representation 
$\textbf{v}_h(\textbf{x},t^n)$ at level $\ell$ simply as $\textbf{v}_h(\mathcal{S}_n^{\ell})$. Let us further denote a generic 
virtual child cell as $\mathcal{C}_v$ and the virtual mother or parent cell as $\mathcal{C}_p$. Then, in general, we need to be able to perform
the two operations 
\begin{eqnarray}
\label{AMR-limiter-projection}
&&\textbf{v}_h(\mathcal{S}_n^{\ell}) \rightarrow \textbf{v}_h(\mathcal{S}_v^{\ell+1}): \ \ \ {\rm DG \ limiter - AMR \ projection },\\
\label{AMR-limiter-averaging}
&&\textbf{v}_h(\mathcal{S}_n^{\ell}) \rightarrow \textbf{v}_h(\mathcal{S}_p^{\ell-1}): \ \ \ {\rm DG \ limiter - AMR \ averaging },
\end{eqnarray}
which we describe below.
%
\subsubsection{DG limiter - AMR projection}

\begin{figure}
	\centering
		\includegraphics[width=0.95\textwidth]{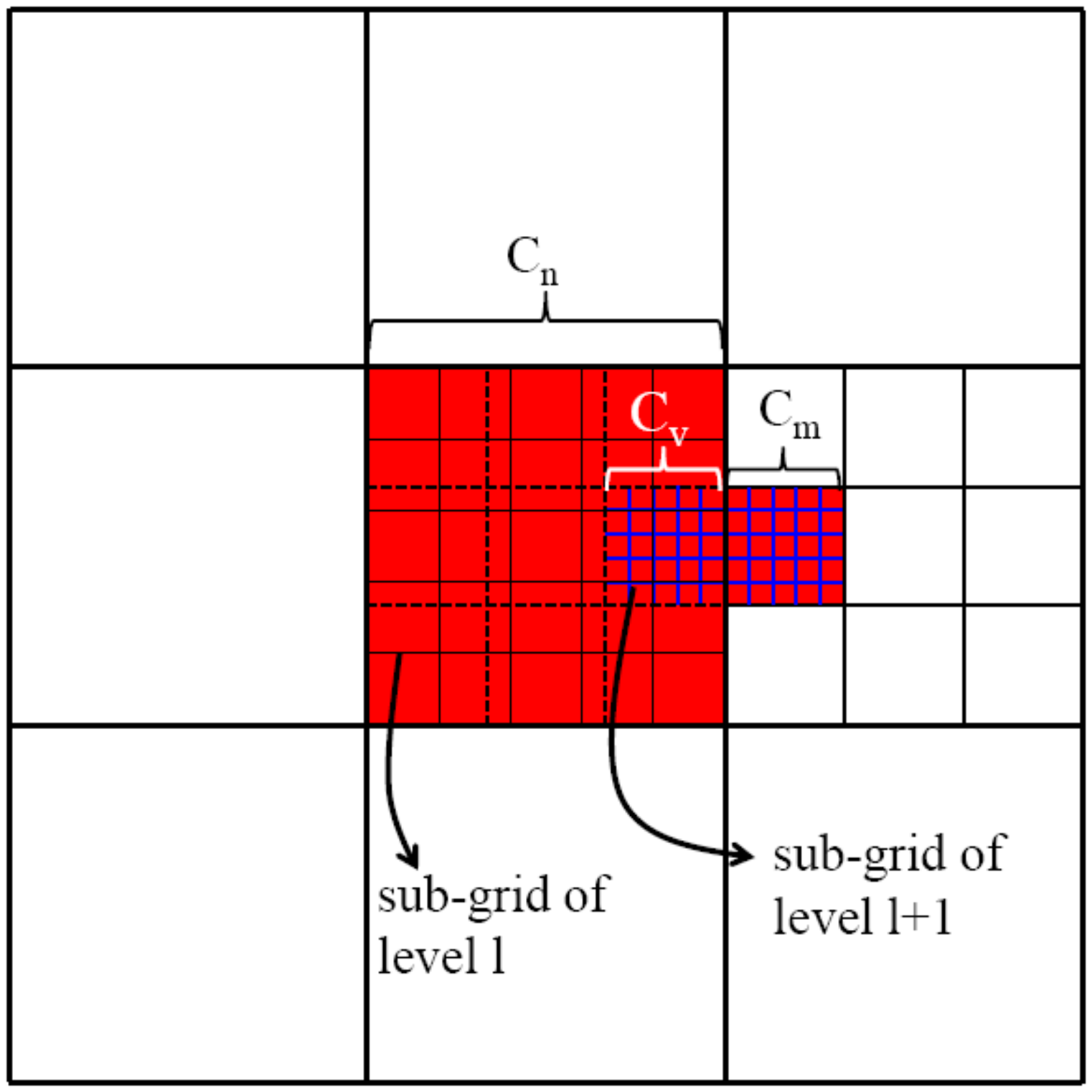}
	\caption{Sketch of the combination of AMR and DG sub-cell reconstruction. The cell $\mathcal{C}_n$ at level $\ell$ and the cell  $\mathcal{C}_m$ at level $\ell+1$ (both in red) have
limiter status $\beta=1$. The cell $\mathcal{C}_n$ must project $\textbf{v}_h$ from the sub-grid of level $\ell$ to the sub-grid of level $\ell+1$ in the virtual cell $\mathcal{C}_v$.
(see color version on-line).}
	\label{fig:AMR_DGsub-cell}
\end{figure}

This operation becomes necessary when an active cell with limiter status $\beta=1$,
namely a troubled cell, has virtual children cells. 
In such circumstances, we need to project the alternative data representation
$\textbf{v}_h(\textbf{x},t^n)$ from the sub-cells of a given level of refinement $\ell$ to the sub-cells of the next level 
$\ell+1$. 
We recall that a pure DG scheme with AMR, but without limiters, would not require any virtual cell (status $\sigma=\pm 1$), because pure DG schemes
do not perform any reconstruction. We also recall that in our implementation virtual children cells are created to allow any cell marked for refinement to perform 
a spatial reconstruction, and more precisely when the stencil corresponding to the specific reconstruction procedure chosen (TVD, WENO, etc.) covers  adjacent cells belonging to different levels of refinement. However, our DG scheme is not pure, because it works in combination with the limiter, and the limiter involves a WENO reconstruction on the sub-grid. Hence, our ADER-DG-AMR scheme still implies the introduction of virtual cells, which must be created when the WENO reconstruction
on the sub-grid of level $\ell+1$ uses a stencil that covers a portion of the grid belonging to level $\ell$. In such circumstances, it is necessary to 
perform the operation expressed by (\ref{AMR-limiter-projection}) above.
A simplified situation is reported in Figure~\ref{fig:AMR_DGsub-cell}, sketching a two-dimensional configuration in which AMR and the sub-cell limiter of the DG scheme 
are interlinked.  
In that figure two AMR refinement levels are involved.  The cell $\mathcal{C}_n$ at level $\ell$ and the cell  $\mathcal{C}_m$ at level $\ell+1$ have
limiter status $\beta=1$, and for this reason they are colored in red. In order to allow $\mathcal{C}_m$ to perform the WENO reconstruction on its sub-grid,
cell $\mathcal{C}_n$ must project $\textbf{v}_h$ from the sub-grid of level $\ell$ to the sub-grid of level $\ell+1$ in the virtual cell $\mathcal{C}_v$.
Hence, the subcell averages on the finer level $\ell+1$ are computed from the condition that 
\begin{equation}
\label{eqn.subcell.r} 
  \int \limits_{S_{v,j}} {\bf v}_h (\mathcal{S}_v^{\ell+1}) d\x = \int \limits_{S_{i,j}} \mathcal{W} \left( {\bf v}_h(\mathcal{S}_n^{\ell}) \right) d\x, 
	\qquad	\forall S_{v,j} \in \mathcal{S}_v^{\ell+1}, 
\end{equation} 
where $\mathcal{W}$ denotes the WENO reconstruction operator applied on the cell averages of the subgrid on level $\ell$. We use a WENO reconstruction to pass 
subgrid data from the coarse level to the finer one, since this projection operation is carried out in troubled cells where typically discontinuities are present.  
Therefore, we need a nonlinear, essentially non-oscillatory reconstruction that is at the same time high order accurate and which is also able to deal with 
shocks and other discontinuities. 

%
\subsubsection{DG limiter - AMR averaging}

Conversely, we also need to 
perform the  averaging of $\textbf{v}_h(\textbf{x},t^n)$ from the sub-cells of a given level $\ell$ to the sub-cells of the previous level $\ell-1$.
Then the averaging operator acting on the degrees of freedom of the sub-grid WENO polynomial can be written in a compact form as 
\begin{equation}
  \int \limits_{S_{p,j}} {\bf v}_h (\mathcal{S}_p^{\ell-1}) d\x = \int \limits_{S_{p,j}}  {\bf v}_h(\mathcal{S}_n^{\ell} ) d\x, 
	\qquad	\forall S_{p,j} \in \mathcal{S}_p^{\ell-1}, 
\label{eqn.average} 
\end{equation}  
%
From an operational point of view, this transformation is most conveniently performed in a dimension-by-dimension fashion. No reconstruction is needed here, 
since the averaging over known cell averages is trivial. 
%

\section{Numerical results} 
\label{sec:numerics}

\subsection{Euler equations of compressible gas dynamics}
The first set of PDEs that we have considered in our numerical tests is represented by the classical Euler equations, which can be written as a system of conservation laws as required by Eq.~\eqref{eqn.pde.nc}, where the conserved variables and the corresponding
fluxes are given by
\begin{equation}
\u=\left[\begin{array}{c}
\rho \\ \rho \mathbf{v} \\ E
\end{array}\right],~~~
{\F}=\left[\begin{array}{c}
 \rho \mathbf{v}\\
 \rho \mathbf{v} \mathbf{v} + p \mathbf{I} \\
 (E+p) \mathbf{v} 
\end{array}\right].
\label{eq:fluxes}
\end{equation}
Here $\mathbf{v}=(u,v,w)$ are the velocity components, $p$ is the pressure, $\rho$ is the mass density,
$E=p/(\gamma-1)+ \rho \mathbf{v}^2/2$ is the total energy density including the thermal and the kinetic 
contributions, $\mathbf{I}$ is the identity matrix, while $\gamma$ is the adiabatic index of the ideal gas, 
which follows the standard equation of state 
\begin{equation}
\label{eq:EOS}
p=\rho\epsilon(\gamma-1)\,,
\end{equation}
where $\epsilon$ is the internal energy per unit mass. 
The Jacobian matrix associated to the Euler equations has eigenvalues that are all real and 
a set of linearly independent eigenvectors \citep{toro-book}, thus allowing for the implementation of a large class of Riemann solvers.
In the six subsections below we discuss a sample of classical test cases that involve the
propagation of linear and non-linear waves admitted by the Euler equations.
For practical purposes, we have represented 
in blue the unlimited cells, which, in the last time step,  have been successfully evolved through the standard ADER-DG-AMR scheme. 
Conversely, we have represented in red the troubled cells, with limiter status $\beta=1$, which required the activation of the subcell limiter.

\subsubsection{Numerical convergence study} 
In order to asses the convergence properties of the ADER-DG-AMR scheme we have considered the solution of the two-dimensional isentropic vortex, which 
admits an analytic solution~\citep{Shu1}. The test consists of the advection  of a vortex with initial conditions given by
a perturbation superposed  to a uniform mean flow as
\begin{equation}
\left( \rho,u,v,w,p \right) =(1+\delta\rho, 1+\delta v_x, 1+\delta v_y, 0, 1+\delta p)\,,
\end{equation} 
with 
\begin{equation}
\left(\begin{array}{c}
\delta \rho \\ \delta v_x \\ \delta v_y \\ \delta p 
\end{array}\right)
=
\left(\begin{array}{c}
(1+\delta T)^{1/(\gamma-1)}-1 \\
-(y-5)\epsilon/2\pi \exp{[0.5(1-r^2)]} \\
\phantom{-}(x-5)\epsilon/2\pi \exp{[0.5(1-r^2)]} \\
(1+\delta T)^{\gamma/(\gamma-1)}-1
\end{array}\right).~~~
\label{eq:pert}
\end{equation}
The perturbation in the temperature is 
\begin{equation}
\delta T=-\frac{\epsilon^2(\gamma-1)}{8\gamma\pi^2}~\exp{(1-r^2)}\,,
\end{equation}
where $r^2=(x-5)^2+(y-5)^2$, while the vortex strength is $\epsilon=5$ and the adiabatic index is $\gamma=1.4$. 
It is easy to check that, under these conditions, the entropy per unit mass $s=p/\rho^\gamma$ is constant everywhere.
The numerical domain is the square $\Omega=[0,10]\times[0,10]$, 
and periodic boundary conditions are used along the four edges. In this way, after setting the final time of the simulation to $t_{\rm final}=10$, the vortex recovers the initial position.
We have solved this problem using the Rusanov flux with reconstruction in characteristic variables. 
Due to the smoothness of the solution, we expect that the sub-cell limiter is never activated, which  is indeed the case. 
We have performed a convergence study by varying $N$ from 2 to 8, with $\ell_{\rm max}=1$ and a refinement factor $\mathfrak{r}=3$, except for the case $N=8$, for which we have used   
$\mathfrak{r}=2$. 
A regular refinement over the moving vortex is better obtained by applying a refinement criterion based on the cell average of the mass density, rather than by applying the standard procedure based on Eq.~\eqref{eq:chi_m}.
In practice, and just for this test,  a cell is marked for refinement if the cell average of the variable $\rho$ is smaller than the threshold $\bar\rho=0.75$.  
Table \ref{tab:Vortex_Error} summarizes the results of this analysis by reporting the $L^1, L^2$ and $L^\infty$ norms of the error, computed with respect to the available analytic solution at time $t=t_{\rm final}$. The second column of the table reports the number of cells, along each direction, of the initial grid at the level zero. When $N\geq 6$, very coarse initial meshes have been adopted, since for larger values of $N_x$ the round-off errors affect negatively the outcome of the test. With this caveat in mind, the computed orders of convergence are in very good agreement with the nominal ones up to $N=8$, thus confirming the high order of accuracy of the proposed ADER-DG scheme even in combination with AMR and time-accurate local time stepping. 
 \begin{table}[!htbp] 
 \centering
 \numerikNine
 \begin{tabular}{|c|c||ccc|ccc|c|}
   \hline
   \multicolumn{9}{|c|}{\textbf{2D isentropic vortex problem --- ADER-DG-$\mathbb{P}_N$ + WENO3 SCL}} \\
   \hline
    & $N_x$ & $L^1$ error & $L^2$ error & $L^\infty$ error & $L^1$ order & $L^2$ order & $L^\infty$ order &
   Theor. \\
   \hline
   \hline
   \multirow{4}{*}{\rotatebox{90}{{DG-$\mathbb{P}_2$}}}
& 15	& 5.5416E-2	& 1.1075E-2	& 1.2671E-2	&---	&---	&--- & \multirow{4}{*}{3}\\
& 30	& 5.7101E-3	& 1.0984E-3	& 1.7374E-3	& 3.28	& 3.33	& 2.87 &\\
& 60	& 8.8511E-4	& 1.8805E-4	& 3.4727E-4	& 2.69	& 2.55	& 2.32 &\\
& 90	& 3.0025E-4	& 6.6257E-5	& 1.3176E-4	& 2.67	& 2.57	& 2.39 &\\
   \cline{2-8}
   \hline
	 \multirow{4}{*}{\rotatebox{90}{{DG-$\mathbb{P}_3$}}}
& 15	& 6.4357E-3	& 1.0325E-3	& 1.0026E-3	&---	&---	&--- & \multirow{4}{*}{4}\\
& 30	& 2.9981E-4	& 4.4304E-5	& 4.2822E-5	& 4.42	& 4.54	& 4.55 &\\
& 60	& 1.1141E-5 & 1.6679E-6	& 2.2108E-6	& 4.75	& 4.73	& 4.27 &\\
& 90	& 1.6787E-6 & 2.9117E-7	& 5.0366E-7	& 4.67	& 4.30	& 3.65 &\\
   \cline{2-8}
   \hline
	 \multirow{4}{*}{\rotatebox{90}{{DG-$\mathbb{P}_4$}}}
& 10	& 5.0587E-3	& 8.2103E-4	& 1.0921E-3	&---	&---	&--- & \multirow{4}{*}{5}\\
& 15	& 6.3888E-4	& 1.0137E-4	& 1.2972E-4	& 5.10	& 5.16	& 5.25 &\\
& 20	& 1.5369E-4 & 2.3219E-5	& 3.5064E-5	& 4.95	& 5.12	& 4.55 &\\
& 25	& 5.1581E-5 & 7.8567E-6	& 1.2824E-5	& 4.89	& 4.86	& 4.51 &\\
   \cline{2-8}
   \hline	
	 \multirow{4}{*}{\rotatebox{90}{{DG-$\mathbb{P}_5$}}}
& 15	& 1.1135E-4	& 1.6708E-5	& 2.5184E-5	&---	&---	&--- & \multirow{4}{*}{6}\\
& 20	& 1.8700E-5	& 2.7597E-6	& 3.4678E-6	& 6.20	& 6.26	& 6.89 &\\
& 25	& 3.9941E-6	& 6.0874E-7	& 9.4323E-7	& 6.92	& 6.77	& 5.83 &\\
& 30	& 1.4623E-6	& 2.1969E-7	& 3.0234E-7	& 5.51	& 5.59	& 6.24 &\\
   \cline{2-8}
   \hline		
	 \multirow{4}{*}{\rotatebox{90}{{DG-$\mathbb{P}_6$}}}
& 5	  & 1.5485E-2	& 2.5835E-3	& 2.6686E-3	&---	&---	&--- & \multirow{4}{*}{7}\\
& 10	& 1.8390E-4	& 2.9877E-5	& 4.1129E-5	& 6.40	& 6.43	& 6.02 &\\
& 15	& 9.8578E-6 & 1.6642E-6	& 2.9090E-6	& 7.22	& 7.12	& 6.53 &\\
& 20	& 1.2041E-6 & 2.0205E-7	& 3.6192E-7	& 7.31	& 7.33	& 7.24 &\\
   \cline{2-8}
   \hline	
	 \multirow{4}{*}{\rotatebox{90}{{DG-$\mathbb{P}_7$}}}
& 5	  & 6.2402E-3	& 1.0963E-3	& 1.4947E-3	&---	&---	&--- & \multirow{4}{*}{8}\\
& 9 	& 6.0168E-5	& 1.0210E-5	& 1.2830E-5	& 7.90	& 7.96	& 8.09 &\\
& 11	& 1.5676E-5 & 2.4524E-6	& 4.0665E-6	& 6.70	& 7.11	& 5.73 &\\
& 13	& 4.8297E-6 & 7.7831E-7	& 1.0593E-6	& 7.05	& 6.87	& 8.05 &\\
   \cline{2-8}
   \hline
	 \multirow{4}{*}{\rotatebox{90}{{DG-$\mathbb{P}_8$}}}
& 7	 & 1.3473E-4 & 2.1259E-5	& 2.3665E-5	&---	&---	&--- & \multirow{4}{*}{9}\\
& 9	 & 1.8066E-5 & 2.8661E-6	& 3.6534E-6	& 7.99	& 7.97	& 7.43 &\\
& 11 & 2.7718E-6 & 4.2166E-7	& 5.2952E-7	& 9.34	& 9.55	& 9.62 &\\
& 13 & 6.2220E-7 & 1.0475E-7	& 1.4401E-7	& 8.94	& 8.34	& 7.79 &\\
   \cline{2-8}
   \hline
 \end{tabular}
 \caption{ \label{tab:Vortex_Error} $L^1, L^2$ and $L^\infty$ errors and convergence rates for the 
   2D isentropic vortex problem for the ADER-DG-$\mathbb{P}_N$ scheme with sub-cell limiter and adaptive mesh refinement. One level of refinement has been used with a refinement factor $\mathfrak{r}=3$,
	except for the case $N=8$, for which we have used $\mathfrak{r}=2$.}
 \end{table}
%

\subsubsection{Riemann problems} 
\begin{figure}
  \begin{center} 
  \begin{tabular}{c} 
    \includegraphics[width=0.85\textwidth]{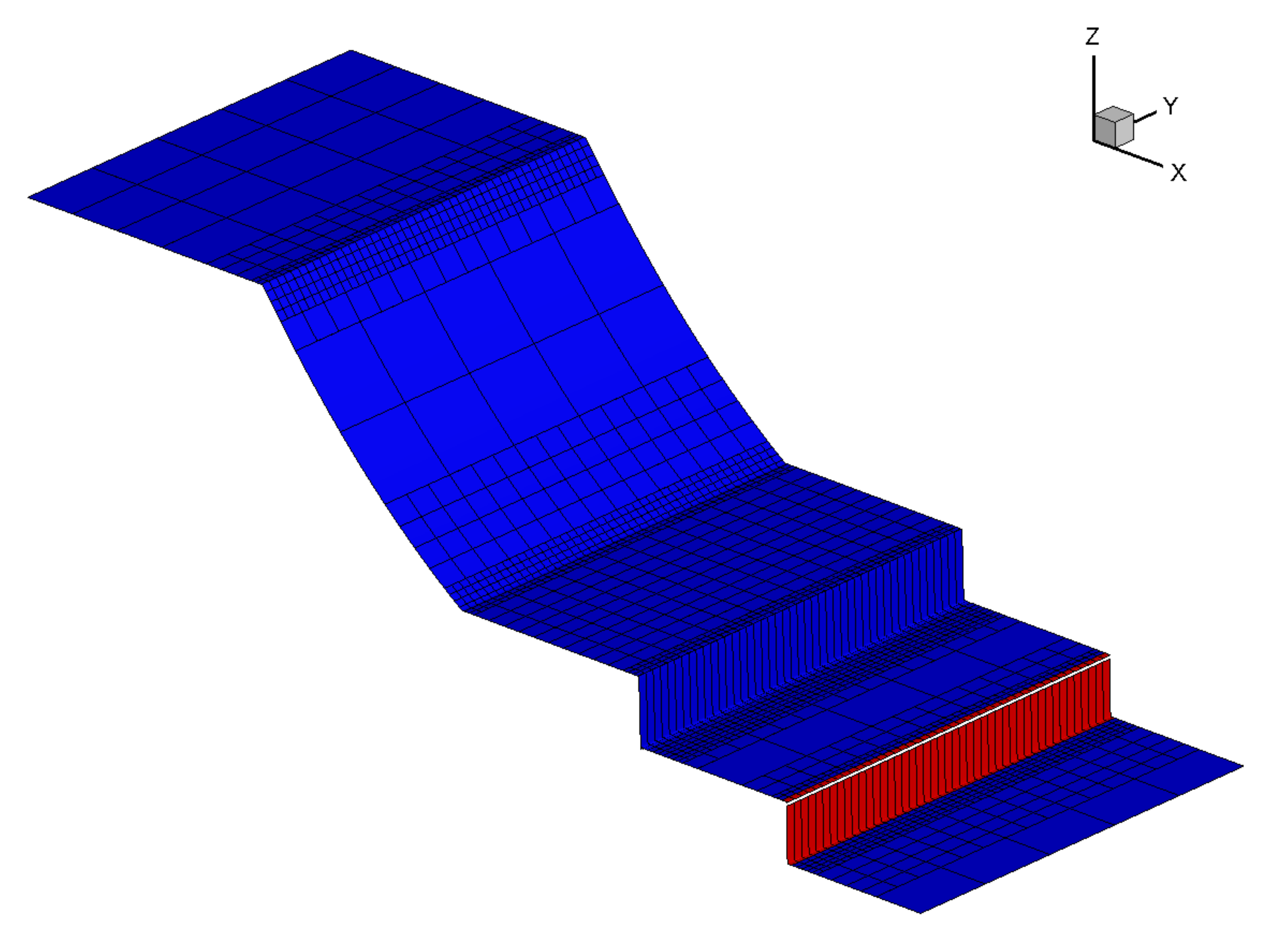} \\
    \includegraphics[width=0.85\textwidth]{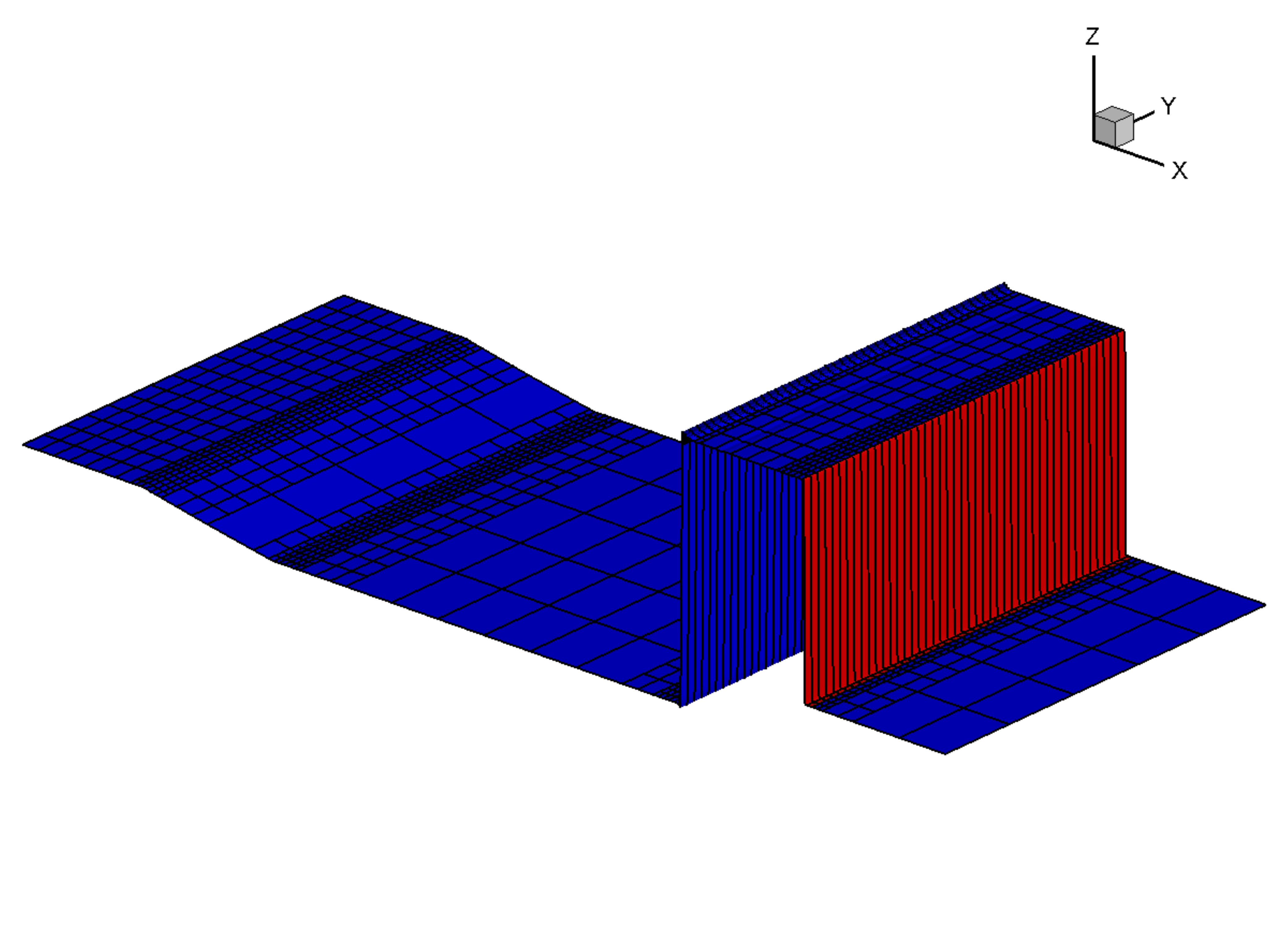}
  \end{tabular}
   \caption{ \label{fig:sodlax3D}
			3D view of the density variable and of the corresponding AMR grid. Top panel: Sod problem at $t_{\text{final}}=0.2$.  Bottom panel: 
      Lax problem at $t_{\text{final}}=0.14$. The limited cells, using the sub-cell ADER-WENO3 finite volume scheme, 
			are highlighted in red, while unlimited DG-$\mathbb{P}_9$ cells are highlighted in blue. 
     }
  \end{center}
\end{figure}

\begin{figure}
  \begin{center} 
  \begin{tabular}{cc} 
    \includegraphics[width=0.47\textwidth]{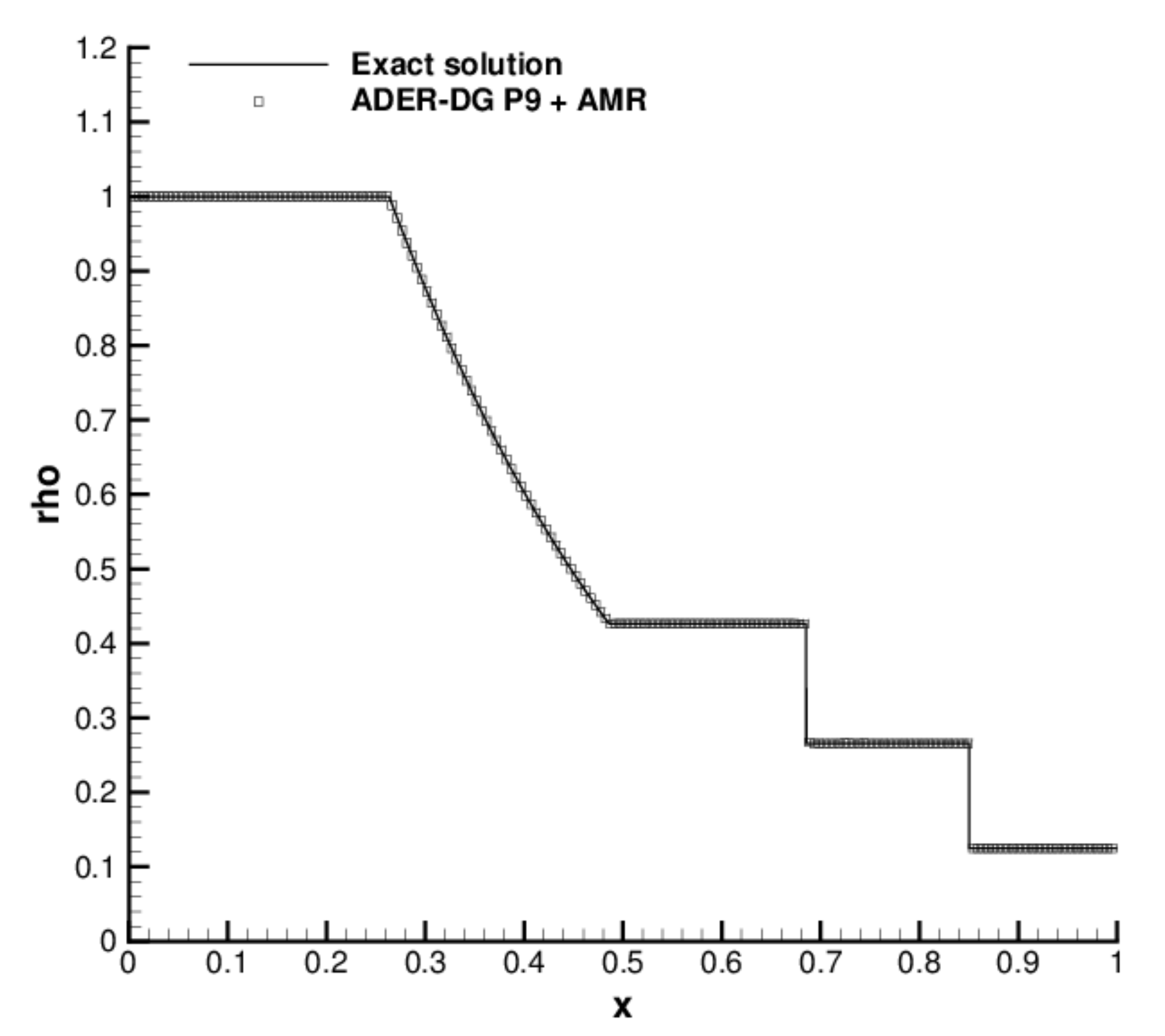}&
    \includegraphics[width=0.47\textwidth]{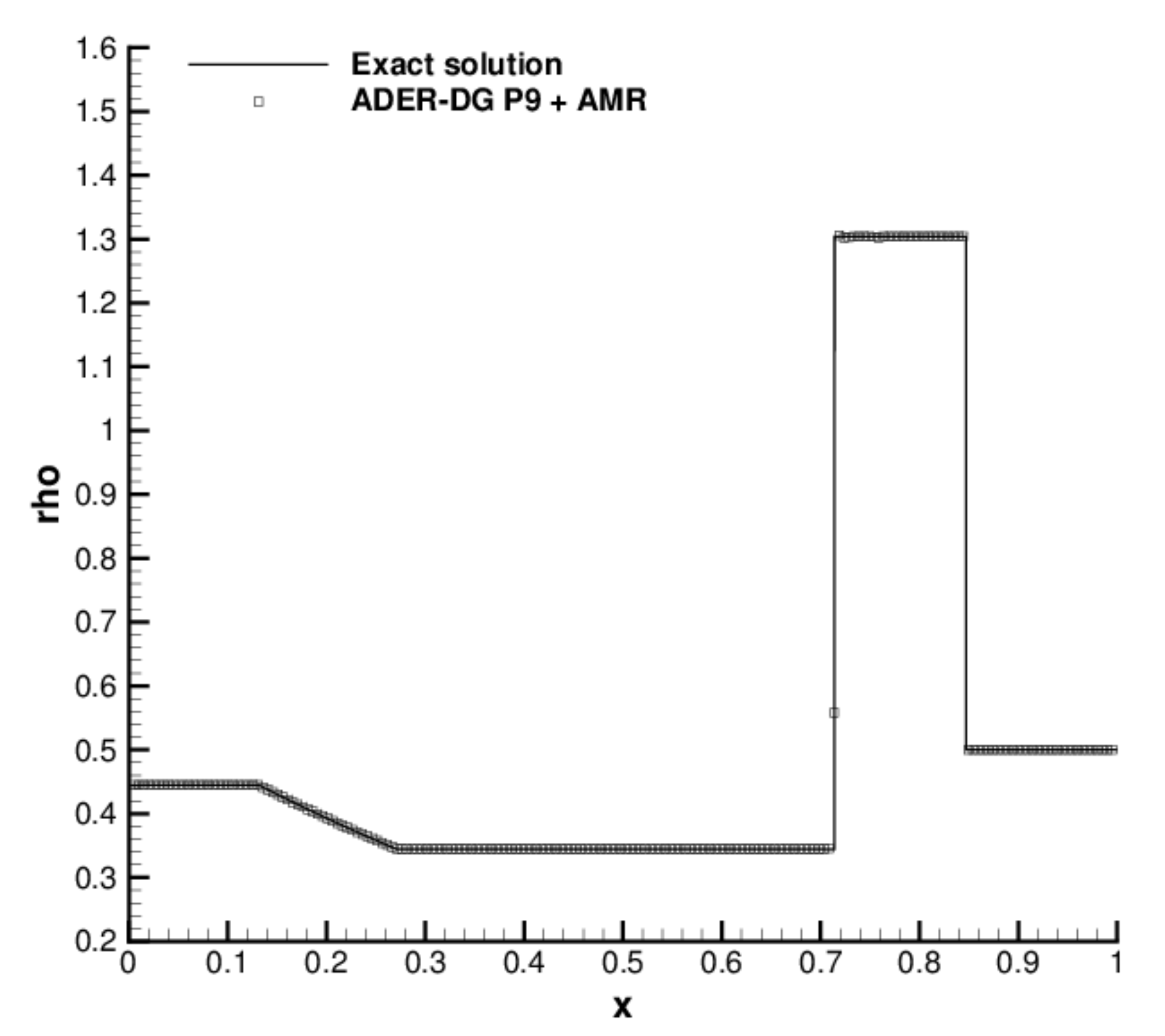}
    \\
    \includegraphics[width=0.47\textwidth]{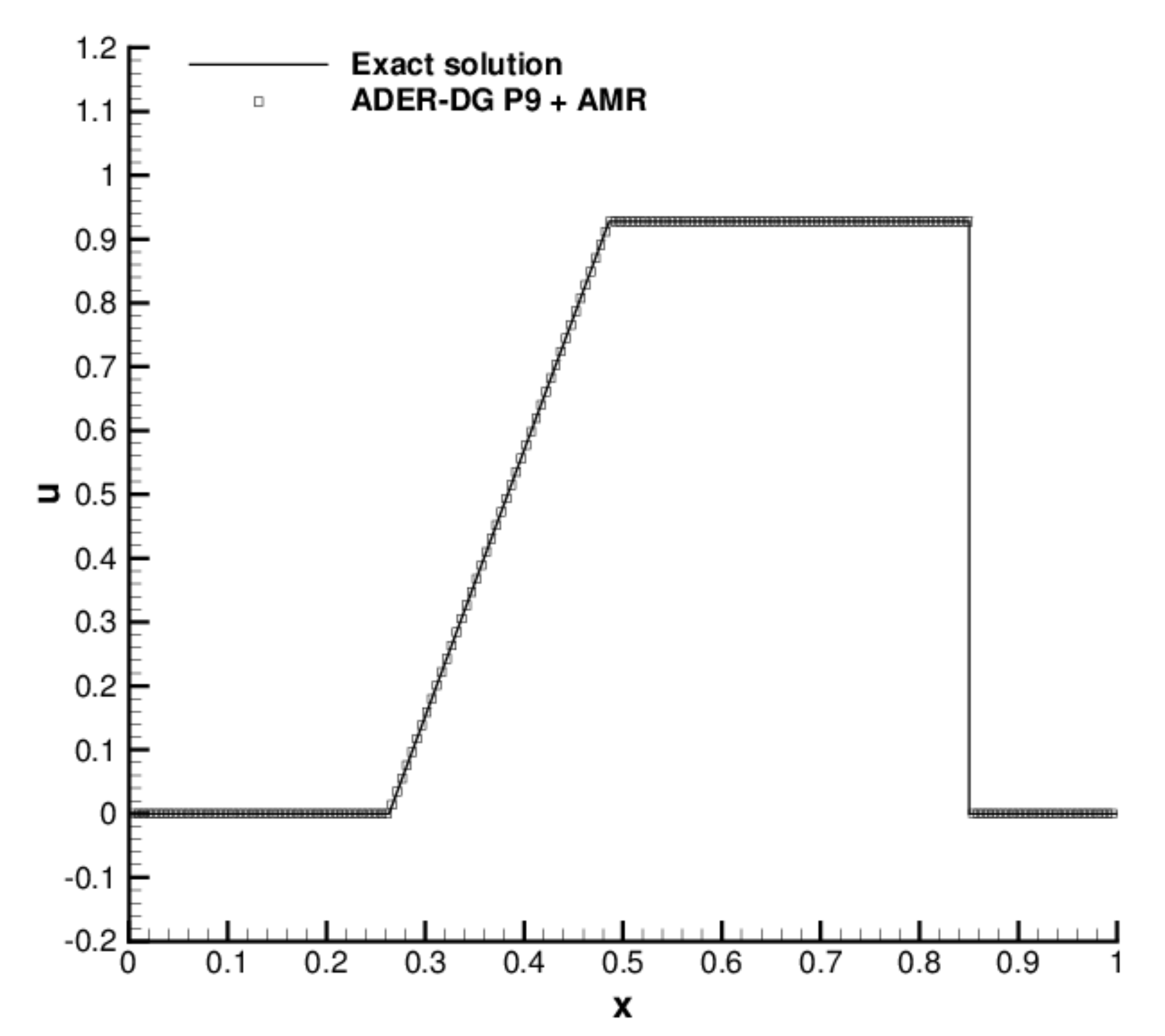}  &
    \includegraphics[width=0.47\textwidth]{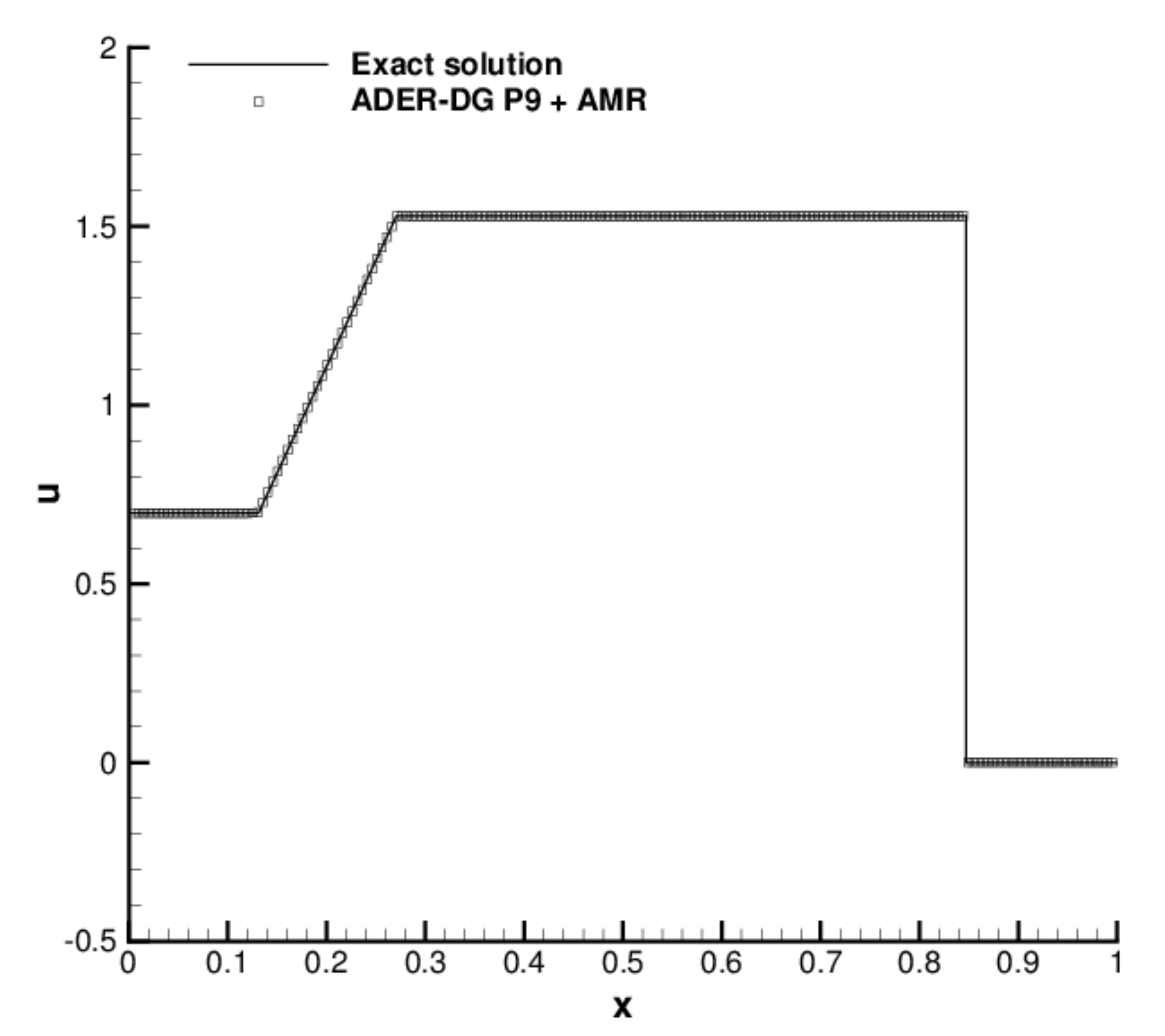} 
    \\
    \includegraphics[width=0.47\textwidth]{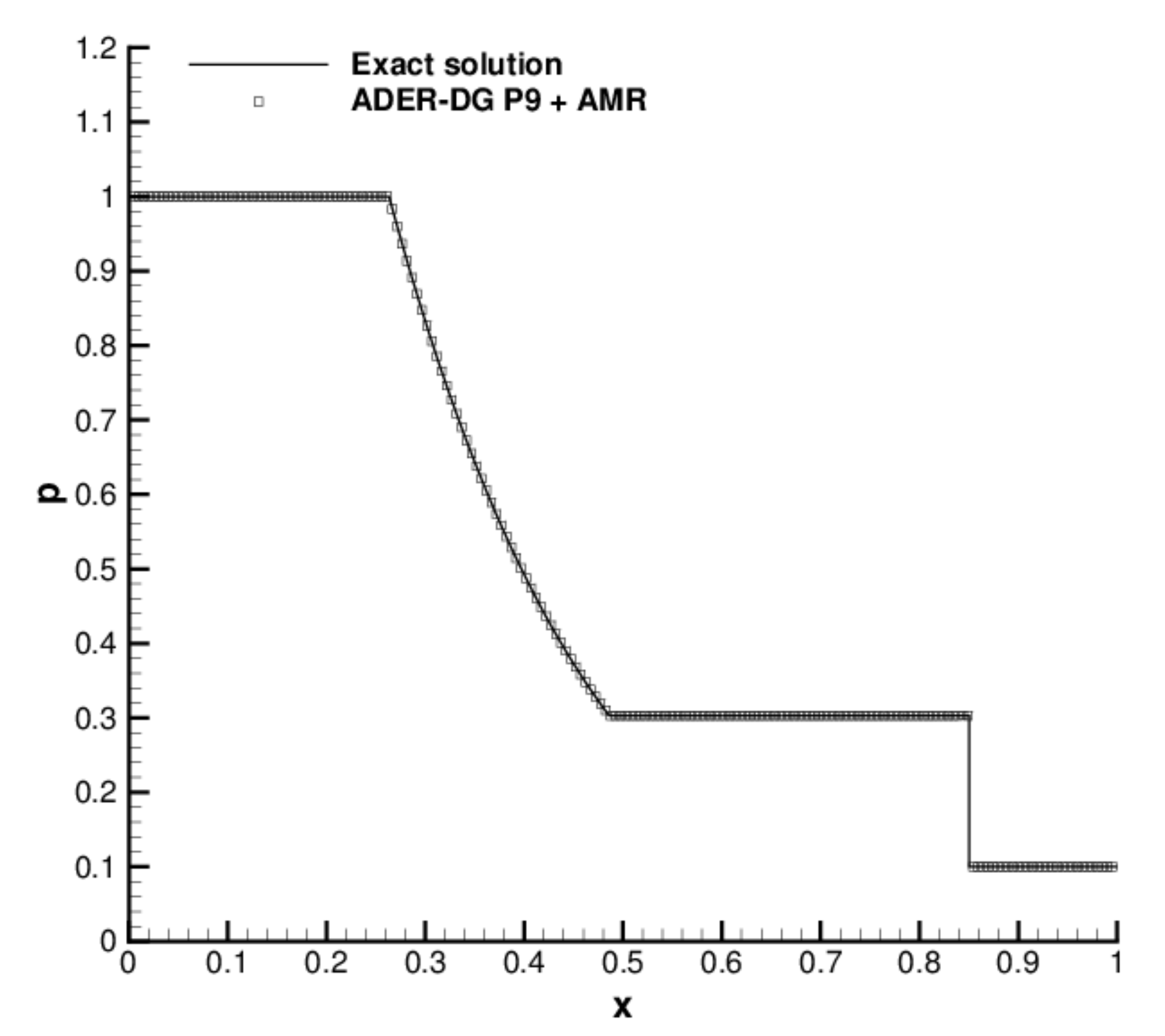}  &
    \includegraphics[width=0.47\textwidth]{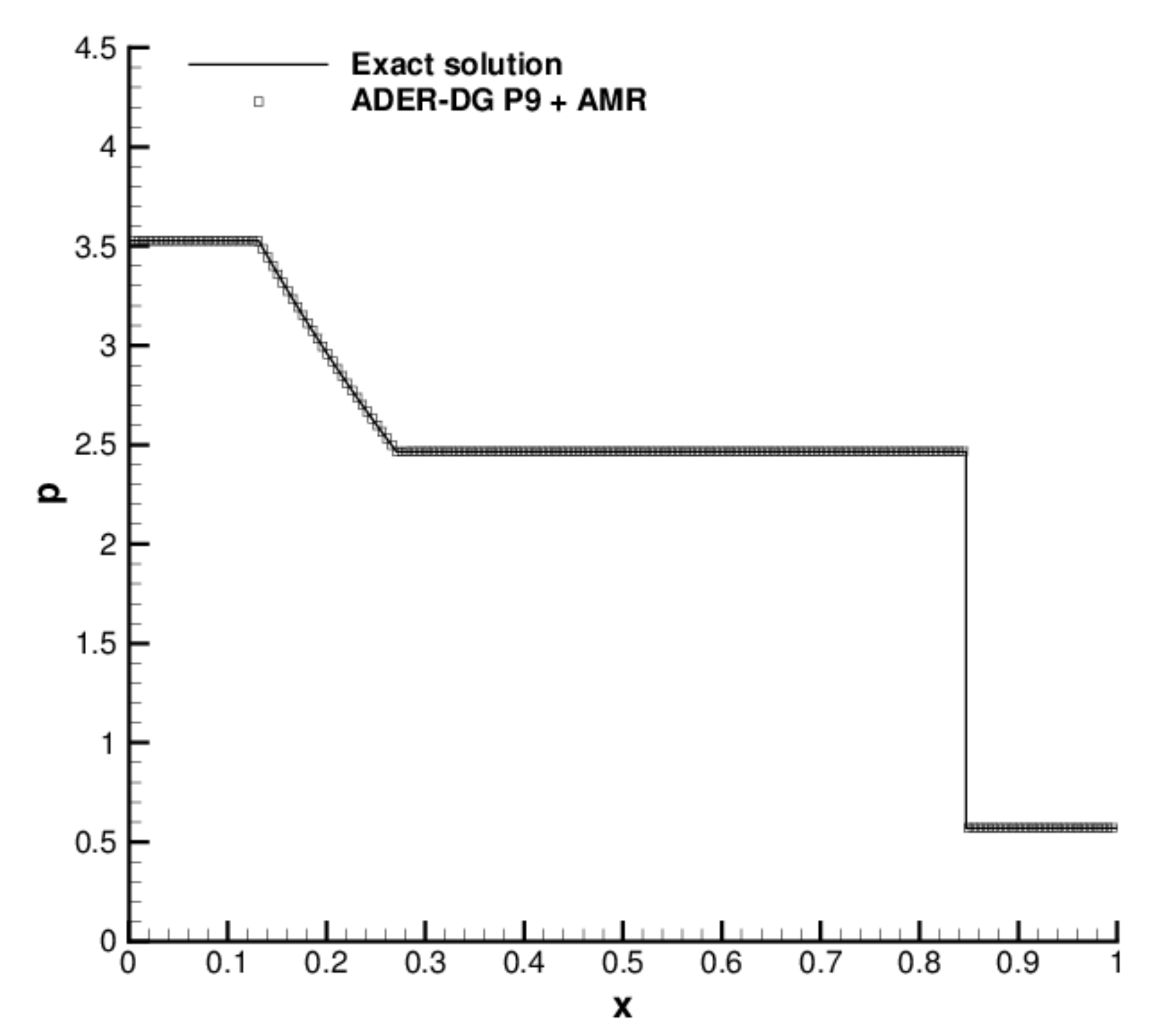}    
  \end{tabular}
   \caption{ \label{fig:sodlax1D}
			Sod shock tube problem (left panels) at $t_{\text{final}}=0.2$ 
      and Lax problem (right panels) at $t_{\text{final}}=0.14$. 
     }
  \end{center}
\end{figure}

Having verified the convergence properties of the ADER-DG-AMR scheme, 
we have considered two classical Riemann problems, proposed by Sod and Lax, with initial conditions given, respectively, by
\begin{equation}
\label{eqn.sod}
(\rho,u,p)_{\rm Sod}= \left\{
\begin{array}{llll}
(1.0,0.0,1.0) &   {\rm if} & x \in [0;0.5] \,, \\ 
(0.125,0.0,0.1) &   {\rm if} & x \in [0.5;1.0] \,,
\end{array} \right.
\end{equation}
and
\begin{equation}
\label{eqn.lax}
(\rho,v,p)_{\rm Lax}= \left\{
\begin{array}{llll}
(0.445,0.698,3.528) &   {\rm if} & x \in [0;0.5] \,, \\ 
(0.5,0.0,0.571) &   {\rm if} & x \in [0.5;1.0] \,.
\end{array} \right.
\end{equation}
The computational domain is actually two-dimensional, but the second direction $y$ acts as a passive one. Moreover, the adiabatic index of the gas is $\gamma=1.4$, and the final time of the simulation is $t_{\text{final}}=0.2$ for Sod's problem, while it is $t_{\text{final}}=0.14$ for Lax's. Both tests have been solved using the ADER-DG-$\mathbb{P}_9$ scheme, supplemented with our \aposteriori 
ADER-WENO3 finite volume sub-cell limiter. The initial grid is composed of $N_x\times N_y=20\times 5$ cells, which are then adaptively refined using $\mathfrak{r}=3$ and $\ell_{\rm max}=2$. The results of our calculations, for which we have used the Osher flux \cite{OsherUniversal}, are reported in Figs.~\ref{fig:sodlax3D}--\ref{fig:sodlax1D}. 
Figure \ref{fig:sodlax3D}, in particular, shows the three-dimensional view of the solution by plotting the corresponding polynomials, highlighted in blue (for the unlimited cells) and in red (for the limited cells) according to our standard convention. We recall that the blue polynomials really represent the DG polynomials within each cell, while in the red cells we visualize the data as a piecewise 
linear interpolation of the subcell averages, produced by the subcell limiter.
As in Fig.~4 of \cite{Dumbser2014}, in both the tests the contact discontinuity is resolved within one single cell, which, due to our AMR algorithm, in the present case is always at the maximum level of refinement. We further note that the contact wave is unlimited (blue). This is due to the fact that after a certain time our ADER-DG scheme recognizes this linear degenerate wave as a 
\textit{smooth feature}, after the initial smoothing of the contact discontinuity by the subcell limiter and the Riemann solver. The right propagating shock, on the contrary, is always 
limited (red), as expected, and it is very  sharply resolved. 
In Fig.~\ref{fig:sodlax1D} we have instead reported the comparison of the exact solution of the Riemann problem \cite{toro-book} with the numerical solution for a few representative variables,  
extracted from the polynomial data representation of the DG scheme -or the subcell limiter- along a 1D line of 200 equidistant sample points. The agreement between numerical and exact solution 
is excellent. Finally, for this test we have also performed a profiling analysis to quantify the relative computational costs of the subcell limiter. 
In a representative simulation using the ADER-DG-$\mathbb{P}_2$ scheme, with approximately $15\%$ of the cells that are limited, 
the overhead with respect to the unlimited DG scheme amounts to a factor $\approx 1.5$ in terms of CPU time.

\subsubsection{Double Mach reflection problem} 

\begin{figure}
  \begin{center} 
  \begin{tabular}{cc} 
    \includegraphics[width=0.47\textwidth]{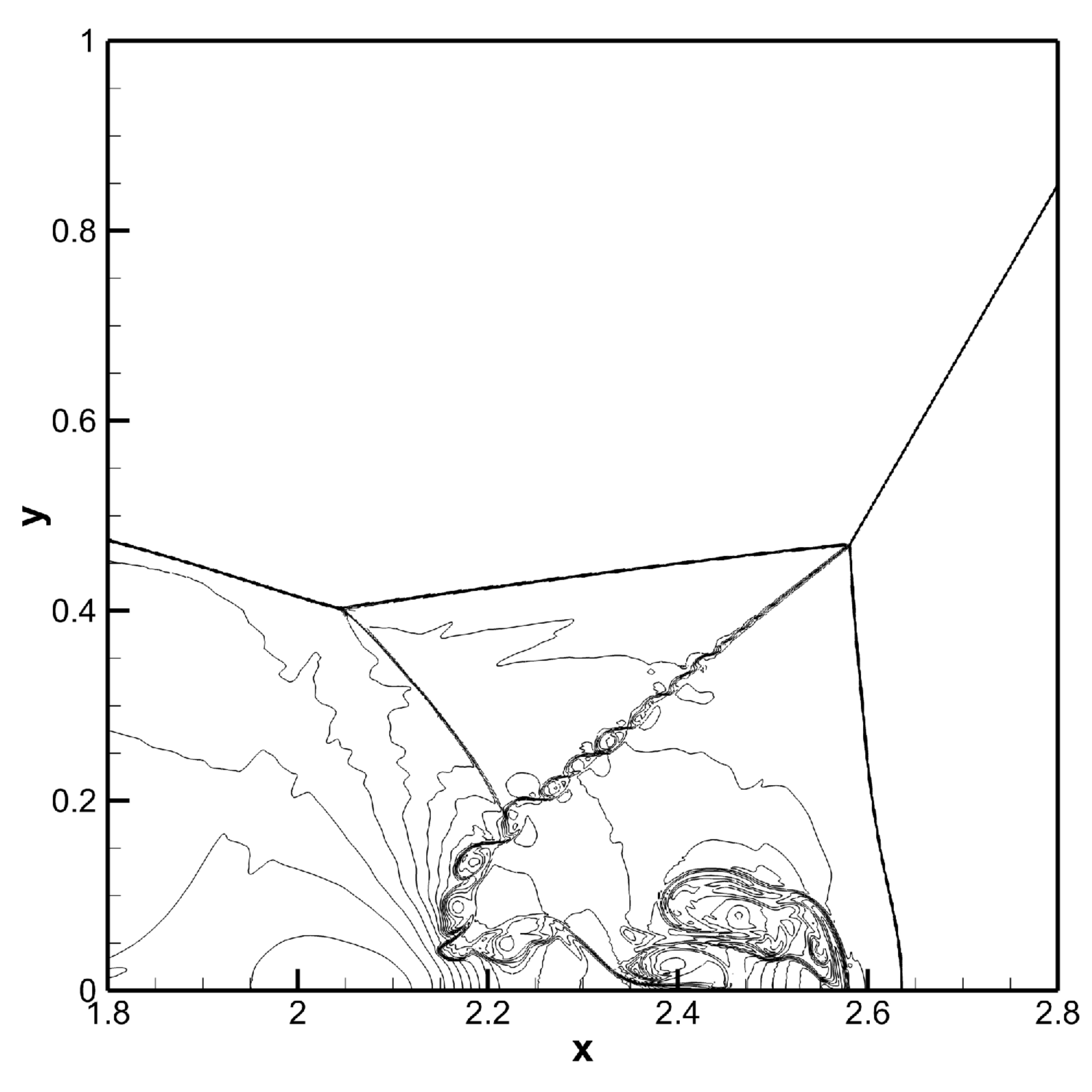} &
    \includegraphics[width=0.47\textwidth]{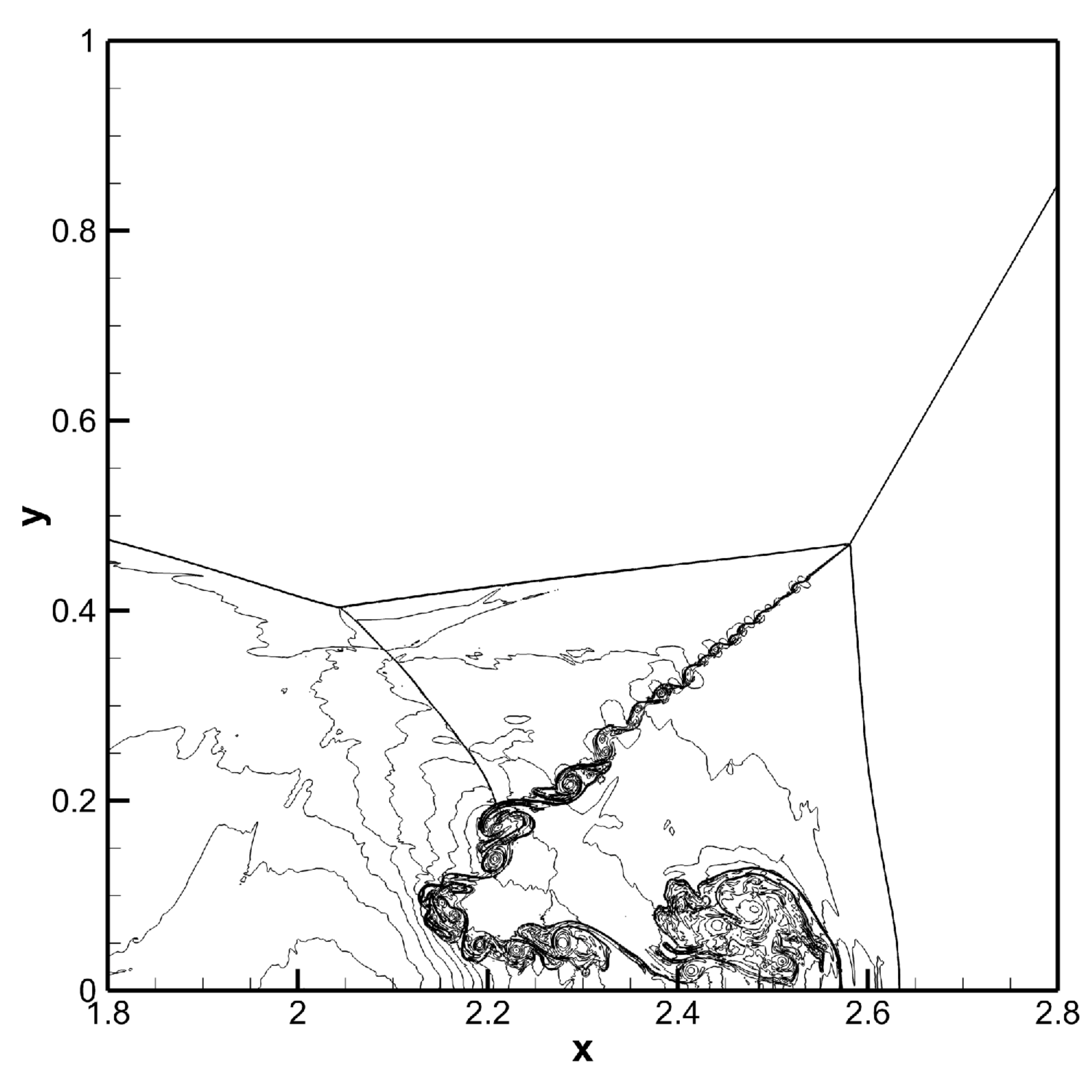} \\
    \includegraphics[width=0.47\textwidth]{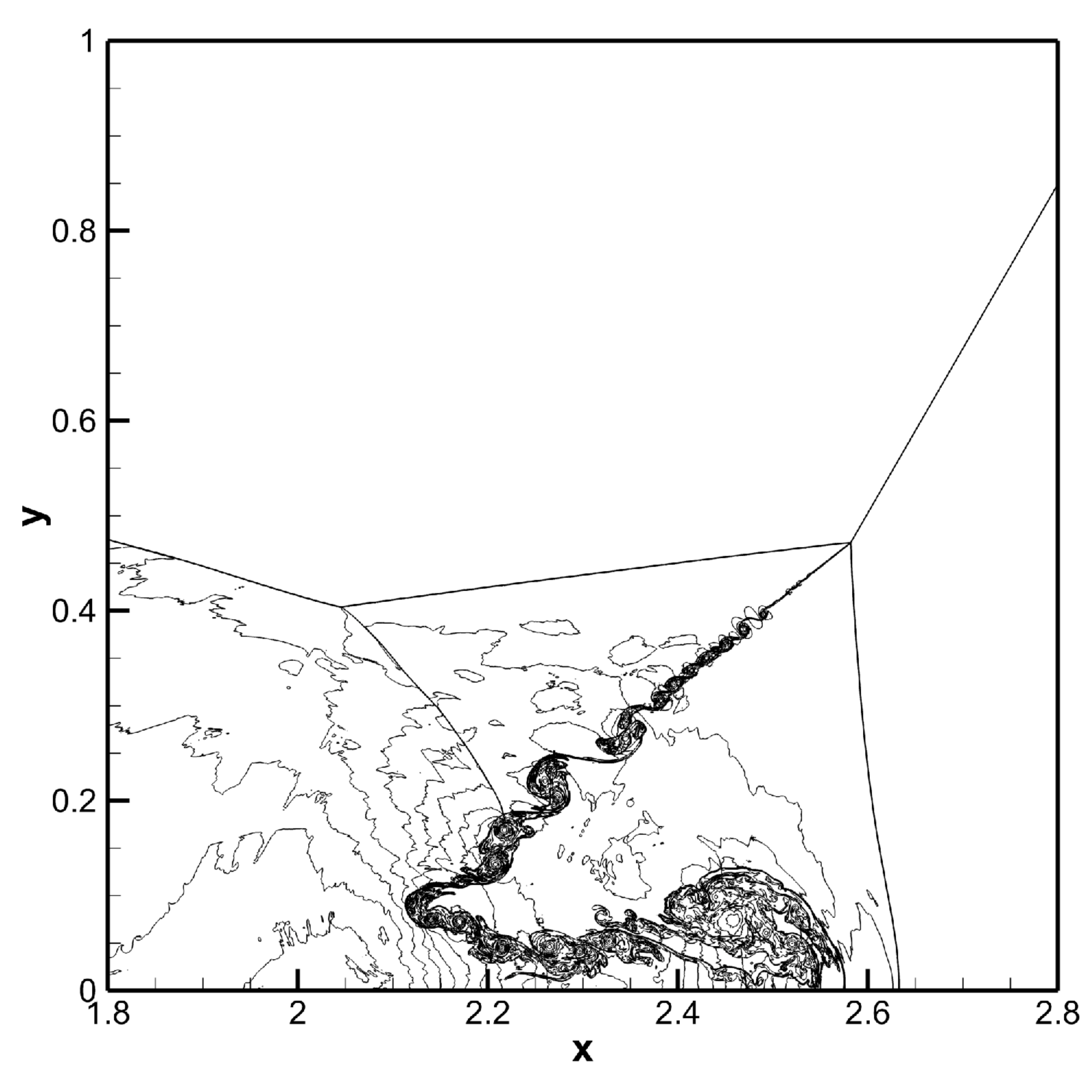}  &
    \includegraphics[width=0.47\textwidth]{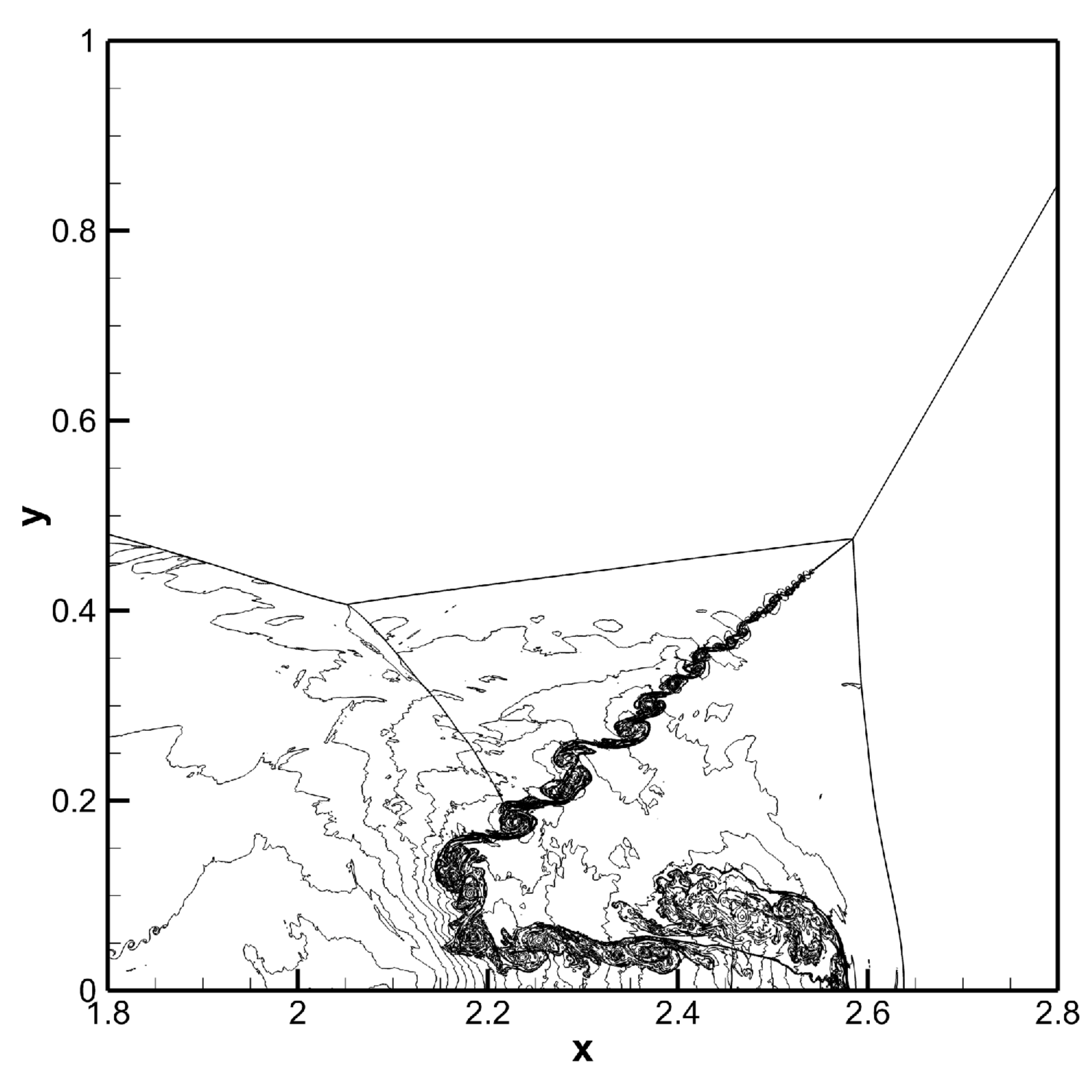}  
  \end{tabular}
   \caption{ \label{fig:DMR2D}
			Zooms of the interaction zone for the double Mach reflection problem at $t=0.2$. Equidistant contour lines of the density variable are shown.
			Top left: AMR-ADER-DG-$\mathbb{P}_2$ with initial $75\times25$ grid. Top right: AMR-ADER-DG-$\mathbb{P}_5$ with initial $75\times25$ grid.
			Bottom left: AMR-ADER-DG-$\mathbb{P}_8$ with initial $75\times25$ grid. Bottom right: AMR-ADER-DG-$\mathbb{P}_5$ with initial $150\times50$ grid.
     }
  \end{center}
\end{figure}

\begin{figure}
  \begin{center} 
  \begin{tabular}{cc} 
    \includegraphics[width=0.47\textwidth]{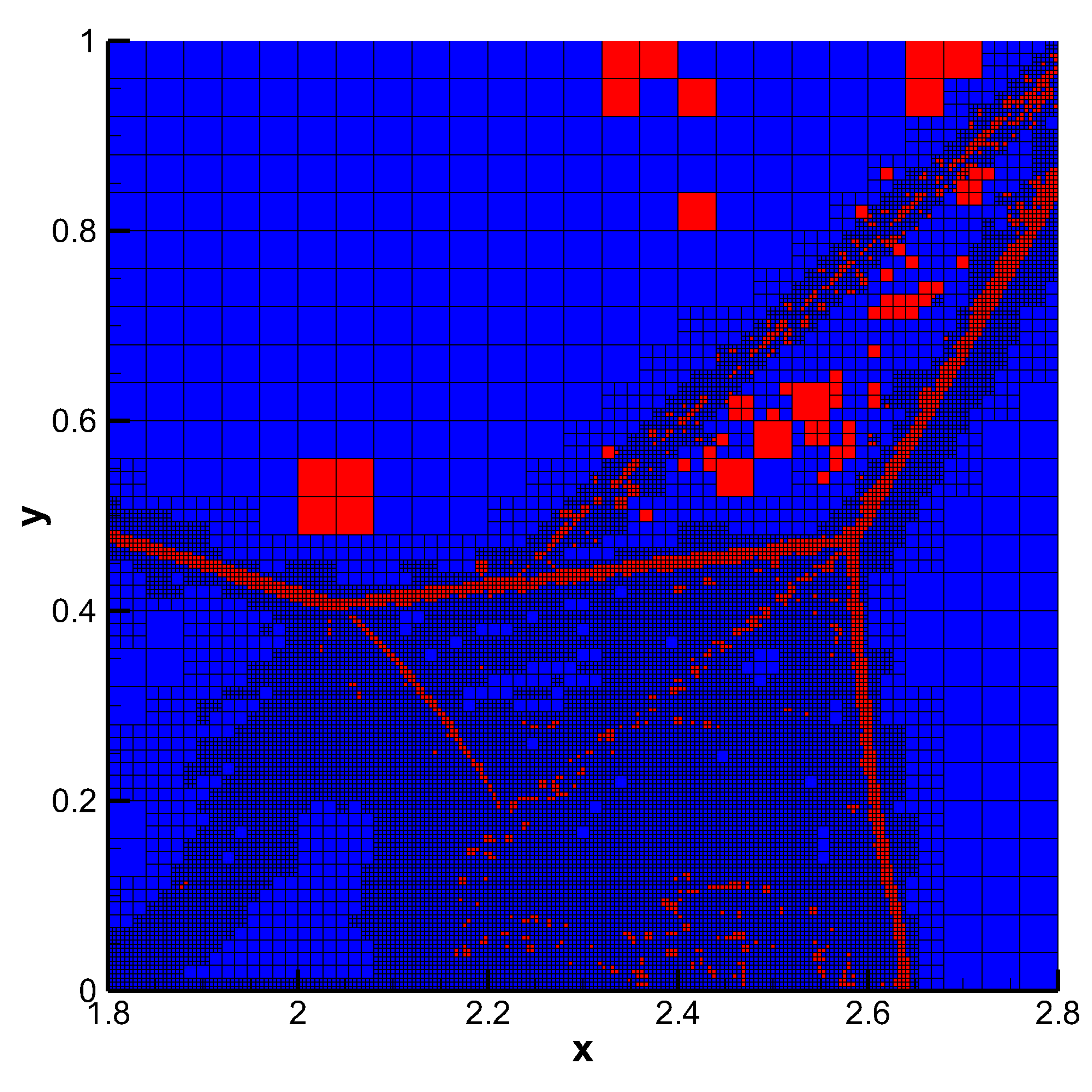} &
    \includegraphics[width=0.47\textwidth]{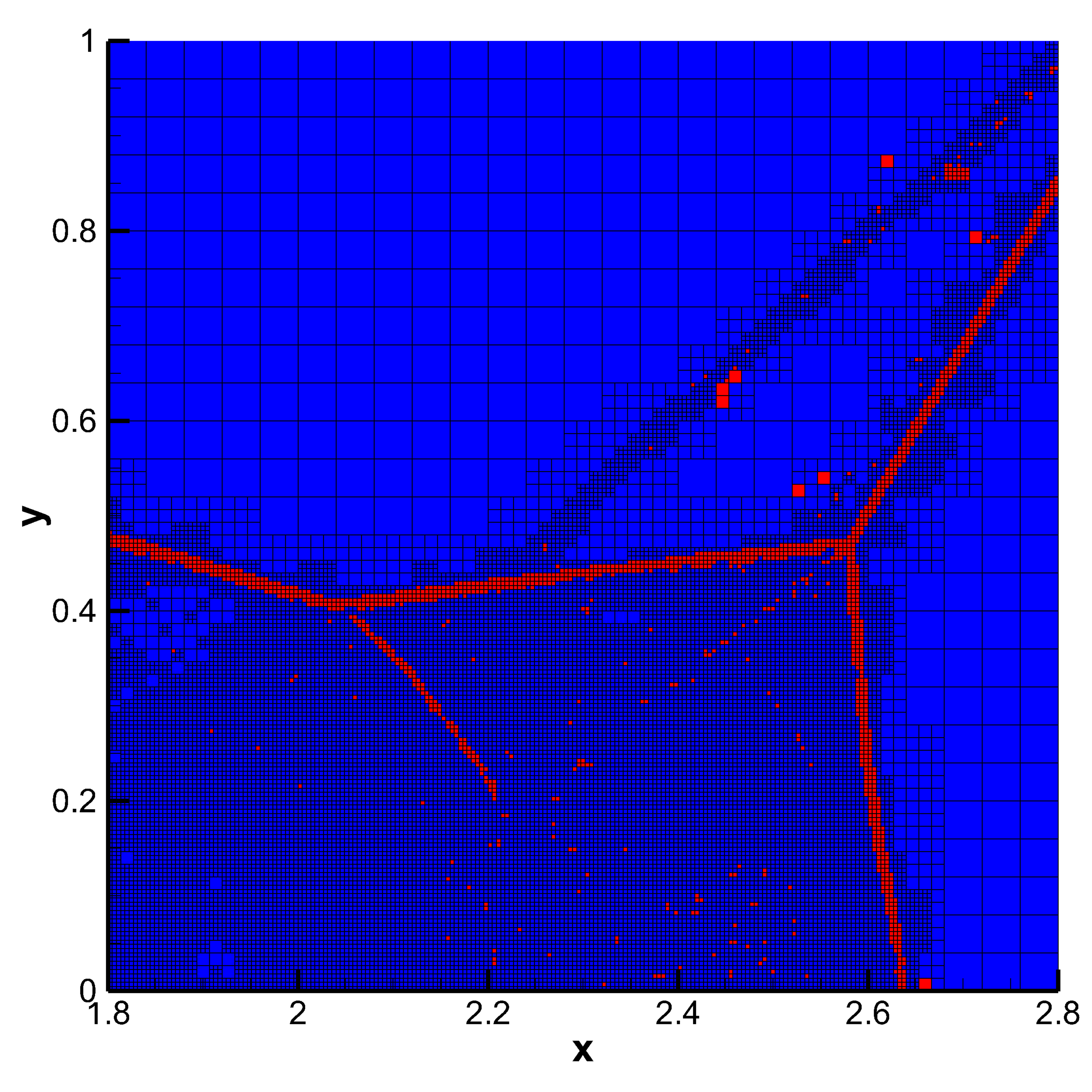} \\
    \includegraphics[width=0.47\textwidth]{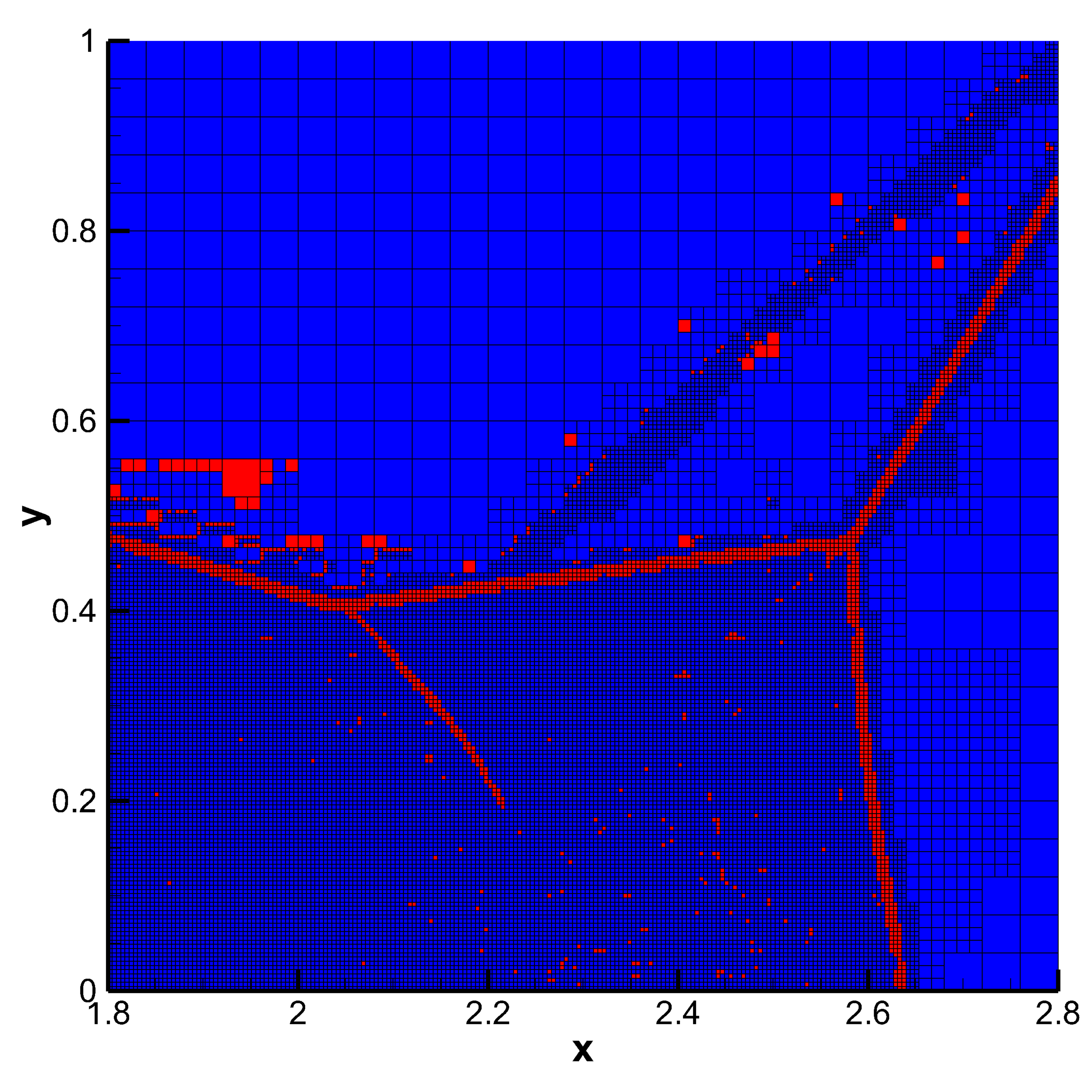}  &
    \includegraphics[width=0.47\textwidth]{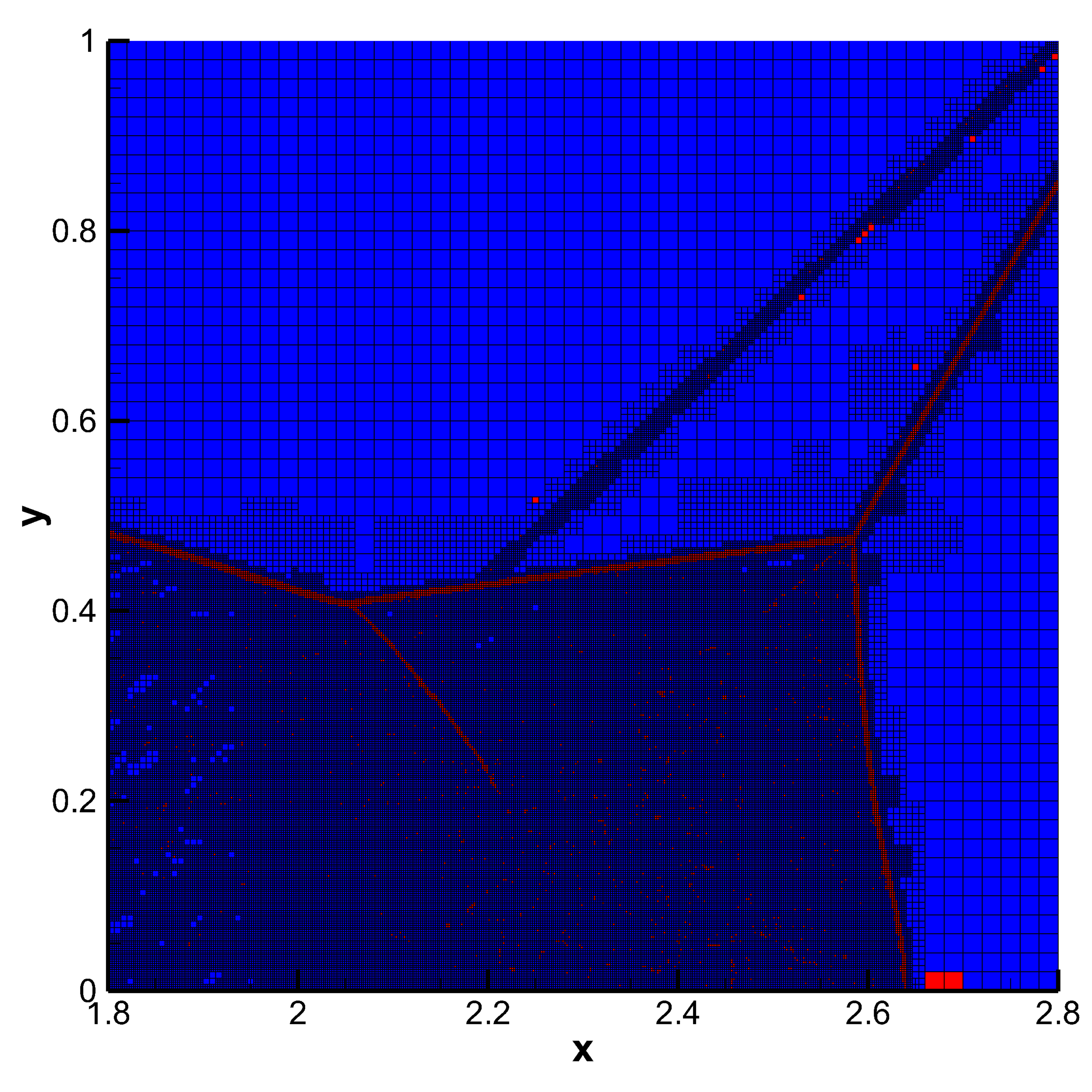}  
  \end{tabular}
   \caption{ \label{fig:DMR2D.grid}
	Zooms of the interaction zone for the double Mach reflection problem at $t=0.2$. The AMR grid and the limited cells (highlighted in red) are shown.
			Top left: ADER-DG-$\mathbb{P}_2$ with initial $75\times25$ grid. Top right: ADER-DG-$\mathbb{P}_5$ with initial $75\times25$ grid.
			Bottom left: ADER-DG-$\mathbb{P}_8$ with initial $75\times25$ grid. Bottom right: ADER-DG-$\mathbb{P}_5$ with initial $150\times50$ grid.
	     }
  \end{center}
\end{figure}

\begin{figure}
  \begin{center} 
  \begin{tabular}{cc} 
    \includegraphics[width=0.48\textwidth]{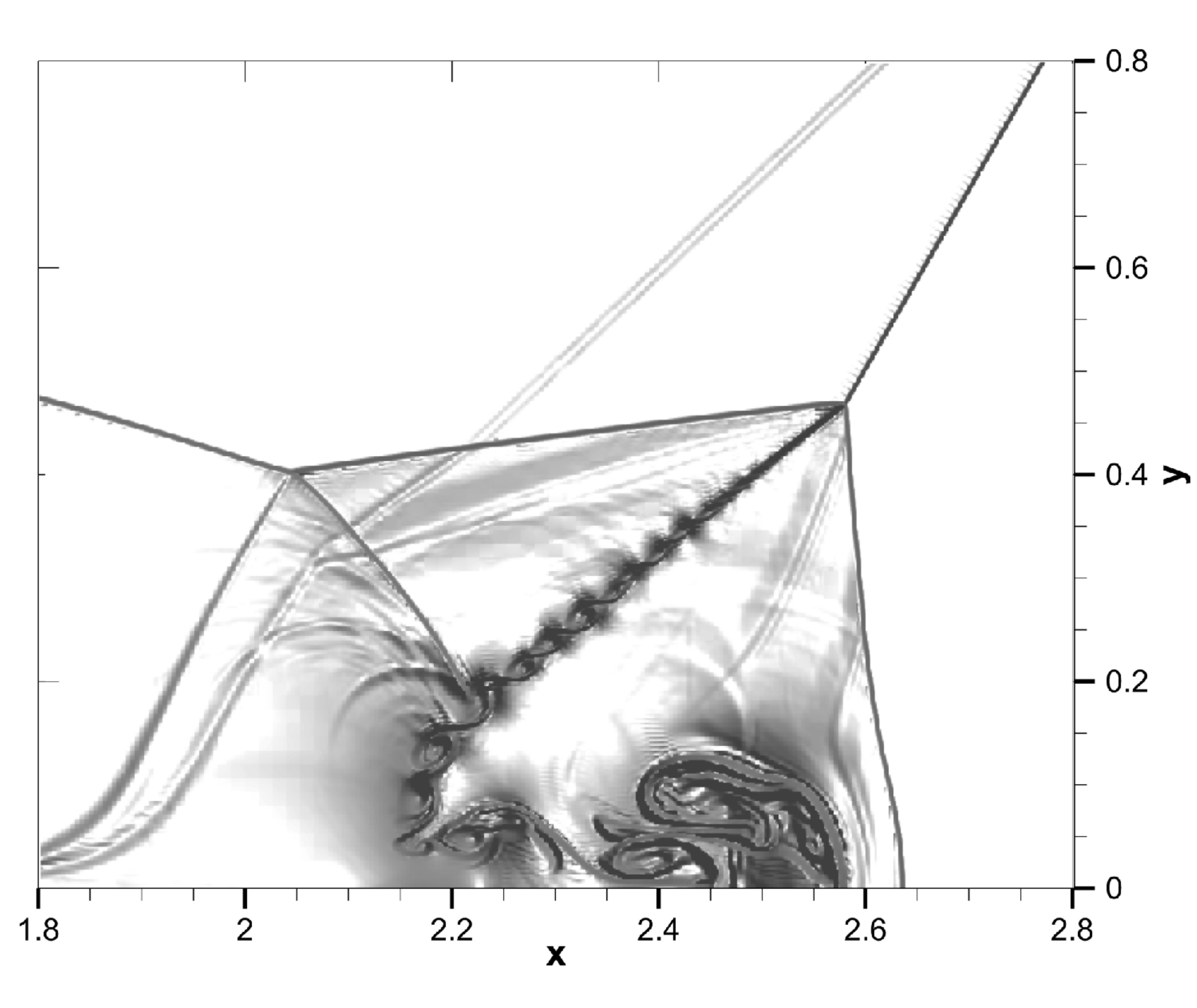} &
    \includegraphics[width=0.48\textwidth]{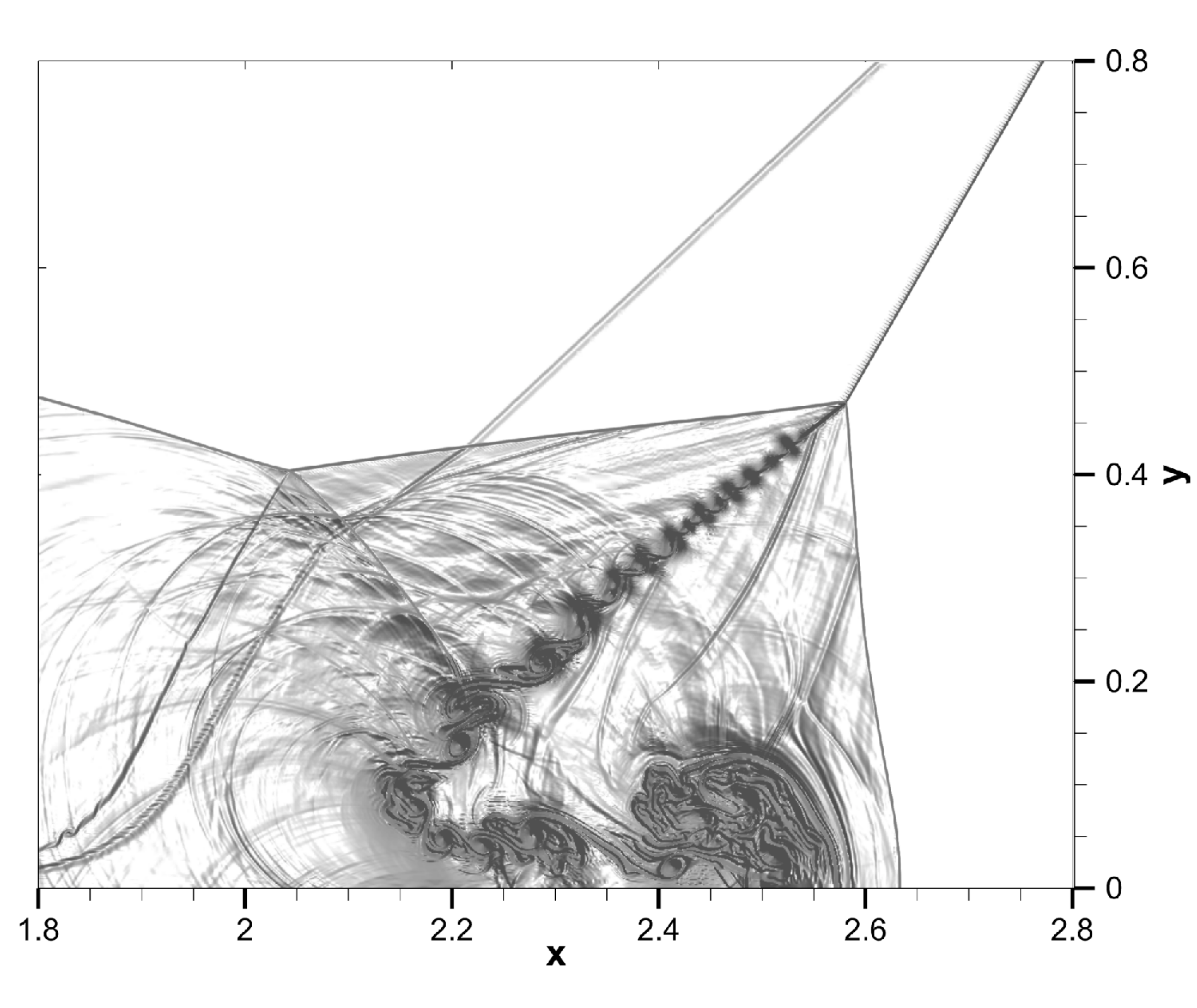} \\
    \includegraphics[width=0.48\textwidth]{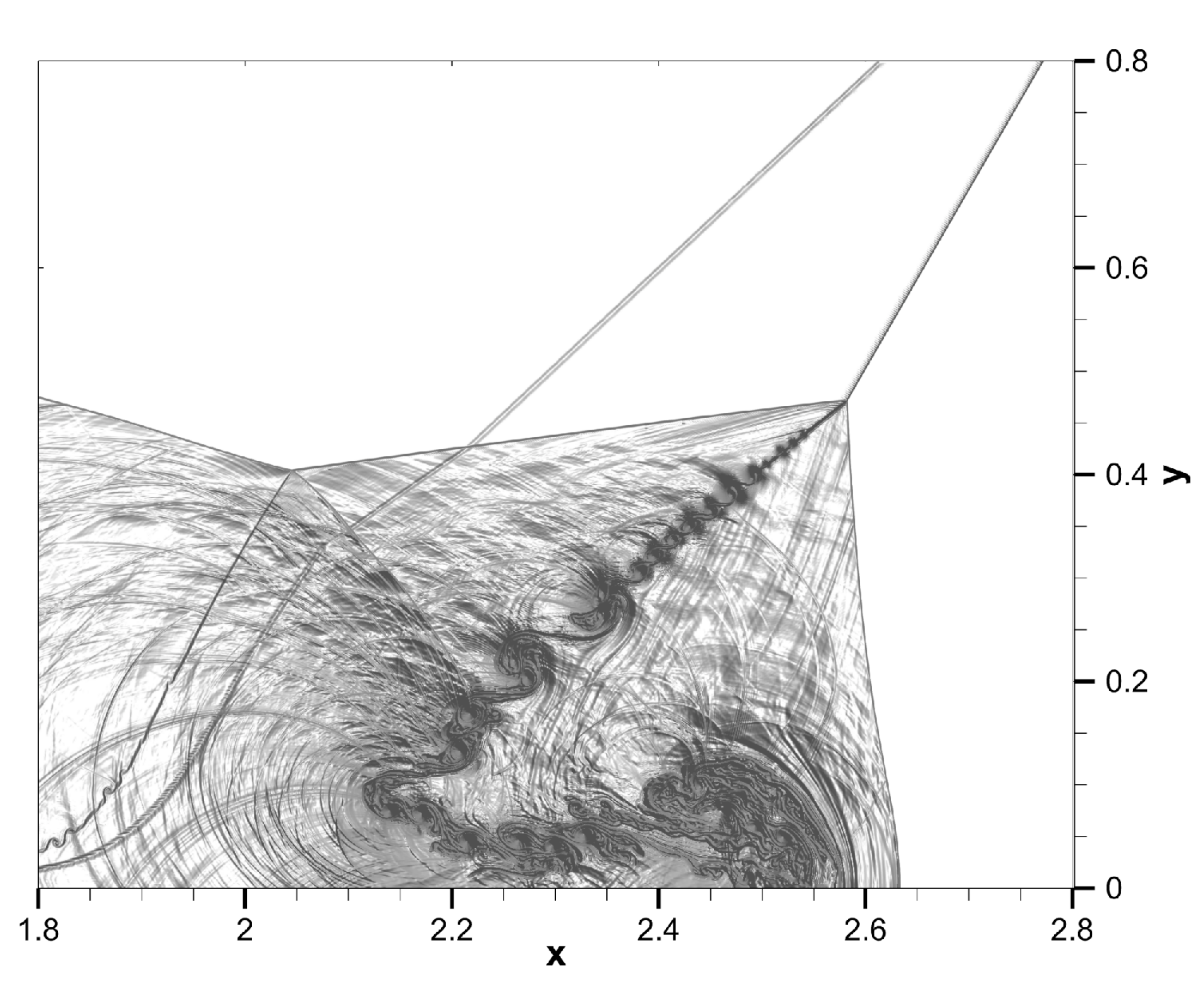}  &
    \includegraphics[width=0.48\textwidth]{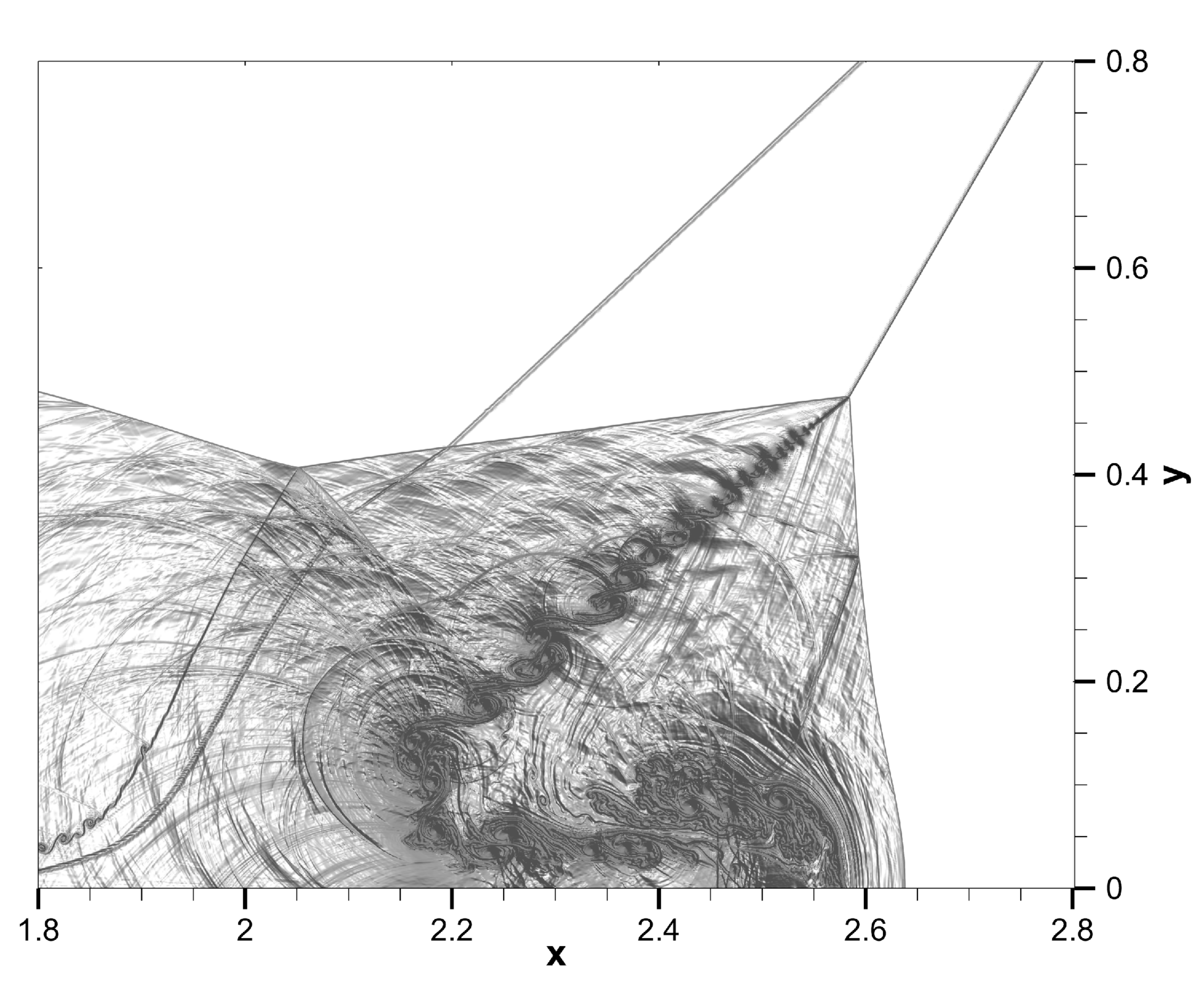}  
  \end{tabular}
   \caption{ \label{fig:DMR3D}
			Schlieren image of the density variable for the double Mach reflection problem at $t=0.2$. 
			Top left: ADER-DG-$\mathbb{P}_2$ with initial $75\times25$ grid. Top right: ADER-DG-$\mathbb{P}_5$ with initial $75\times25$ grid.
			Bottom left: ADER-DG-$\mathbb{P}_8$ with initial $75\times25$ grid. Bottom right: ADER-DG-$\mathbb{P}_5$ with initial $150\times50$ grid.
	     }
  \end{center}
\end{figure}

A complex test problem in two space dimensions which contains a variety of 
waves such as strong shock waves, contact waves and shear waves, 
we have considered the so called {\em double Mach reflection problem}, which was first proposed   
in \cite{woodwardcol84}. The initial conditions are given by a right-moving shock wave with a Mach number $M=10$,
which intersects the $x-$ axis at  
$x=1/6$ with an inclination angle of $\alpha=60^{\circ}$. In order to provide the physical states ahead and behind the shock, it is necessary to solve the Rankine--Hugoniot conditions, which provide
\begin{eqnarray}
(\rho, u, v, p)( \x,t=0) =
\left\{
\begin{array}{cll}
  \frac{1}{\gamma}(8.0, 8.25, 0.0, 116.5), \quad & \text{ if } & \quad x'<0.1, \\
  (1.0, 0.0, 0.0, \frac{1}{\gamma}),       \quad & \text{ if } & \quad x'\geq 0.1, 
\end{array}
\right.
\end{eqnarray}
where $x' = (x - 1/6) \cos\alpha - y \sin\alpha$ is the coordinate in the rotated frame, while $\gamma=1.4$. 
The boundary conditions on the left side and on the right side are just given by
inflow and outflow, while on the bottom we have used reflecting boundary conditions. 
On the other hand, the boundary conditions on the top require some more attention,
since we need to impose the exact solution of an isolated moving oblique shock wave with the same shock Mach 
number $M_s=10$.
The computational domain is given by $\Omega = [0;3.0] \times [0;1]$, which is covered by an initial uniform grid composed of
$75\times25$ cells. For our simulations, the Rusanov flux has been used and AMR is activated with $\ell_{\rm max}=2$ and $\mathfrak{r}=3$. 
The results of our calculations at time $t=0.2$
are reported in Figs.~\ref{fig:DMR2D}-\ref{fig:DMR3D}, for which we have used 
three different schemes: ADER-DG-$\mathbb{P}_N$ with $N=2,5,8$. In all these figures we have zoomed into the interaction zone with $1.8\leq x\leq 2.8$ in order to 
highlight the differences among the orders of accuracy. Moreover, the bottom right panel in each of these figures refers to a configuration with 
a finer initial grid, composed of $150\times50$ cells. Fig.~\ref{fig:DMR2D}, in particular, shows the contour lines of the density. Fig.~\ref{fig:DMR2D.grid} shows the AMR grid and the troubled cells, highlighted in red, which required 
the activation of the limiter. Finally, Fig.~\ref{fig:DMR3D} reports the Schlieren images of the density.
There are a number of comments that can be made about these results. First, and mostly obvious, all DG schemes can detect the shock waves very well.
On the other hand, by increasing the order of accuracy, the vortex-type flow structures manifest a larger and richer rolling-up, especially in the transition from ADER-DG-$\mathbb{P}_2$
to ADER-DG-$\mathbb{P}_5$. 
Secondly, the largest number of troubled cells, including false-positive troubled cells, is present for the lowest order scheme, i.e. the ADER-DG-$\mathbb{P}_2$, and it is concentrated along the shocks, while leaving the vortex-type flow structures unaffected.
This is reassuring, since it indicates that
higher order DG schemes have \textit{better subcell resolution} capabilities. Last but not least we would like to note that the vortices generated by the rolling of the shear waves create \textit{sound waves}, which travel through the computational domain. Although these simulations do not contain physical viscosity, and as such the vortex generation and rolling is only controlled by numerical viscosity, 
we can deduce from our numerical results that the novel scheme is able to resolve shock waves properly, as well as shear waves, vortex structures and sound waves.  


\subsubsection{Forward facing step} 

\begin{figure}[!t] 
  \begin{center} 
      \includegraphics[width=0.85\textwidth]{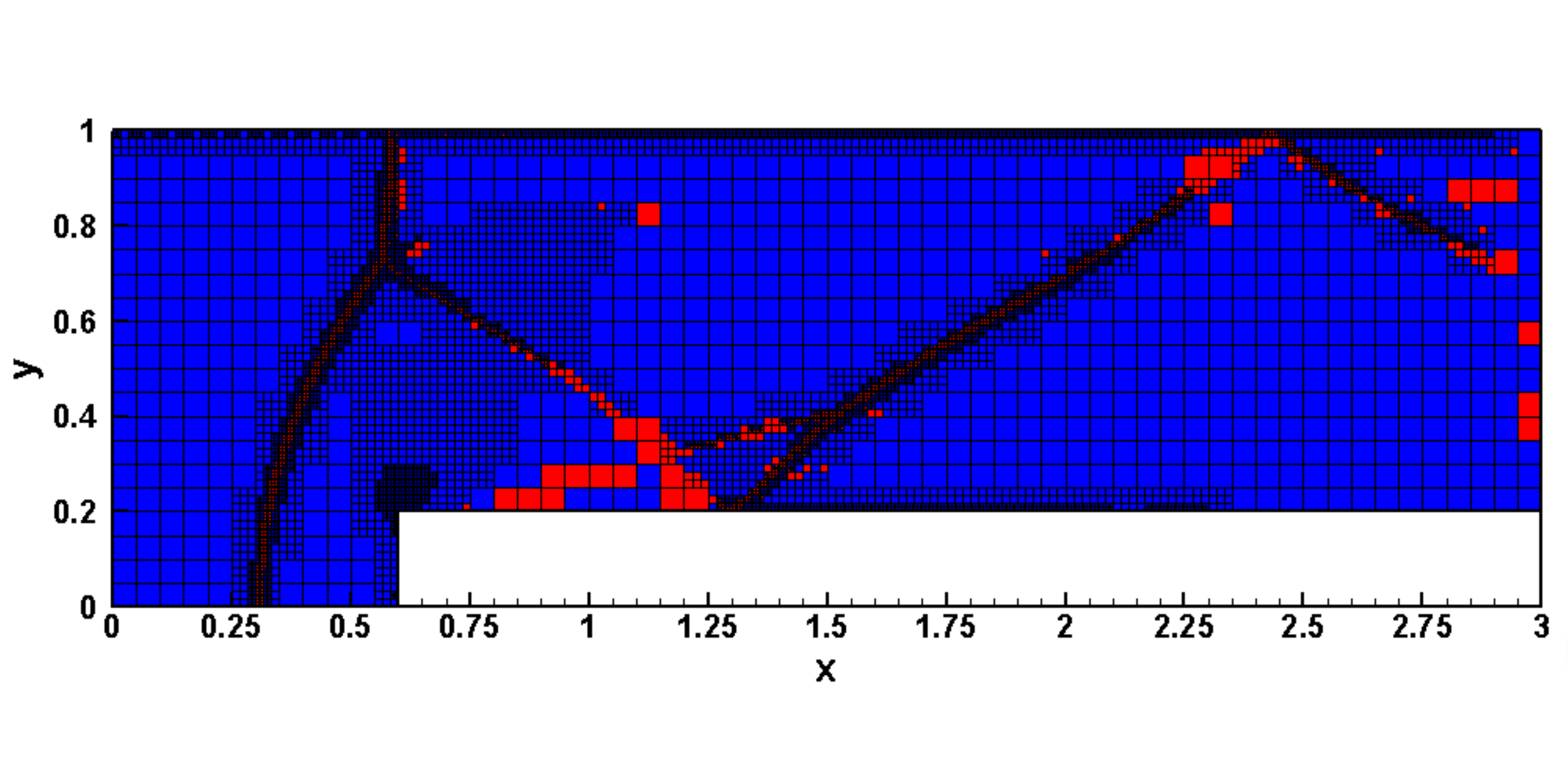} \\ 
      \includegraphics[width=0.85\textwidth]{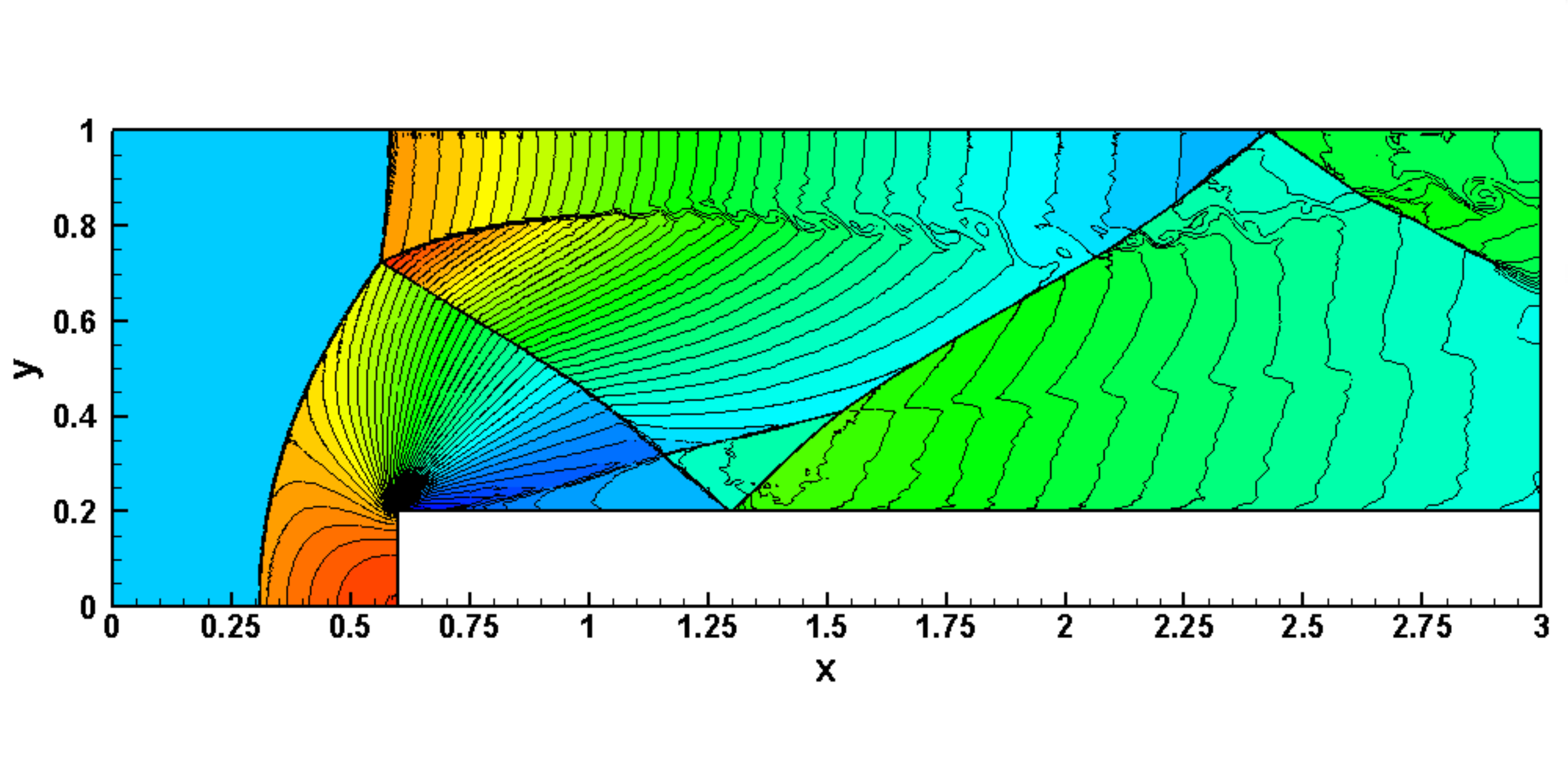}  
    \caption{ \label{fig:FFS}
      Forward facing step problem using ADER-DG-$\mathbb{P}_5$ with \aposteriori ADER-WENO3 sub-cell limiter. Top: 2D view of the AMR grid together with limited 
			cells (red) and unlimited cells (blue). 
			Bottom: 41 equidistant density contour levels in the interval $[0.1;4.5]$.} 
  \end{center}
\end{figure}

The forward facing step problem
is a classical test, often referred to as the \emph{Mach 3 wind tunnel test}, which was proposed for the first time in \cite{woodwardcol84}. We take as computational domain $\Omega=[0;3]\times[0;1]\backslash[0.6;3]\times[0;0.2]$. The initial conditions are given by 
a uniform flow moving to the right with Mach number $M=3$, $\rho=1$, $p=1/\gamma$, $u=3$, $v=0$, and adiabatic index $\gamma=1.4$.
The final time of simulation is $t=4.0$. Regarding the boundary conditions, 
we have used reflecting boundaries at the lower and upper parts of the numerical domain,
while inflow boundary conditions are imposed at the entrance and outflow boundary conditions at the exit. 
Figure~\ref{fig:FFS} represents the numerical solution obtained using the ADER-DG-$\mathbb{P}_5$ scheme 
with \aposteriori ADER-WENO3 sub-cell limiter. 
The panel on the top is a 2D view of the AMR grid showing, as usual, in red the limited cells and in blue the unlimited ones. 
The bottom panel, on the other hand, is a contour plot with 
41 equidistant density contour levels in the interval $[0.1;4.5]$. 
The mesh at the coarsest level has $150\times 50$ cells, which is subsequently refined using 
AMR parameters $\ell_{\rm max}=2$ and $\mathfrak{r}=4$, corresponding to a uniform grid composed of
$2400\times 800$ cells. 
It can be appreciated that there is a very good resolution of the physical instability and also it can be observed that both AMR and sub-cell limiter act where they are needed.
%
\subsubsection{2D Riemann problems} 
\begin{figure}
\begin{center}
\begin{tabular}{lr}
\includegraphics[width=0.35\textwidth]{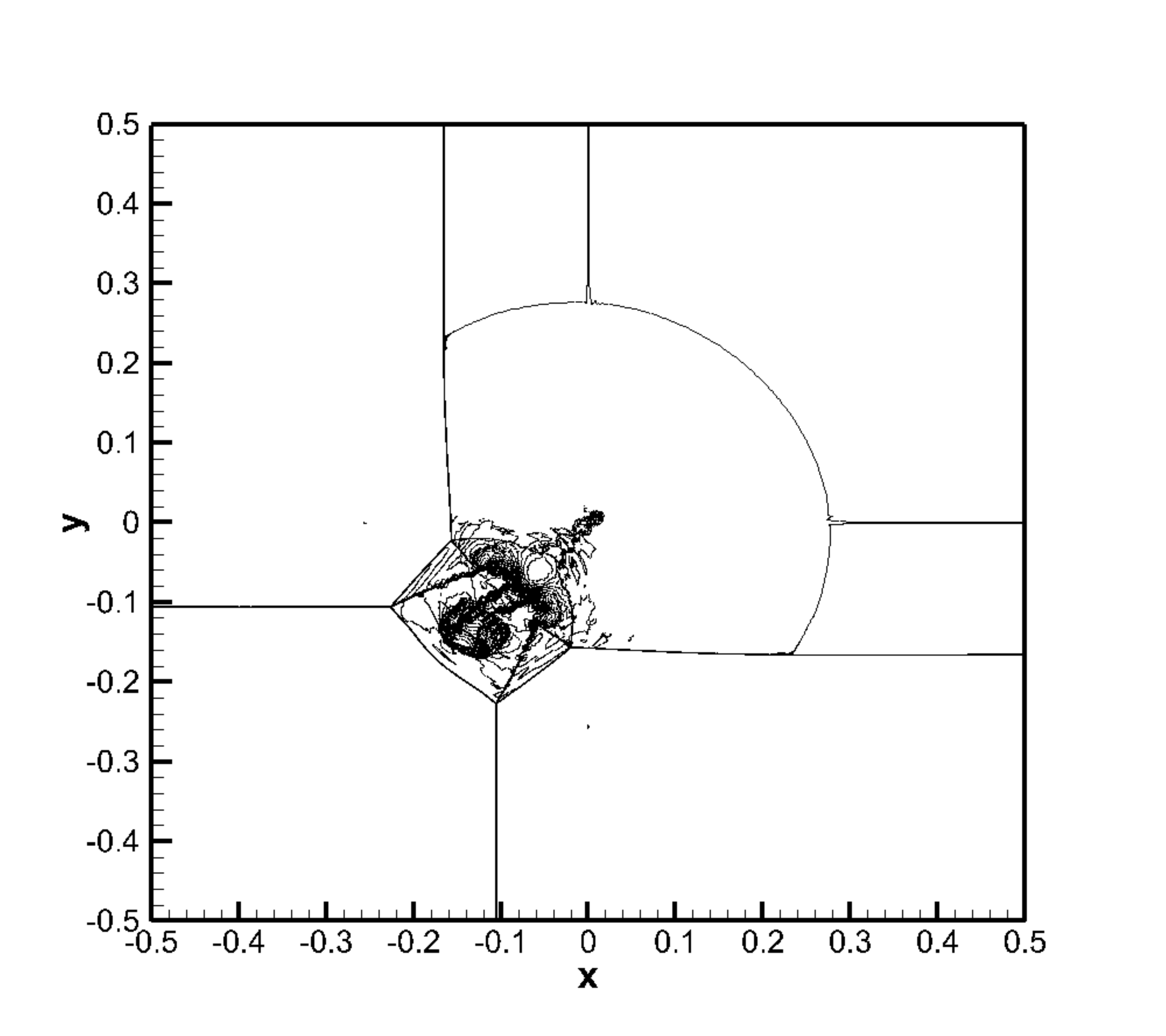}     &  
\includegraphics[width=0.35\textwidth]{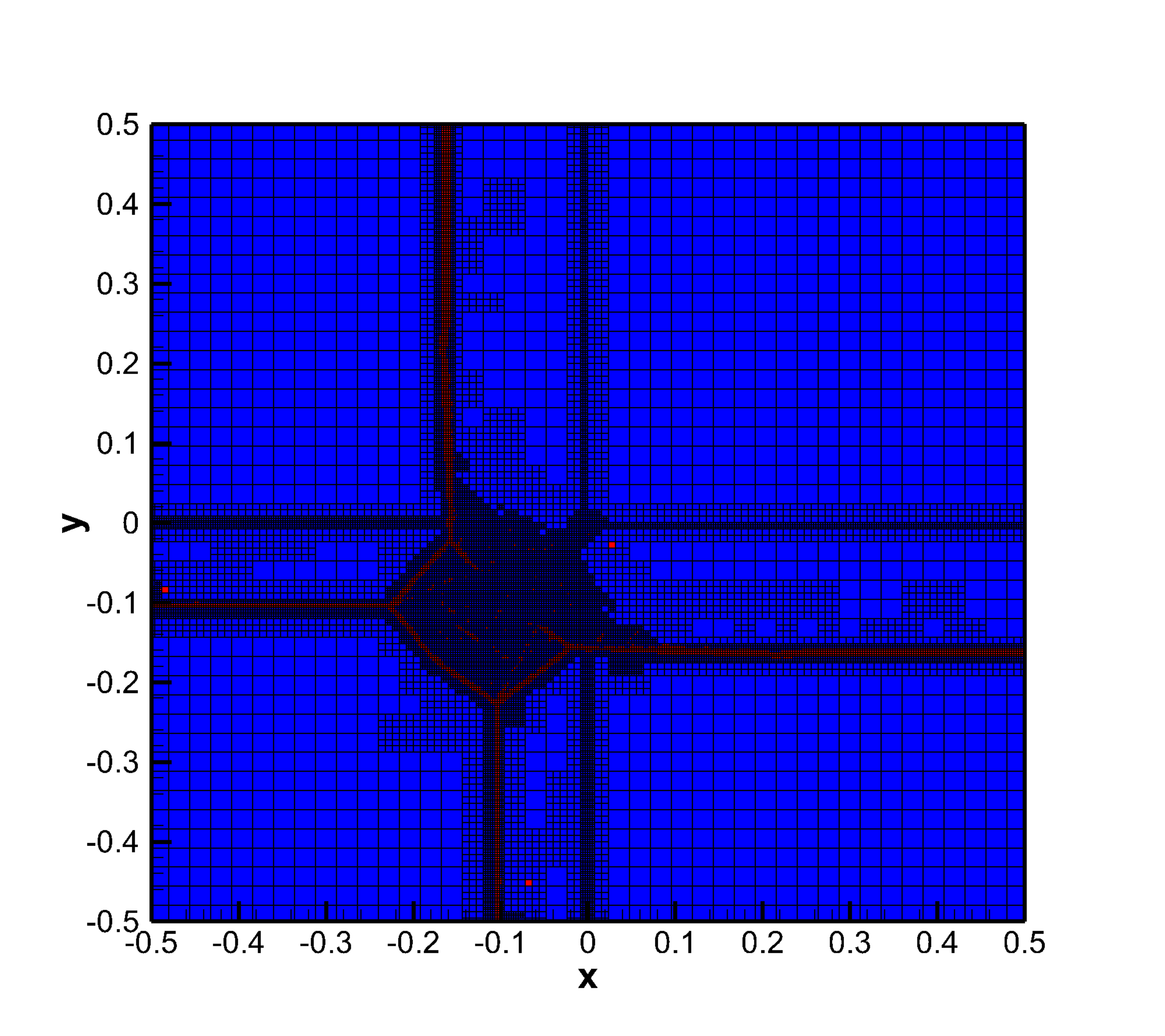} \\
\includegraphics[width=0.35\textwidth]{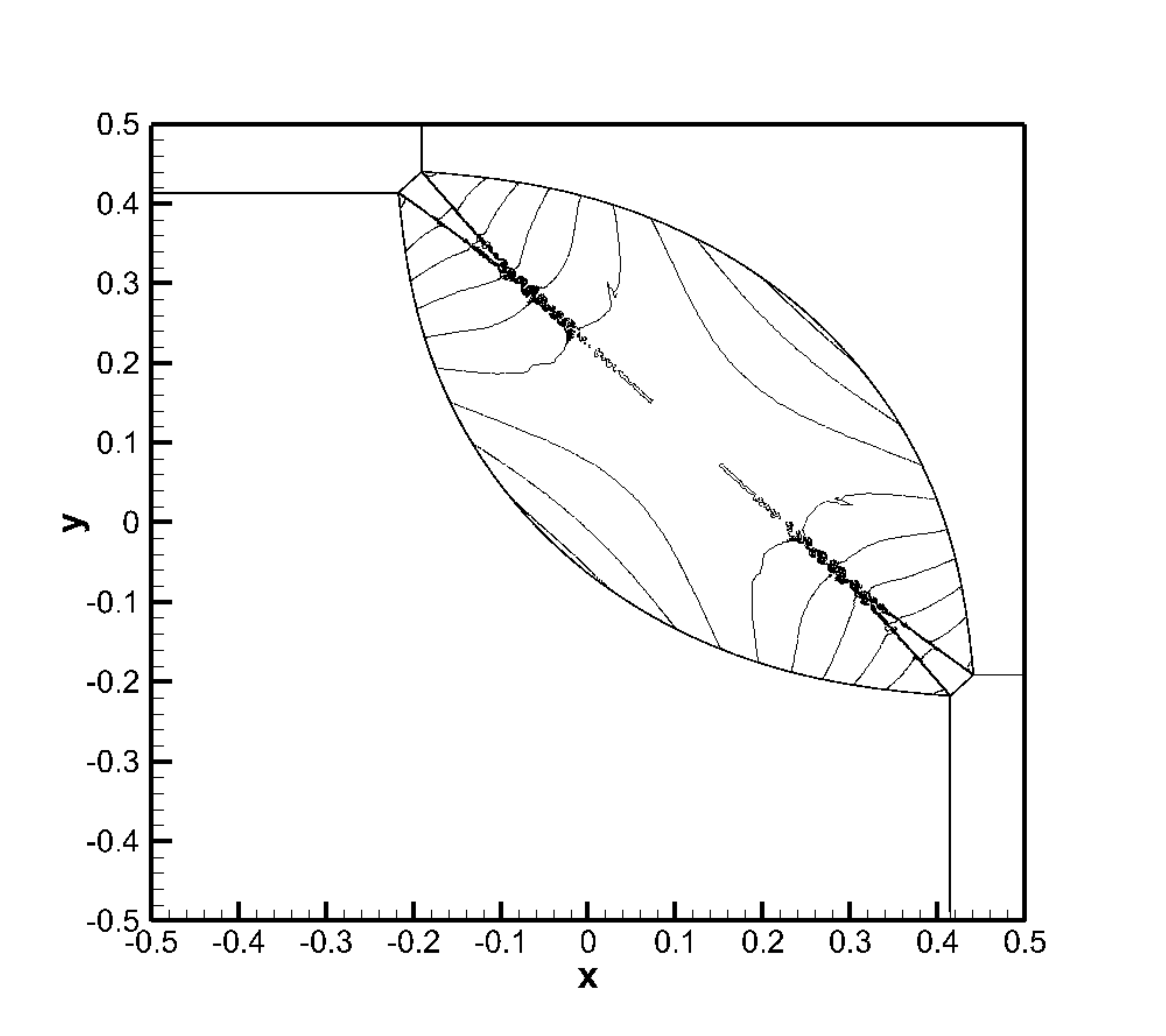}     &  
\includegraphics[width=0.35\textwidth]{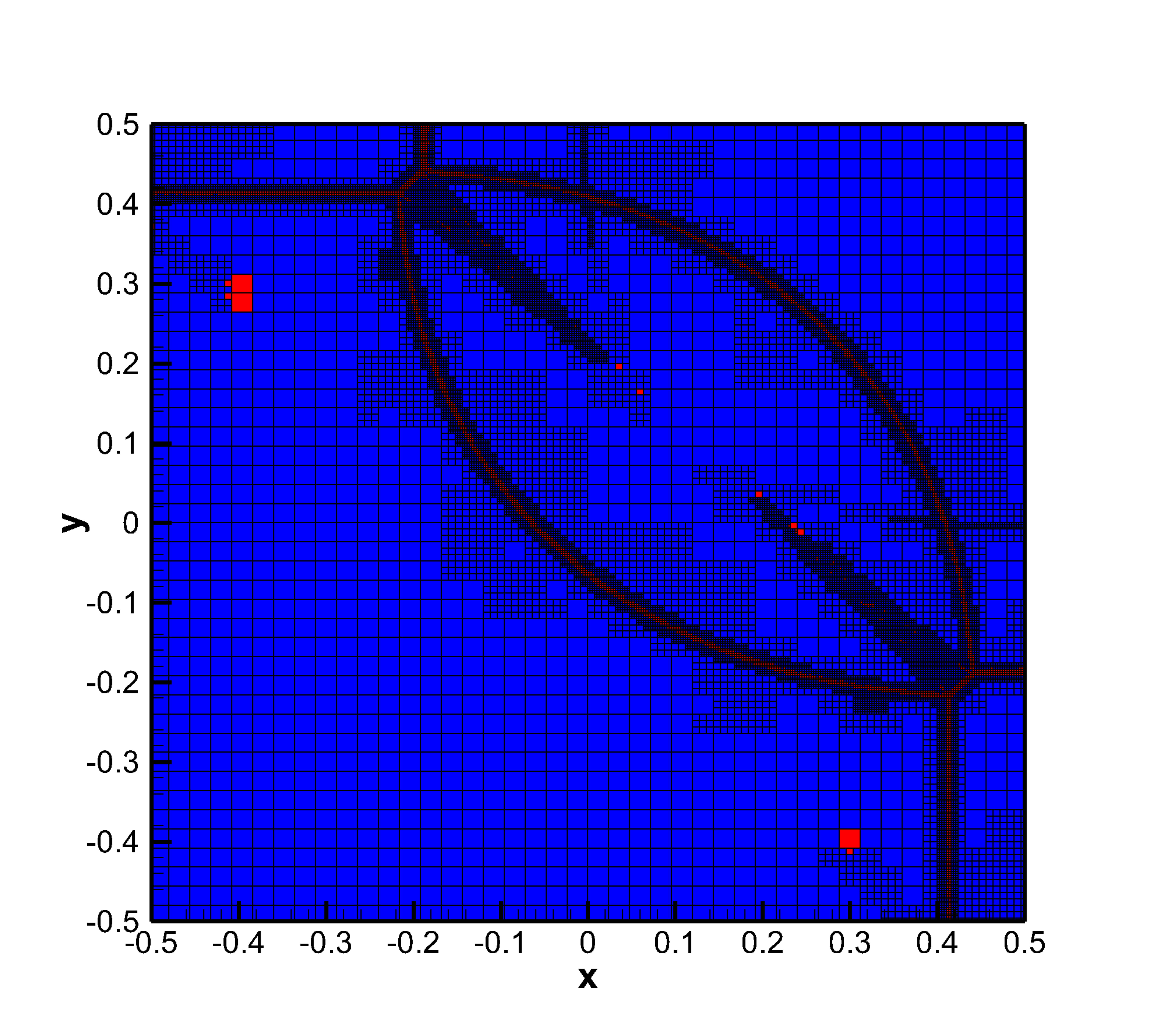} \\
\includegraphics[width=0.35\textwidth]{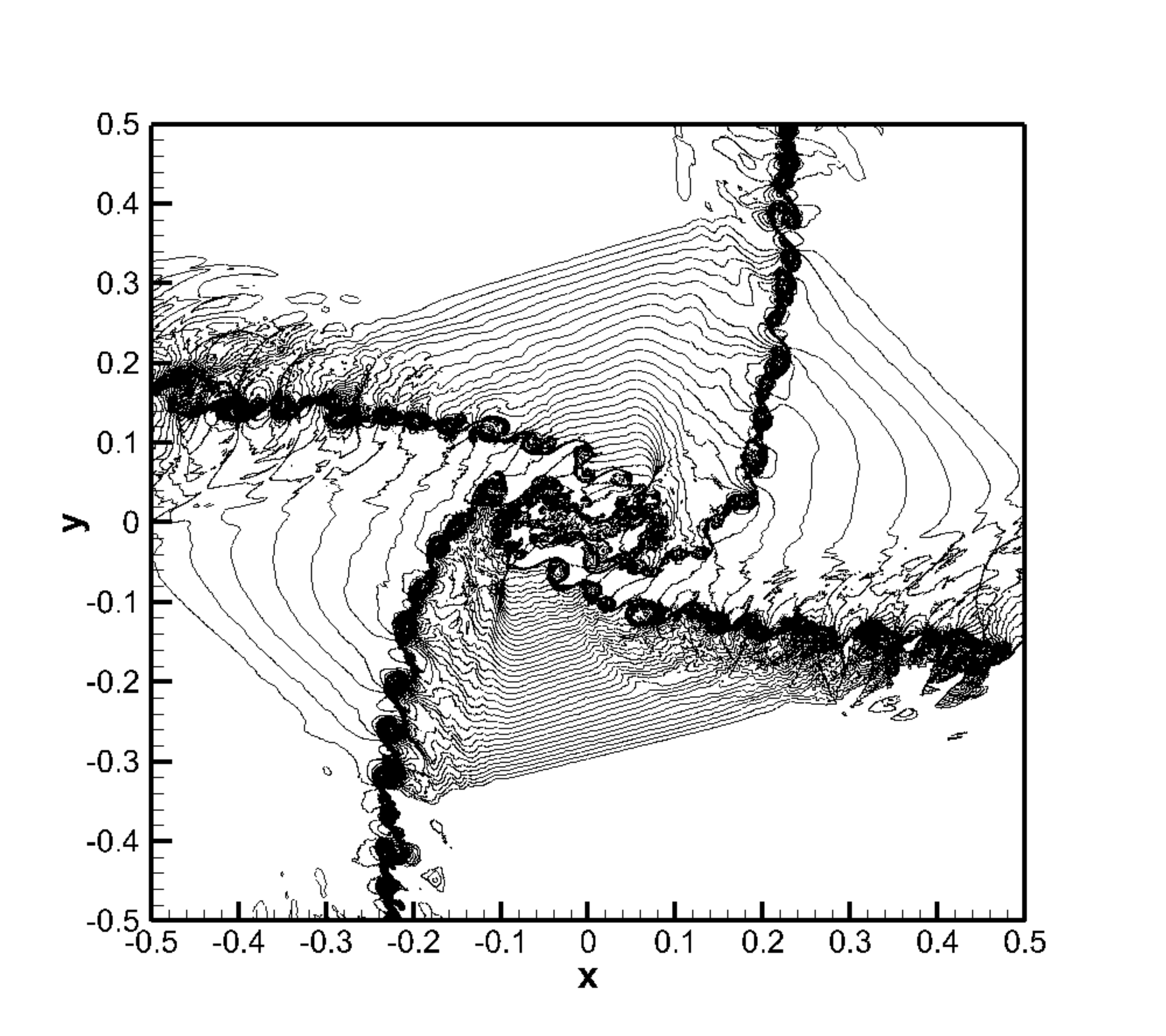}     &  
\includegraphics[width=0.35\textwidth]{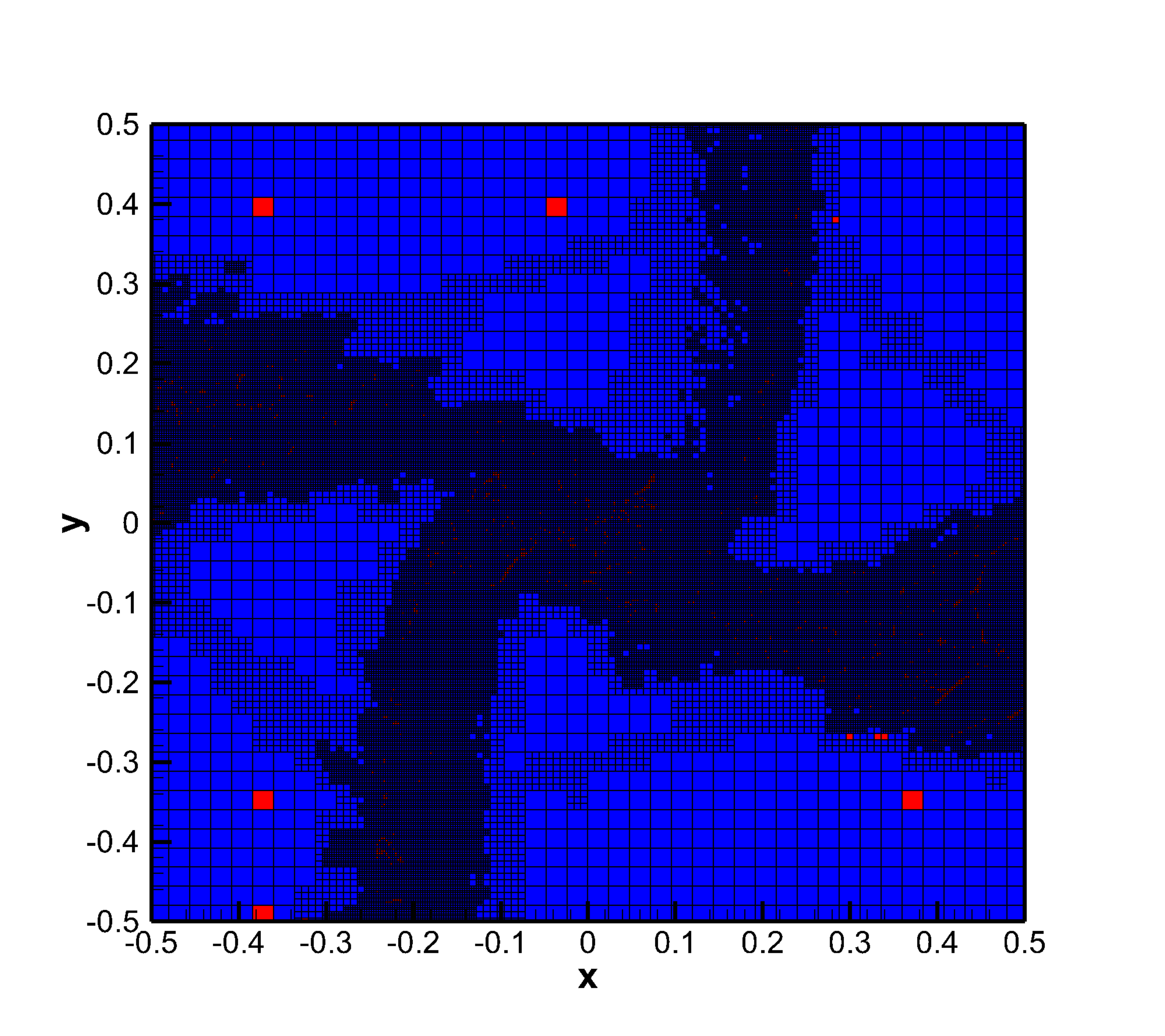} \\
\includegraphics[width=0.35\textwidth]{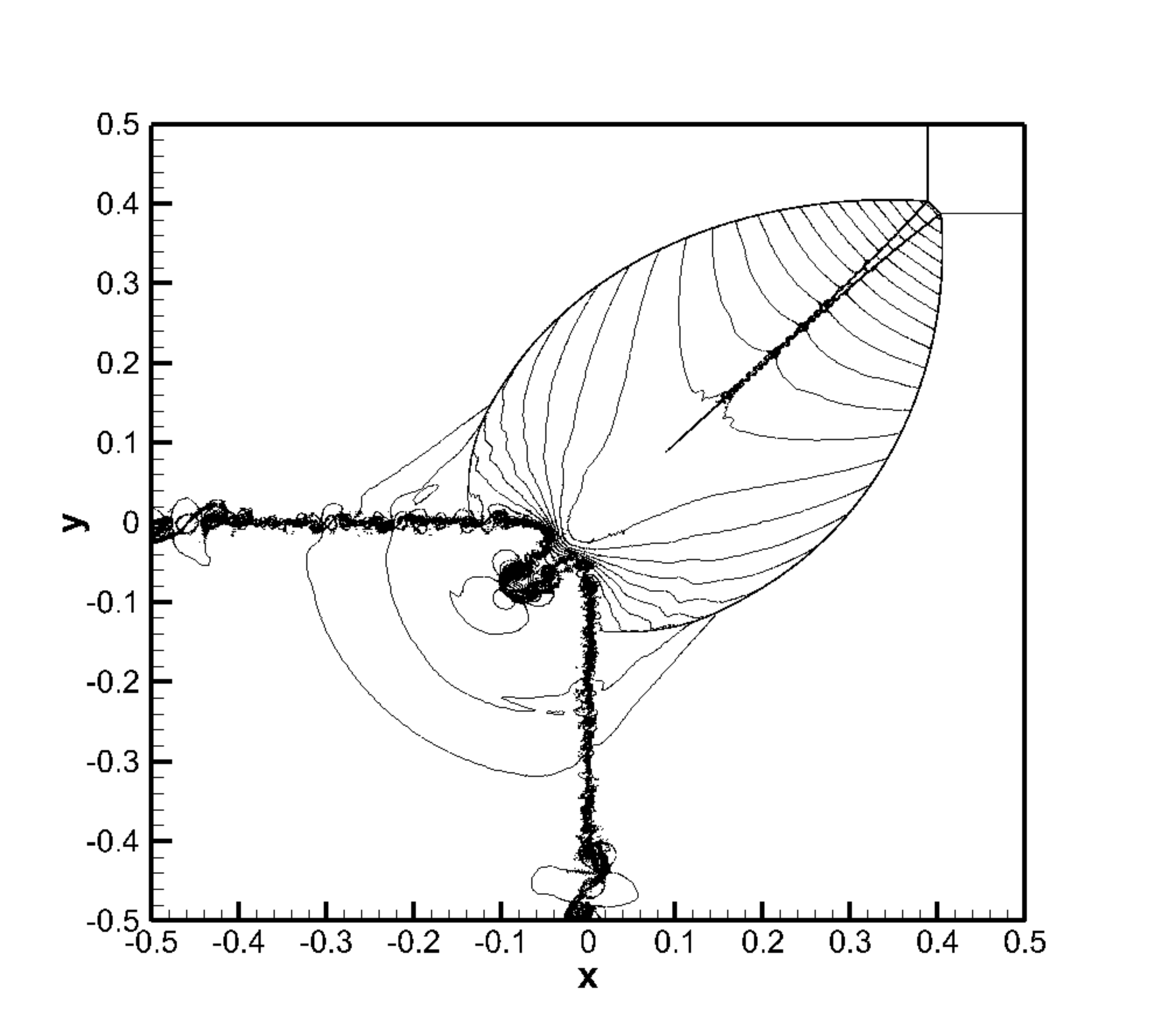}     &  
\includegraphics[width=0.35\textwidth]{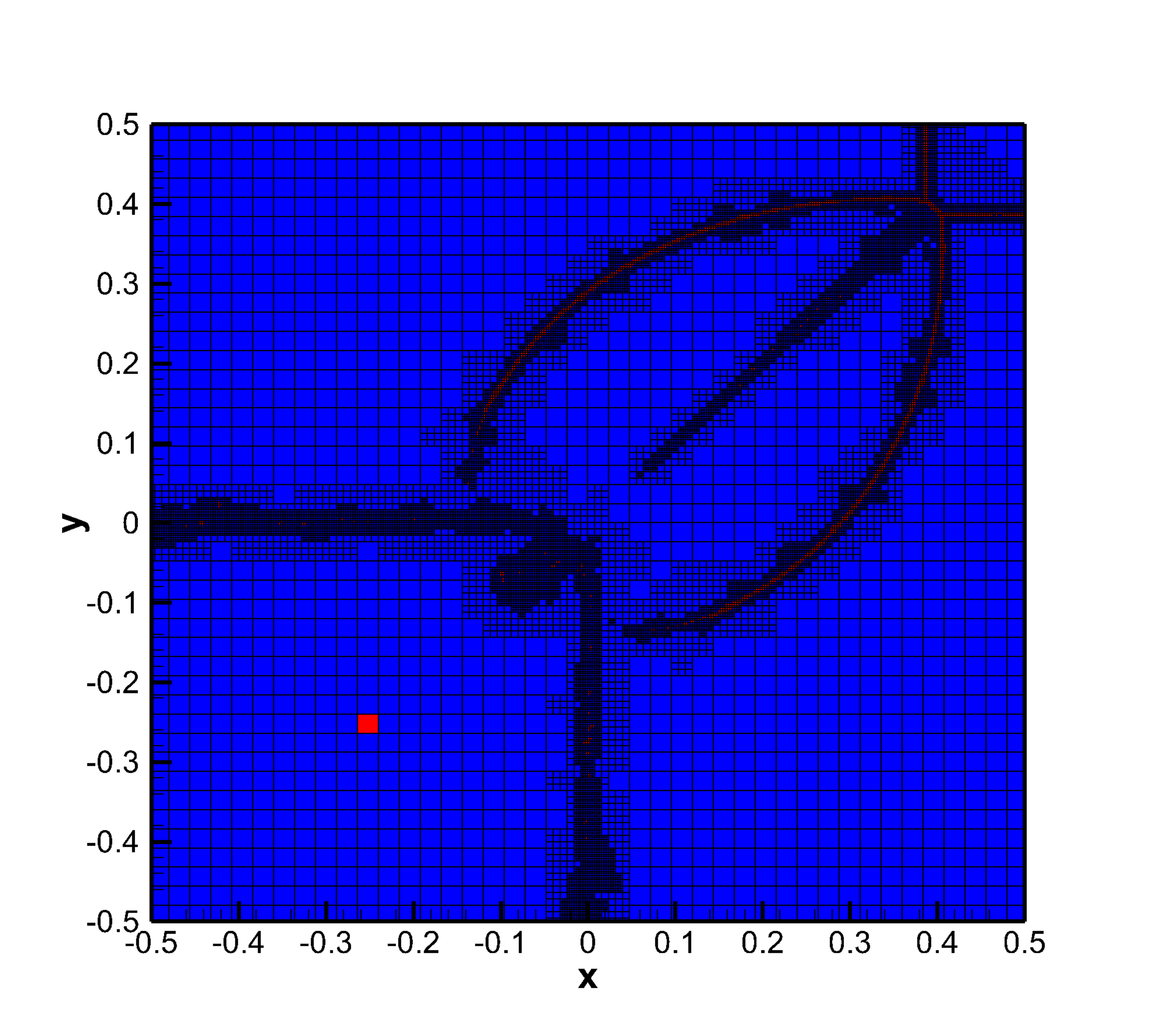}  
\end{tabular} 
\caption{ Two-dimensional Riemann problems solved with the 
 AMR-ADER-DG-$\mathbb{P}_5$ method with sub-cell limiter on an initial uniform grid with $50 \times 50$  cells.
Two levels of refinement have been adopted, with refinement factor $\mathfrak{r}=3$.
Left panels: isolines of the density. Right panels:
AMR grid (black),  limited cells (red) and unlimited cells (blue). 
}
\label{fig.rp2d}
\end{center}
\end{figure}

The two dimensional Riemann problems first proposed in \cite{kurganovtadmor} have become a classic benchmark for any numerical scheme solving the
Euler equations.
The  initial conditions are represented by constant states in each of the four quadrants of
the computational domain $\Omega = [-0.5;0.5] \times [-0.5;0.5]$, namely
\begin{equation}
 \mathbf{u}(x,y,0) = \left\{ \begin{array}{ccc} 
 \mathbf{u}_1 & \textnormal{ if } & x > 0 \wedge y > 0,    \\ 
 \mathbf{u}_2 & \textnormal{ if } & x \leq 0 \wedge y > 0, \\ 
 \mathbf{u}_3 & \textnormal{ if } & x \leq 0 \wedge y \leq 0, \\ 
 \mathbf{u}_4 & \textnormal{ if } & x > 0 \wedge y \leq 0\,. 
  \end{array} \right.
\end{equation} 
The data of the four configurations that we have considered are reported in Table \ref{tab.rp2d.ic}. We emphasize that the adiabatic index
is $\gamma=1.4$ in all cases.
%
\begin{table}[!t] 
  \numerikNine  
\begin{center} 
\begin{tabular}{|c|c||cccc|cccc|c|} 
\hline
Problem    && $\rho$ & $u$ & $v$  & $p$ & $\rho$ & $u$ & $v$  & $p$ & \multirow{2}{*}{$t_{\text{final}}$} \\ 
\cline{3-10}
    && \multicolumn{4}{c|}{$x \leq 0$} & \multicolumn{4}{c|}{$x>0$} & \\
\hline
\hline
\multirow{1}{*}{\rotatebox{0}{\textbf{RP1}} } 
                                 &$y > 0$    & 0.5323 & 1.206   & 0.0     & 0.3   & 1.5    & 0.0 &  0.0    & 1.5 &  \multirow{2}{*}{0.25} \\ 
(Case 3 in KT)&$y \leq 0$ & 0.138  & 1.206   & 1.206   & 0.029 & 0.5323 & 0.0 &  1.206  & 0.3 & \\ 
\hline
\hline
\multirow{1}{*}{\rotatebox{0}{\textbf{RP2}}}
                                 &$y > 0$    & 0.5065 &  0.8939 & 0.0     & 0.35 & 1.1    & 0.0 &  0.0    & 1.1  &\multirow{2}{*}{0.25}\\ 
(Case 4 in KT)&$y \leq 0$ & 1.1    &  0.8939 & 0.8939  & 1.1  & 0.5065 & 0.0 &  0.8939 & 0.35 &\\ 
\hline
\hline
\multirow{1}{*}{\rotatebox{0}{\textbf{RP3}}}
                                 &$y > 0$    & 2.0    &  0.75  & 0.5   & 1.0  & 1.0  &  0.75 &  -0.5  & 1.0  &\multirow{2}{*}{0.30}\\ 
(Case 6 in KT)&$y \leq 0$ & 1.0    & -0.75  & 0.5   & 1.0  & 3.0  & -0.75 &  -0.5  & 1.0  &\\ 
\hline
\hline
\multirow{1}{*}{\rotatebox{0}{\textbf{RP4}}}
                                 &$y > 0$    & 1.0    &  0.7276 & 0.0     & 1.0  & 0.5313 & 0.0 &  0.0    & 0.4  &\multirow{2}{*}{0.25}\\ 
(Case 12 in KT)&$y \leq 0$ & 0.8    &  0.0    & 0.0     & 1.0  & 1.0    & 0.0 &  0.7276 & 1.0  &\\ 
\hline 
\end{tabular} 
\caption{\label{tab.rp2d.ic} Initial conditions for the two--dimensional Riemann problems. The "Case No. in KT" refers to the classification of \cite{kurganovtadmor}.} 
\end{center}
\end{table} 
%
The simulations have been performed over a level zero grid of $50 \times 50$ elements, adopting $\ell_{\rm max}=2$ and
$\mathfrak{r}=3$. On the other hand, the numerical scheme is 
the ADER-DG-$\mathbb{P}_5$, with the Rusanov Riemann solver and reconstruction in characteristic variables.
Fig.~\ref{fig.rp2d} shows the result of the simulations at the final time $t_{\text{final}}$ for each model. The left panels report the isolines of the density, while the right panels show, as usual, 
the AMR mesh and the cells updated through the sub-cell limiter, which have been highlighted in red. Due to the unprecedented high order of accuracy adopted, which reduces drastically the numerical
dissipation of the numerical scheme, several small-scale features appear in the solution, typically attributed to the Kelvin--Helmholtz instability but
not visible in the original versions shown by \cite{kurganovtadmor}. A similar effect was already 
noticed by \cite{Dumbser2014} for the test $\textbf{RP3}$, even in the absence of AMR. However, when adaptive mesh refinement is activated, the effects of the 
Kelvin--Helmholtz instability emerge
clearly also in model $\textbf{RP2}$ (along the diagonal of the cocoon structure),  and in model $\textbf{RP4}$ (along the boundary of the bottom-left quadrant). Moreover, we emphasize that the use of AMR makes the sub-cell limiter operate only along strong discontinuities, which are resolved within very few cells at the maximum level of refinement. 

\subsubsection{Cylindrical and spherical explosion problem} 

In multiple space dimensions, a conceptually simple but interesting extension of the one-dimensional Riemann problem is represented by the  
cylindrical and by the spherical explosion problem, both of them described with great detail in \cite{titarevtoro} and \cite{toro-book}. These two tests are indeed very relevant, 
since they involve the propagation of a shock wave that is not aligned with the coordinates, and they can therefore be used to check the ability of the numerical scheme 
in preserving the physical symmetries of the problem. 
As initial conditions, we assume the flow variables to be constant for $r\leq R$ and for $r\geq R$, namely
\begin{equation}\label{RP3D}
  \big(\rho, u, v, w, p\big) = \left\{\begin{array}{ll}
\big(1, 0, 0, 0,   1\big) & \;\textrm{for}\quad r \leq R\,, \\
 \noalign{\medskip}
 \big(0.125,   0, 0, 0,   0.1\big) & \;\textrm{for}\quad r>R    \,,  
 \end{array}\right.
\end{equation} 
where $r = \sqrt{\mathbf{x}^2}$ is the radial coordinate, $\mathbf{x}$ is the vector of spatial coordinates, while $R=0.5$ denotes the radius of the initial discontinuity.  
The computational domain is $\Omega=[-1; 1]^d$, whereas the adiabatic index of the ideal-gas equation of state has been set to $\gamma=1.4$. 
As suggested by \cite{toro-book}, a reference solution can be computed after solving an equivalent one dimensional problem in the radial direction $r$, in which the additional 
geometric terms arising from the choice of curvilinear coordinates can be moved to the right hand side of the governing PDEs as source terms. 

\begin{figure}[!htbp]
  \begin{center}
    \includegraphics[width=0.85\textwidth]{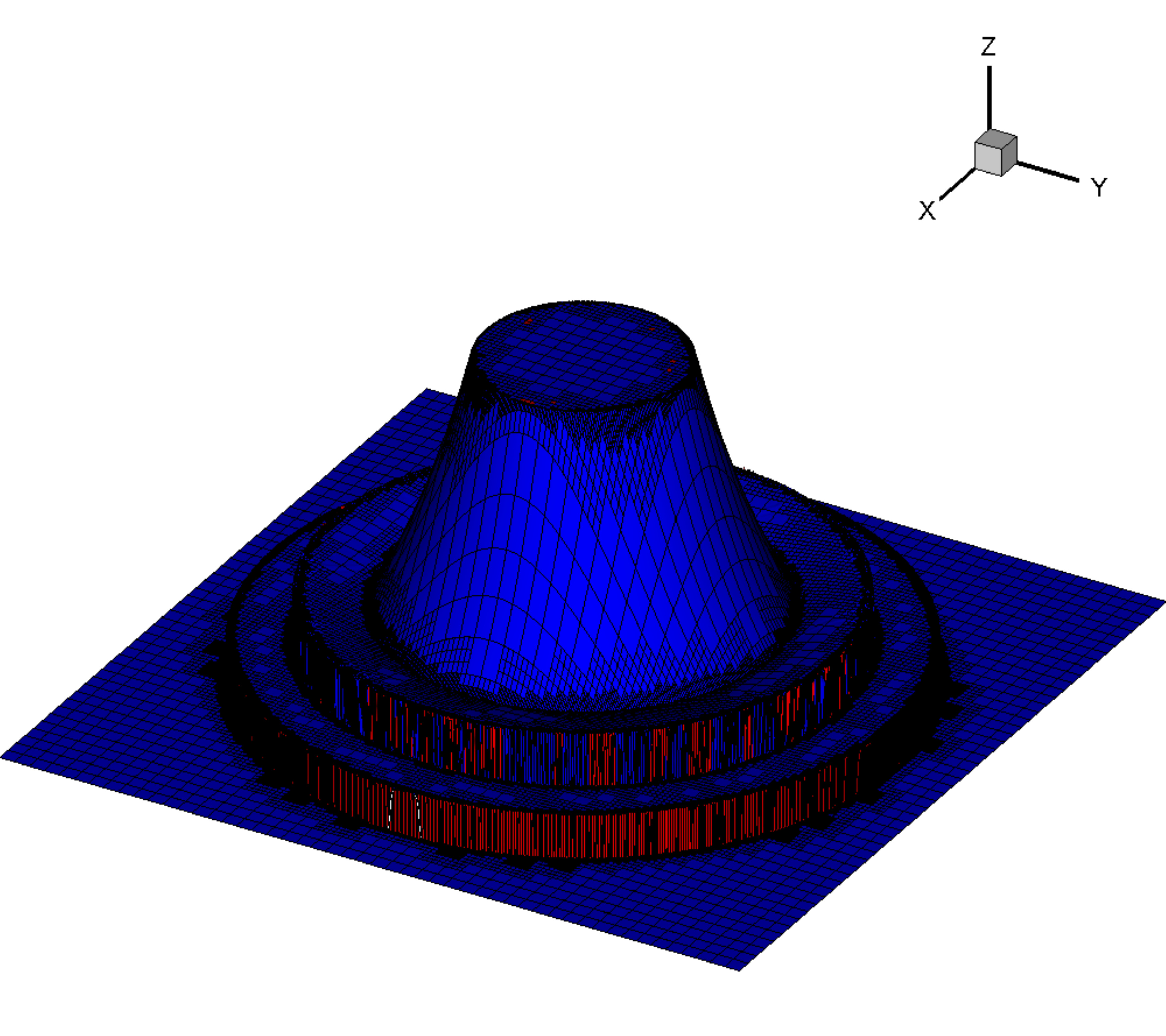}   
    \caption{ \label{fig:EP2D-limiter}
      Three-dimensional view of the density variable and the AMR grid for the two dimensional explosion problem at $t_{\text{final}}=0.20$. 
			Limited cells (red) updated with the sub-cell ADER-WENO3 finite volume scheme and unlimited cells (blue) with
			the ADER-DG-$\mathbb{P}_9$ scheme. 
			The level zero AMR grid uses $50\times50$ elements. } 
  \end{center}
\end{figure}

\begin{figure}[!htbp] 
  \begin{center} 
    \begin{tabular}{cc}
      \includegraphics[width=0.45\textwidth]{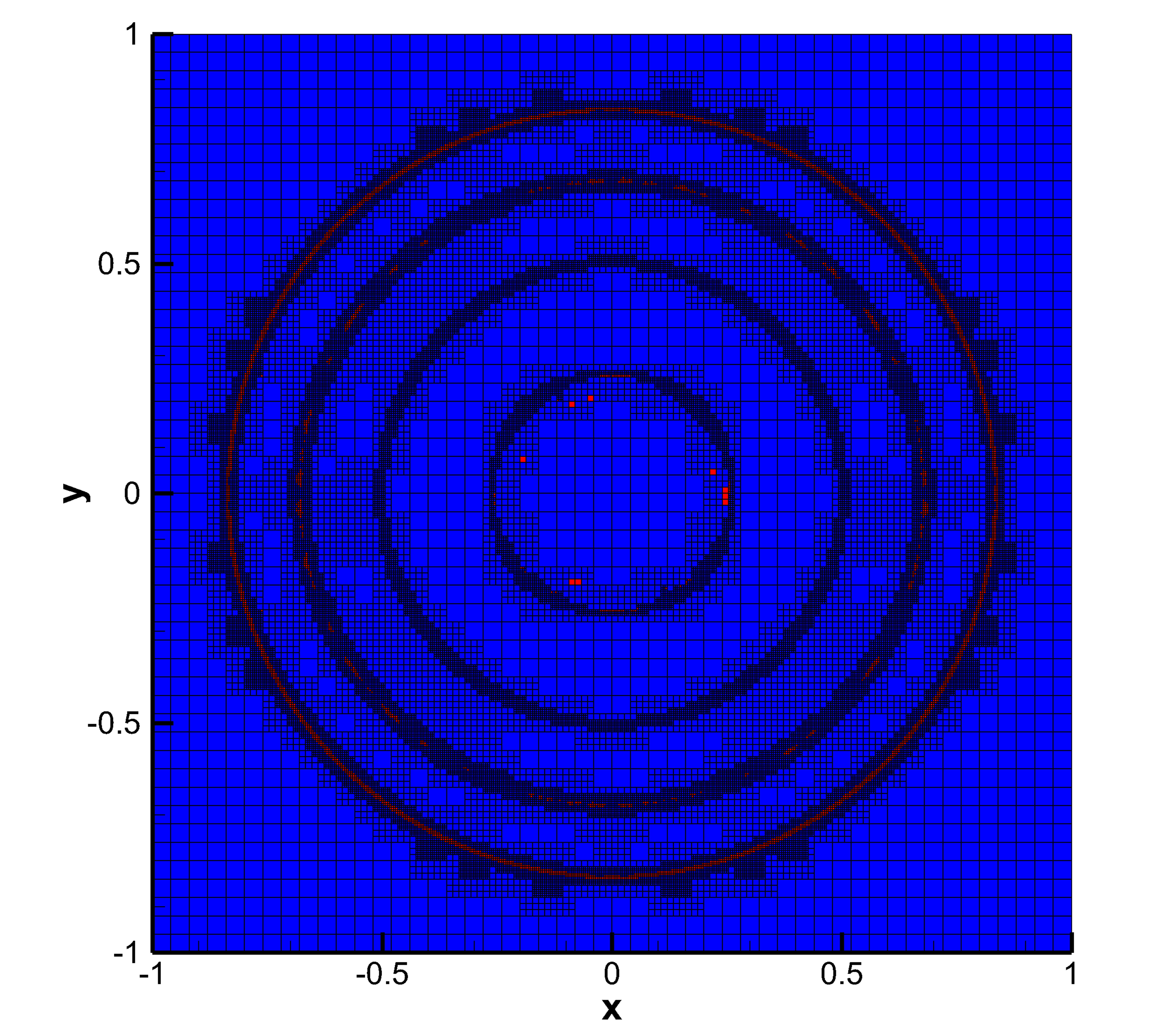} & 
      \includegraphics[width=0.45\textwidth]{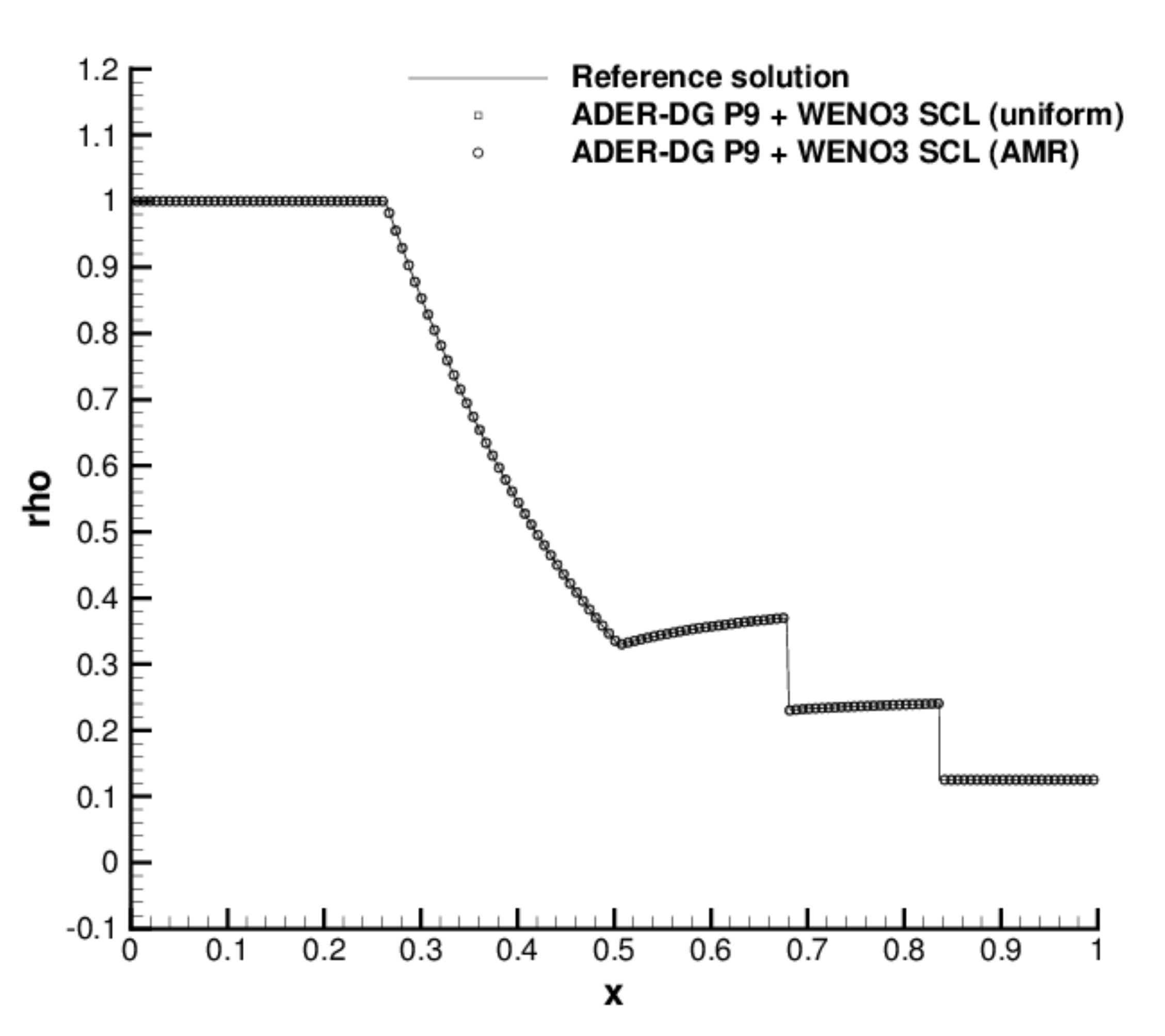}  \\
      \includegraphics[width=0.45\textwidth]{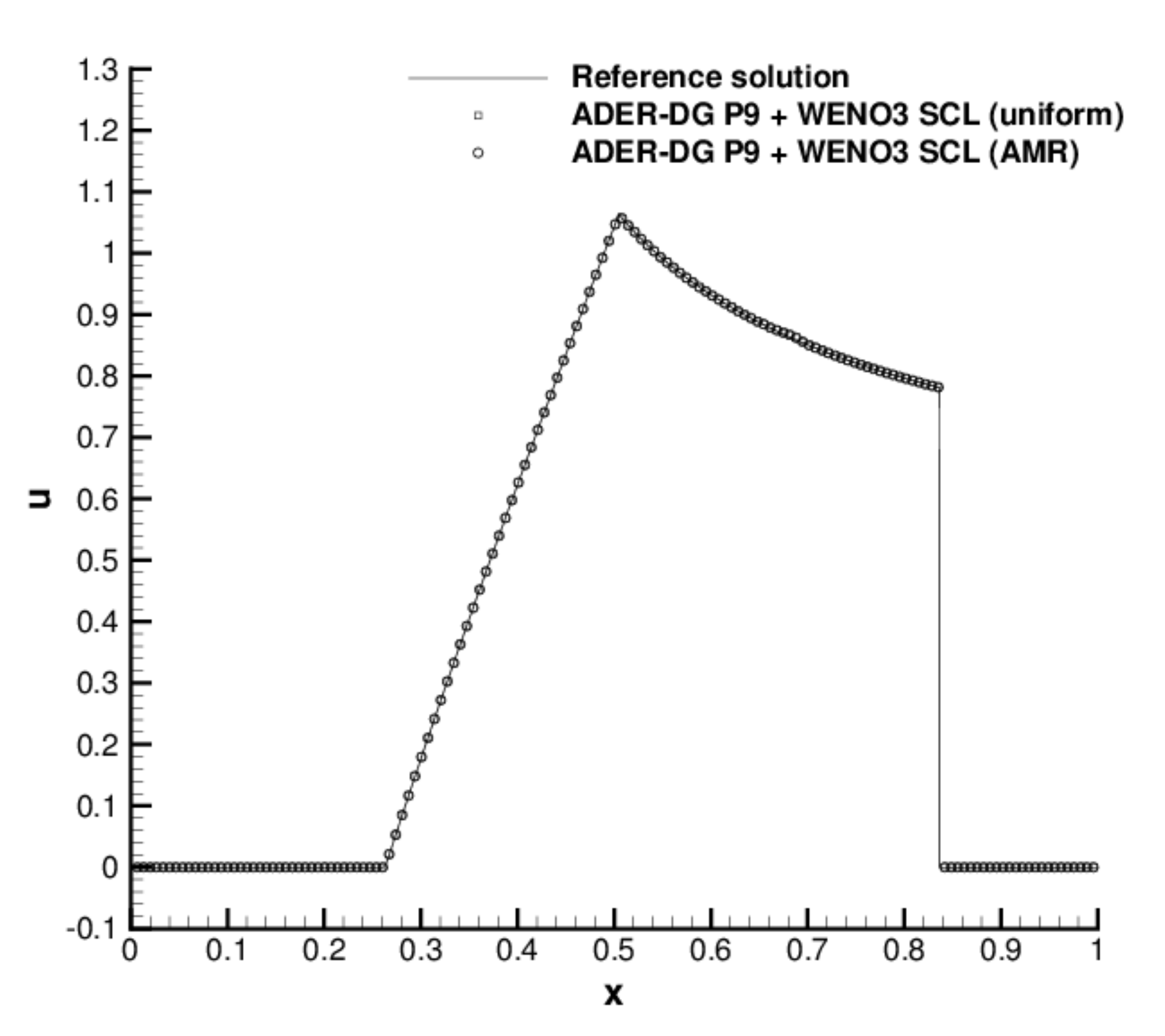}    &
      \includegraphics[width=0.45\textwidth]{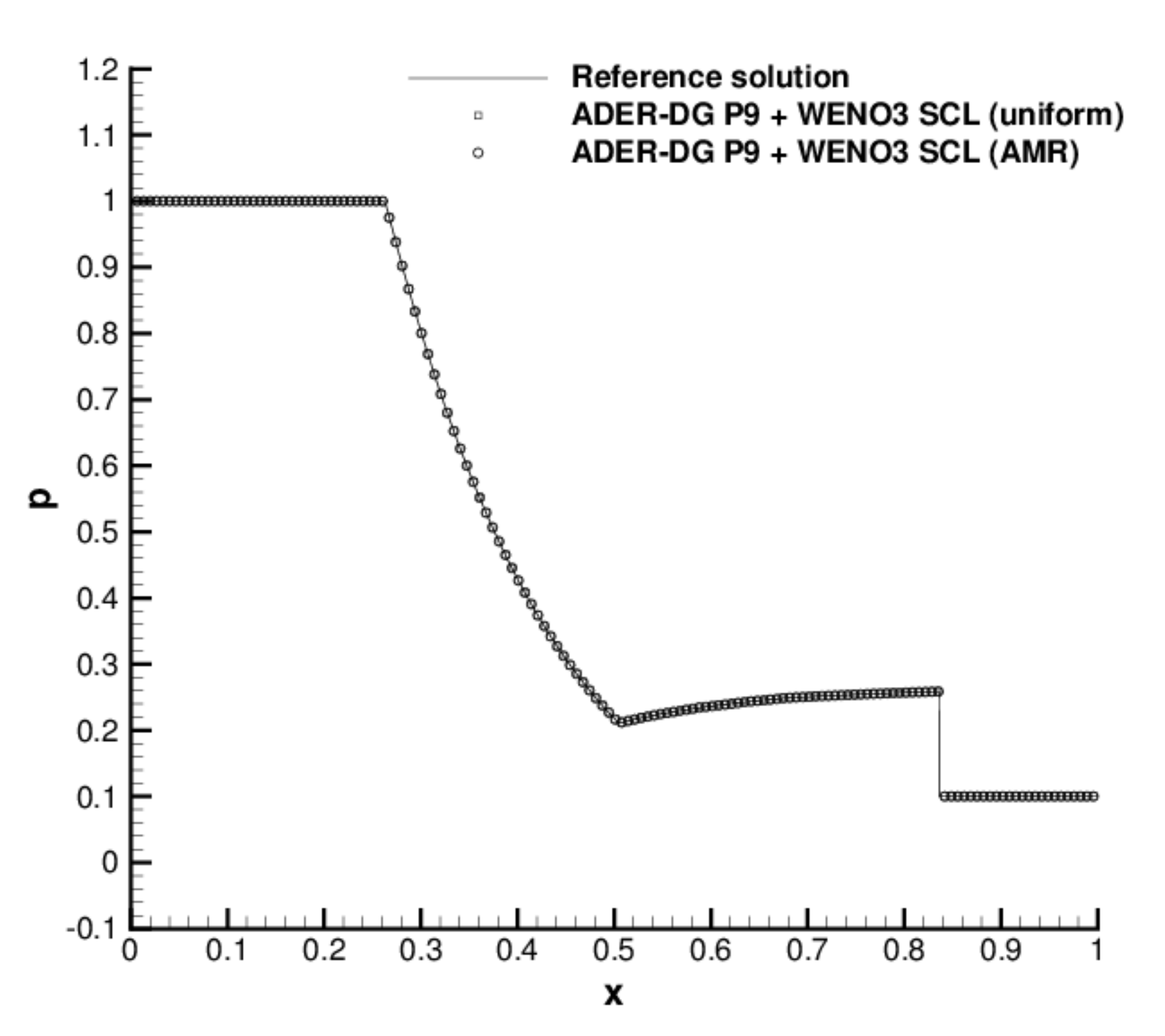}    
    \end{tabular} 
    \caption{ \label{fig:EP2D-1D-cuts}
      Two-dimensional explosion problem. 2D view of the AMR grid together with limited and unlimited cells (top left). 
			One dimensional cuts of the numerical solution for density $\rho$ (top right), velocity $u$ (bottom left) and 
			fluid pressure $p$ (bottom right) on 150 equidistant sample points along the positive $x-$axis obtained at 
			$t_{\text{final}}=0.20$ with the space-time adaptive ADER-DG-$\mathbb{P}_9$ scheme, supplemented with \aposteriori ADER-WENO3 
			sub-cell limiter. For comparison, the solution computed on a uniform fine mesh corresponding to the finest AMR grid 
			level and the 1D reference solution are also reported.} 
  \end{center}
\end{figure}

We have solved the two-dimensional, cylindrical, explosion problem with the ADER-DG-$\mathbb{P}_9$ scheme in combination with our usual \aposteriori sub-cell WENO finite volume limiter, 
the Osher-type flux of \cite{OsherUniversal} and the reconstruction in characteristic variables. 
On the level zero grid, the mesh consists of $50 \times 50$ elements, which are then refined using a refinement factor of $\mathfrak{r}=3$ and $\ell_{\max}=2$. 
This leads to an equivalent resolution on a uniform fine grid of $450 \times 450 = 202,500$ elements. Considering that each P9 element uses 10 degrees of freedom 
per space dimension, this corresponds to a total resolution of $20,250,000$ spatial degrees of freedom on a uniform fine grid. Fig.~\ref{fig:EP2D-limiter} shows a 3D 
plot of the density distribution obtained for the cylindrical explosion case, as well as the AMR grid configuration at the final time $t=0.20$.
Moreover, a 2D view of the 
AMR grid together with 1D cuts through the numerical solution on 150 equidistant sample points along the $x$-axis are depicted in Fig. \ref{fig:EP2D-1D-cuts}. For 
comparison, Fig. \ref{fig:EP2D-1D-cuts} also contains the 1D reference solution as well as the numerical solution obtained with the ADER-DG-$\mathbb{P}_9$ scheme on the uniform 
fine grid. First of all, we observe that the numerical results coincide perfectly well with the reference solution. Secondly, one can note that the uniform fine grid
solution as well as the result obtained with AMR are essentially identical.  

\begin{figure}[!htbp] 
  \begin{center} 
    \begin{tabular}{cc}
      \includegraphics[width=0.45\textwidth]{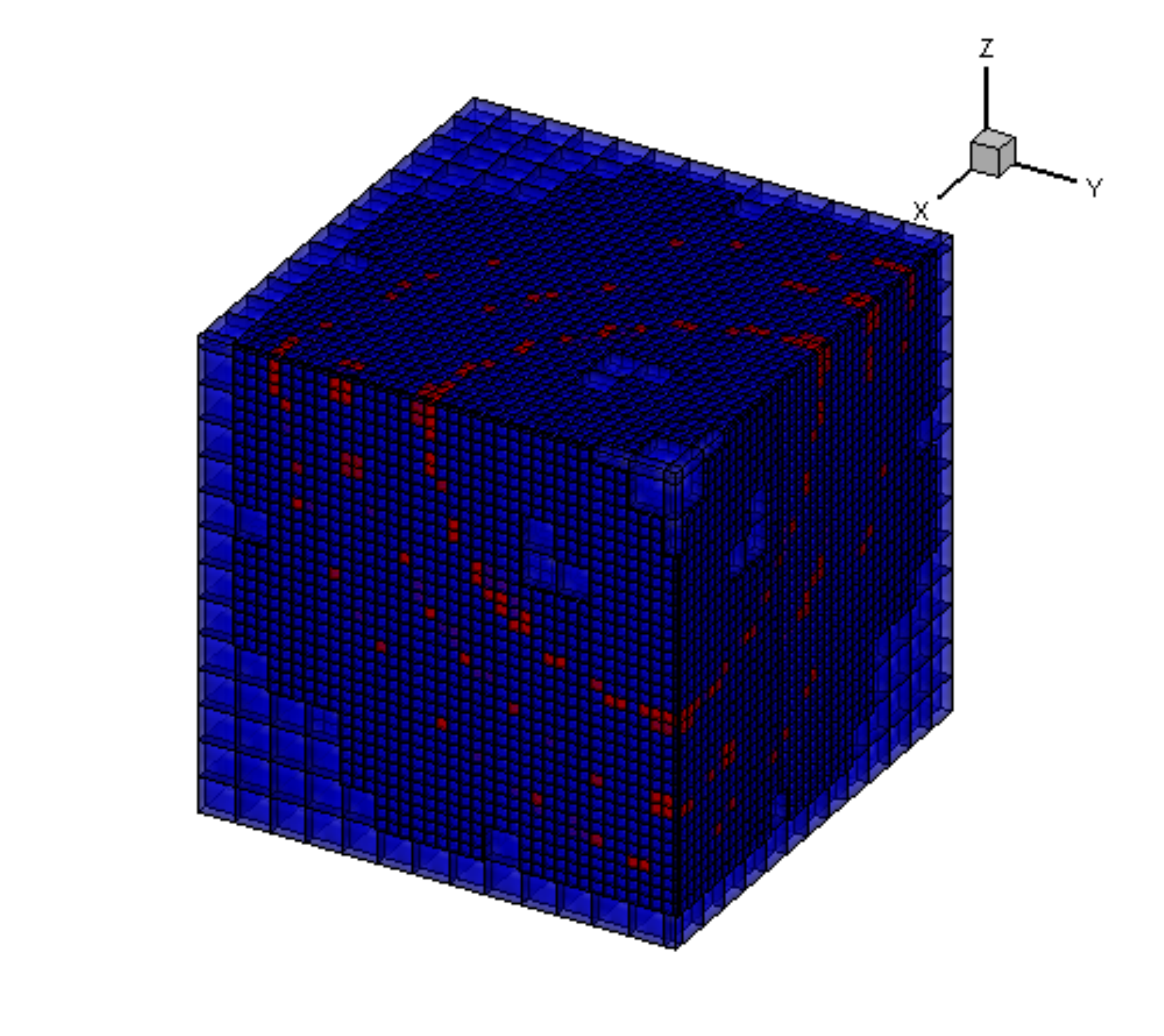} & 
      \includegraphics[width=0.45\textwidth]{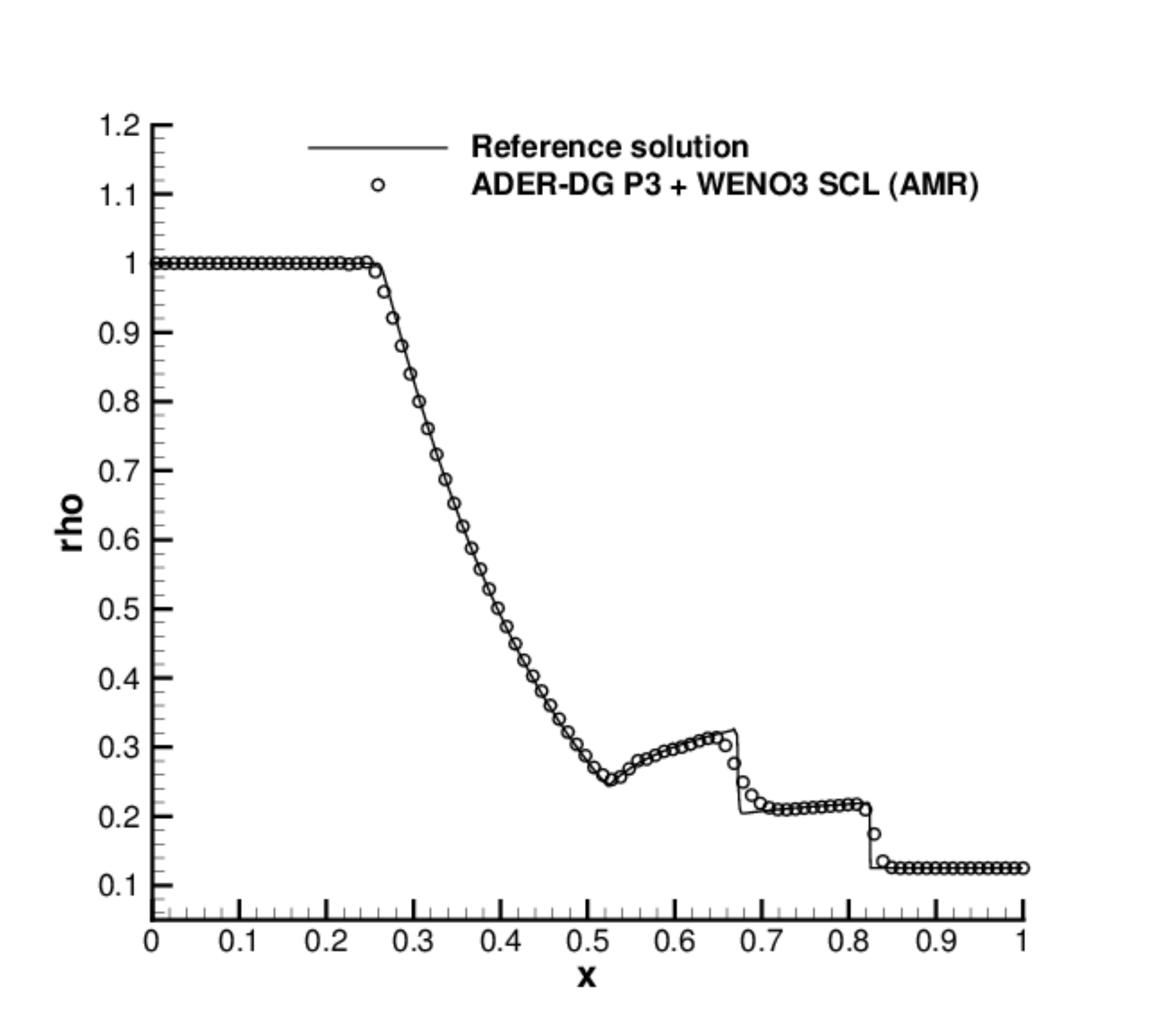}  \\
      \includegraphics[width=0.45\textwidth]{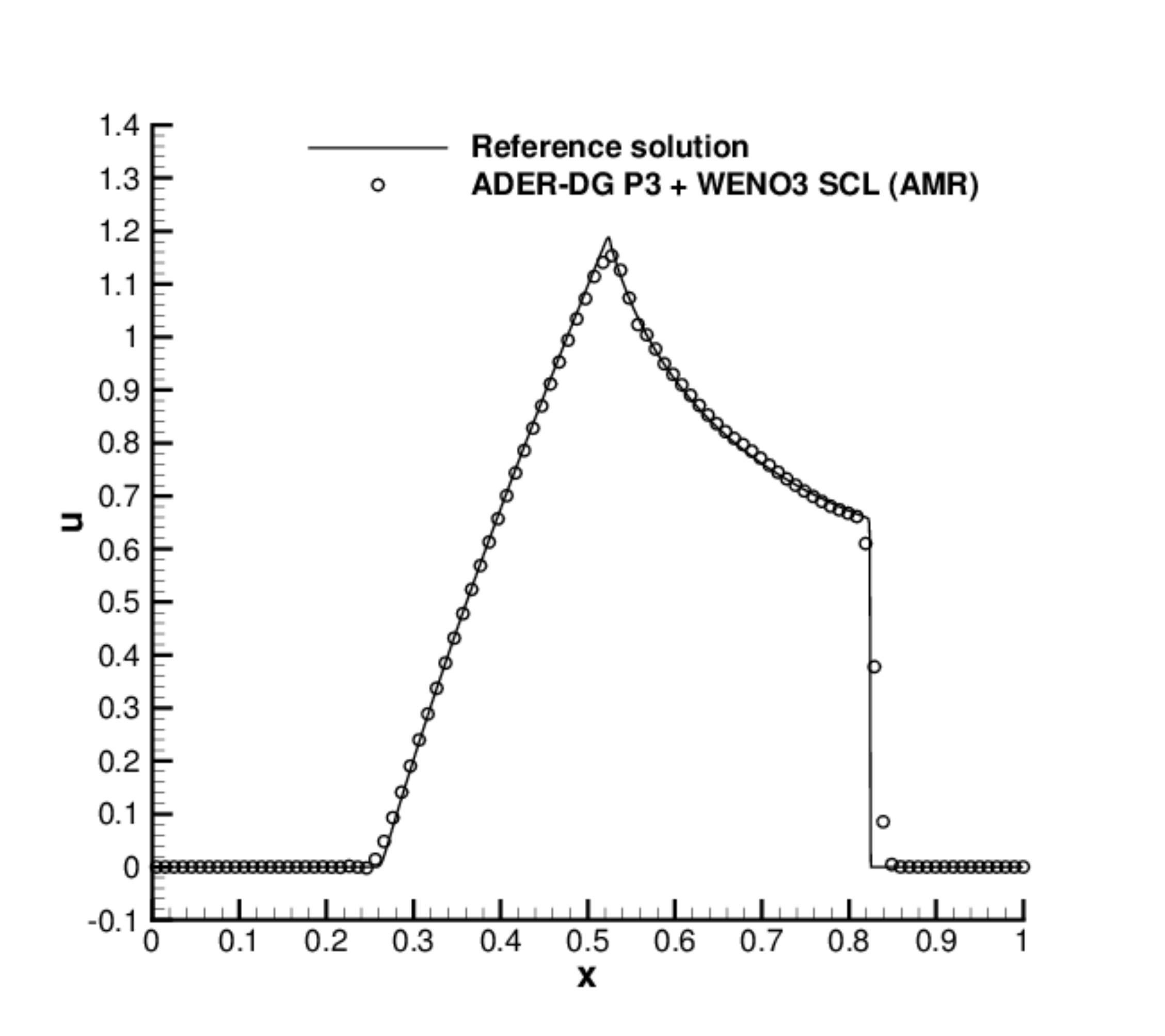}    &
      \includegraphics[width=0.45\textwidth]{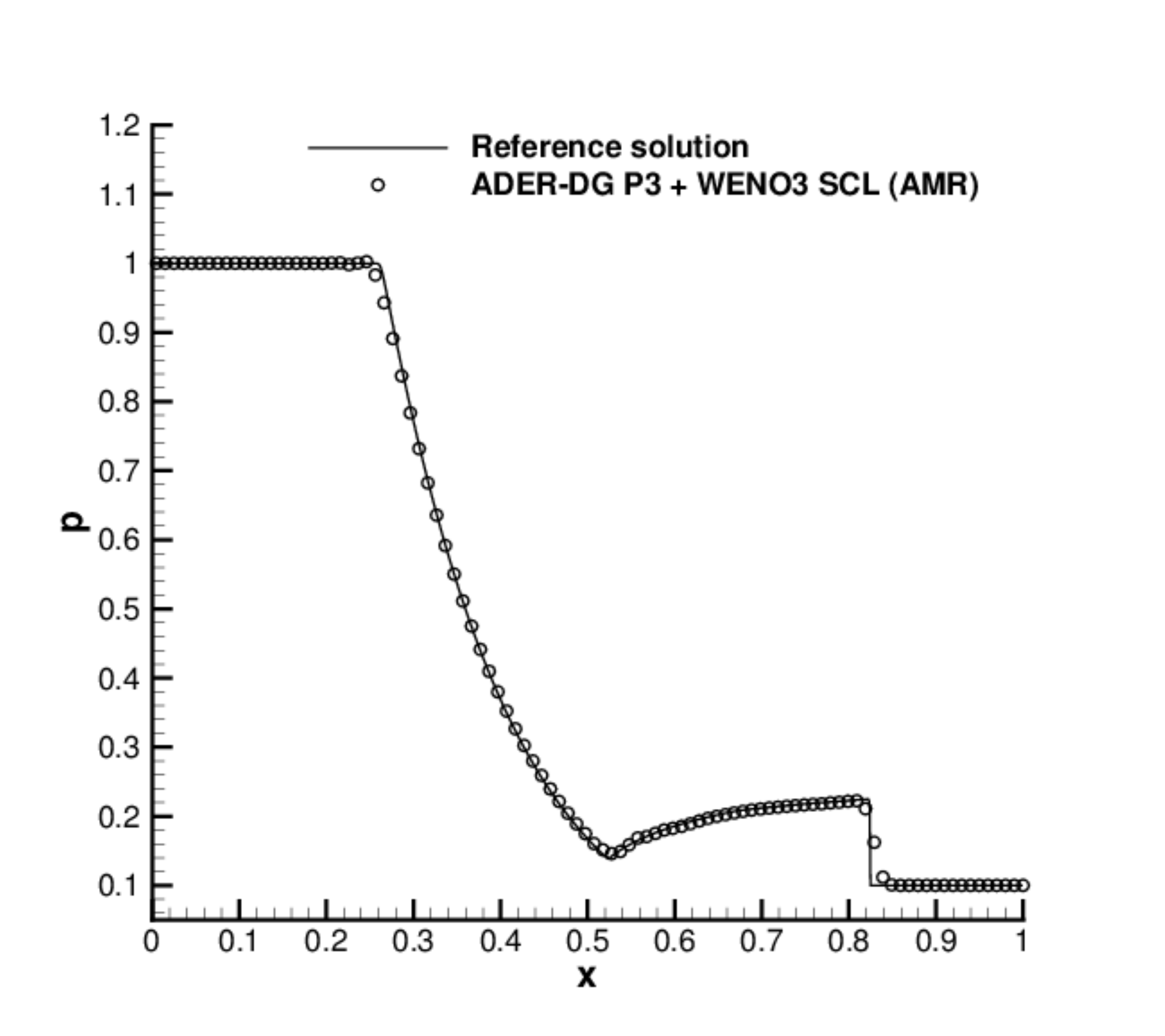}    
    \end{tabular} 
    \caption{ \label{fig:EP3D-1D-cuts}
      Three-dimensional explosion problem. 3D view of the AMR grid together with limited and unlimited cells (top left). 
			One dimensional cuts of the numerical solution for density $\rho$ (top right), velocity $u$ (bottom left) and 
			fluid pressure $p$ (bottom right) on 120 equidistant sample points along the positive $x-$axis obtained at 
			$t_{\text{final}}=0.20$ with the space-time adaptive ADER-DG-$\mathbb{P}_5$ scheme, supplemented with \aposteriori ADER-WENO3 
			sub-cell limiter. } 
  \end{center}
\end{figure}

In addition to the two dimensional case, we have also solved the spherical explosion problem in three spatial dimensions. In this case a very coarse initial mesh has been adopted, consisting of $13\times13\times13$ cells, which is subsequently refined using $\mathfrak{r}=3$ and $\ell_{\max}=2$. The problem has been solved with 
the ADER-DG-$\mathbb{P}_3$ scheme, Rusanov flux and reconstruction in characteristic variables.
The results are shown in Fig.~\ref{fig:EP3D-1D-cuts}. As it is apparent from the top-left panel of this figure, the limiter has
bee activated only at the shock front, at the contact discontinuity and at the head of the rarefaction wave. The comparison with the reference solution is also good.

\subsection{Equations of ideal magnetohydrodynamics} 
A more complex and interesting hyperbolic system is represented by the classical equations of ideal magnetohydrodynamics (MHD). 
These equations are often used to model the dynamics of an electrically conducting fluid in which the hydrodynamic
and the electromagnetic forces are comparable. Unlike the previous hyperbolic system for the classical Euler equations, an additional difficulty is that the numerical scheme must
guarantee that the magnetic field $\B$ remains \emph{locally} divergence-free, assuming that $\nabla\cdot \B=0$ at the initial time. 
While this is theoretically implied by the Maxwell equations, from the numerical point of view specific procedures must be adopted to prevent 
significant deviations from $\nabla\cdot \B=0$
due to accumulation of the numerical error. 
Over the years, several approaches have been adopted to solve this problem (see the review by \cite{Toth2000}).
In our work we have chosen the so-called {\em divergence-cleaning} procedure, which uses the hyperbolic version of the generalized 
Lagrangian multiplier (GLM) divergence cleaning method of \cite{Dedneretal}. By defining an additional auxiliary variable $\Phi$, a coupling term and a linear scalar PDE
are introduced into the MHD system in order to allow the resulting \emph{augmented} system to transport any possible divergence error (or numerical magnetic monopole) 
out of the numerical domain by itself, with an established cleaning velocity $c_h$. 
In this way, the augmented MHD system can be written in conservative form by defining 
the state vector $\u$ and the flux tensor $\F$ as
\begin{equation}
\u=\left[\begin{array}{c}
\rho \\ \rho v_j \\ E \\ B_j \\ \Phi
\end{array}\right],~~~
{F^i}=\left[\begin{array}{c}
 \rho v_i \\
 \rho v_i v_j + \left( p + \B^2 / 8\pi \right) \delta_{ij} -  B_i B_j / 4\pi   \\
 \left( E + p + \B^2 / 8\pi \right) v_i  - \left( \v \cdot \B  \right)   B_i  / 4\pi \\
 \epsilon_{jik} E_k + \Phi \delta_{ij} \\
 c^2_h B_i
\end{array}\right], \, \hspace{1cm}i,j,k=x,y,z.
\label{eq:MHDfluxes}
\end{equation}
where $\B = (B_x,B_y,B_z)$ is the magnetic field vector, $\mathbf{E} = (E_x,E_y,E_z)$ is the electric field vector, $\delta_{ij}$ is the Kronecker delta and $\epsilon_{ijk}$ is the Levi-Civita symbol. 
The equation of state is again that of an ideal gas [cf. Eq.~(\ref{eq:EOS})], while the 
total energy density is given by $E=p/(\gamma-1)+ \rho \v^2/ 2  + \B^2 / 8\pi $.
Moreover, for ideal MHD holds $\mathbf{E}  = -\v \times \B$, so that in practice the electric field does not appear into the equations.
In the following, we consider two nontrivial well-known problems of classical ideal MHD, by adopting
the ADER-DG-$\mathbb{P}_5$ scheme, supplemented with our \aposteriori WENO3 sub-cell limiter, with the Rusanov Riemann solver. 

\subsubsection{MHD rotor problem} 
Our  first test is the MHD rotor problem sketched in \cite{BalsaraSpicer1999}. The computational domain is 
$\Omega = [-0.6,0.6]\times[-0.6,0.6]$, with an initial mesh on the coarsest level composed 
of $50\times50$ elements. The AMR framework is activated with $\mathfrak{r}=4$ 
and $\ell_\text{max}=2$. In this problem a high density fluid is rotating rapidly with angular velocity $\omega$, embedded in a low density fluid at rest. 
More specifically, the initial conditions are given by
\begin{equation}
\rho=\left\{\begin{array}{cl}
10 & \text{for}\;\; 0\le r\le 0.1; \\ 1 & \text{otherwise};
\end{array}\right.,~~~
\omega=\left\{\begin{array}{cl}
10 & \text{for}\;\; 0\le r\le 0.1; \\ 0 & \text{otherwise};
\end{array}\right.,~~~
{\B} = \left(\begin{array}{c}
2.5 \\ 0 \\ 0
\end{array}\right),~~~
p = 1\,.
\label{eq:MHDrotor_ic}
\end{equation}
Torsional Alfv\'en waves are generated by the spinning rotor and launched into the ambient medium. 
As a consequence, the angular momentum of the rotor is diminishing. 
In order to validate the accuracy of the method, the AMR computation is compared with the maximally refined uniform grid 
composed of $800\times800=640,000$ elements, corresponding to a total resolution of $4800 \times 4800 = 23,040,000$ 
spatial degrees of freedom on the uniform grid for the augmented MHD equations.
Transmissive boundary  conditions are applied at the borders.
Following \cite{BalsaraSpicer1999}, a linear taper is applied in the range $0.1 \le r \le 0.105$ to allow
continuity of the physical variables between the internal rotor and the fluid at rest at $r=0.105$. The divergence cleaning velocity is set equal to $c_h=4$, 
while the adiabatic index is $\gamma=1.4$.

\begin{figure}
\begin{center}
\begin{tabular}{lr}
\includegraphics[width=0.33\textwidth]{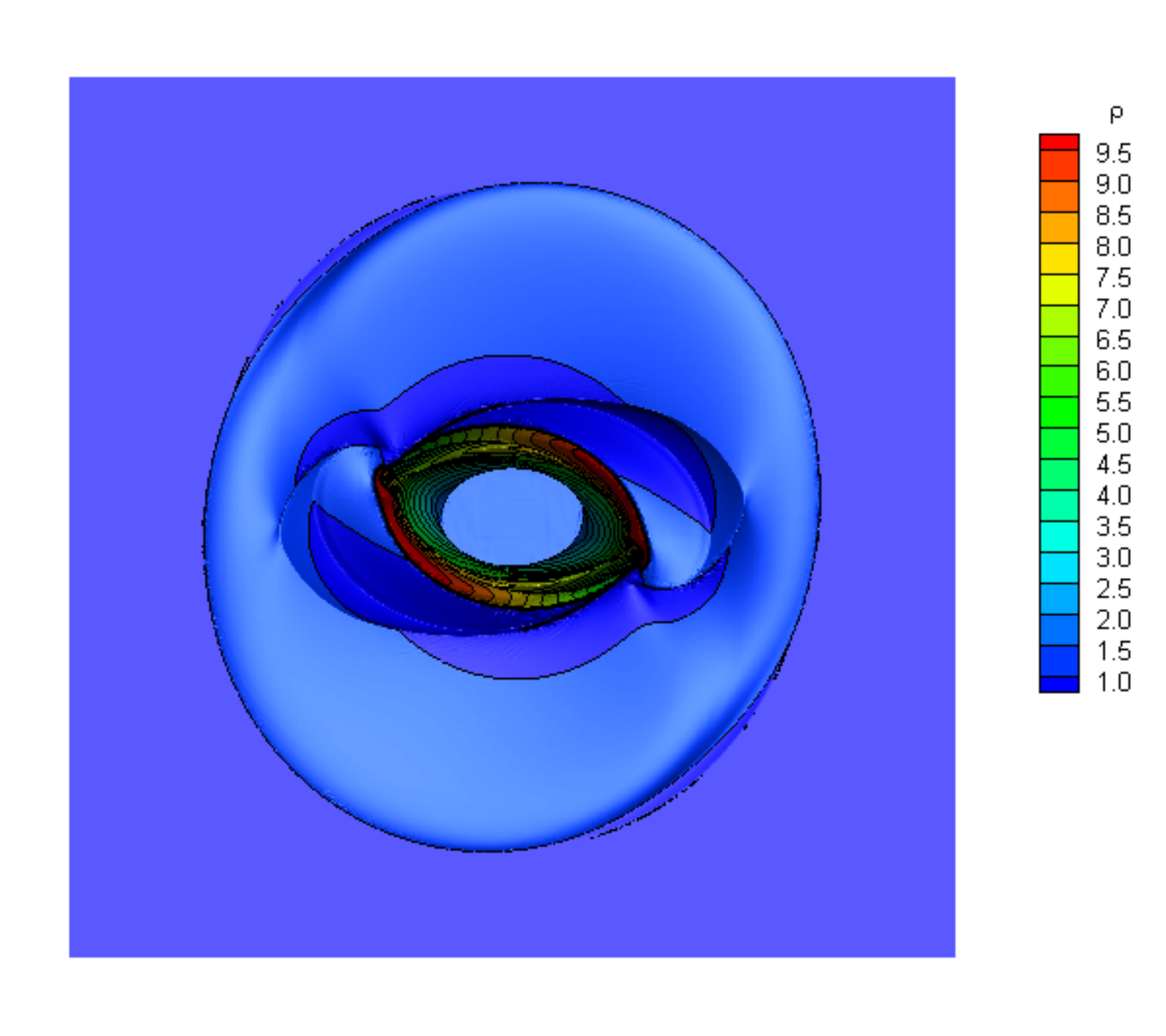}   &  
\includegraphics[width=0.33\textwidth]{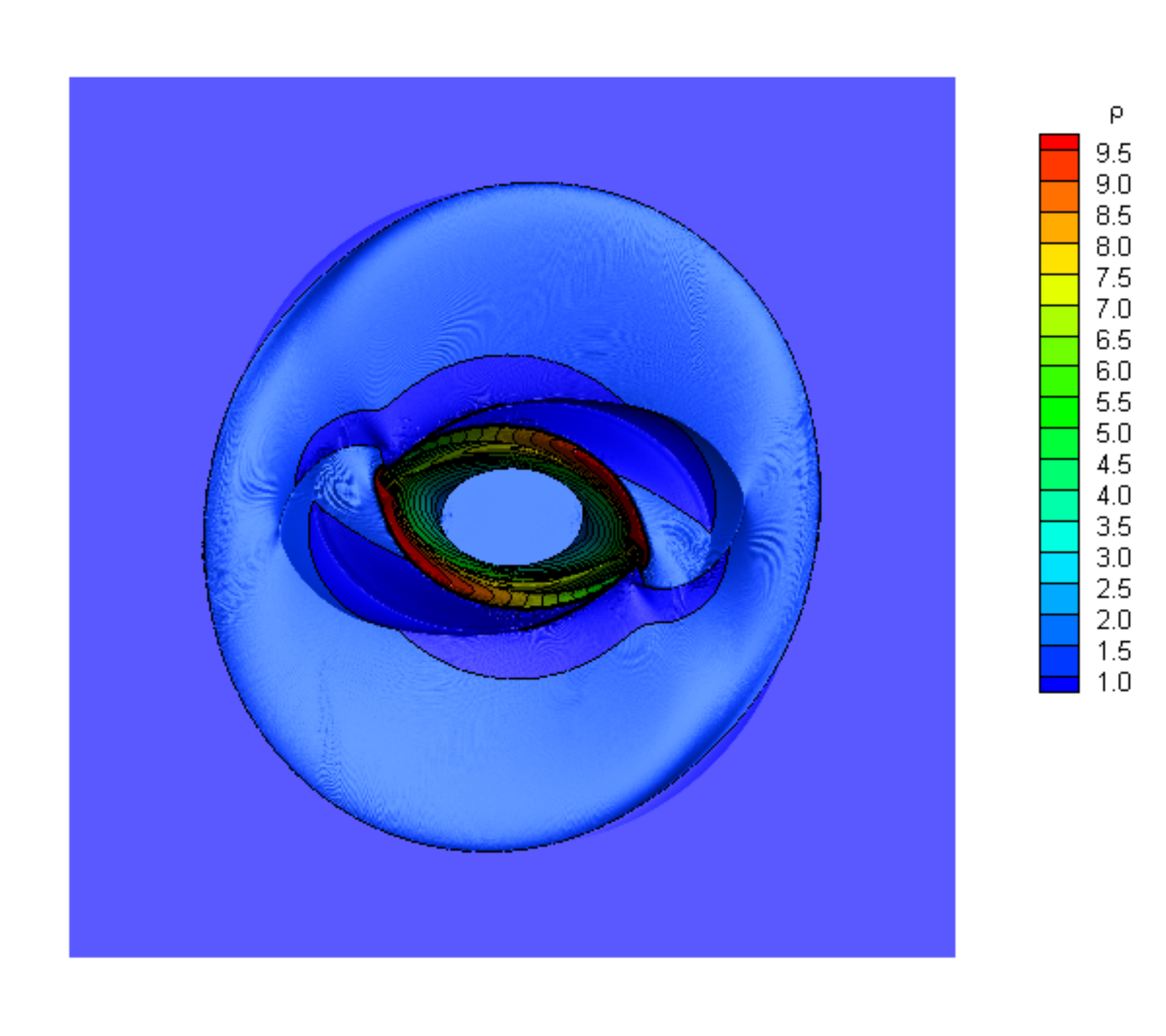}\\
\includegraphics[width=0.33\textwidth]{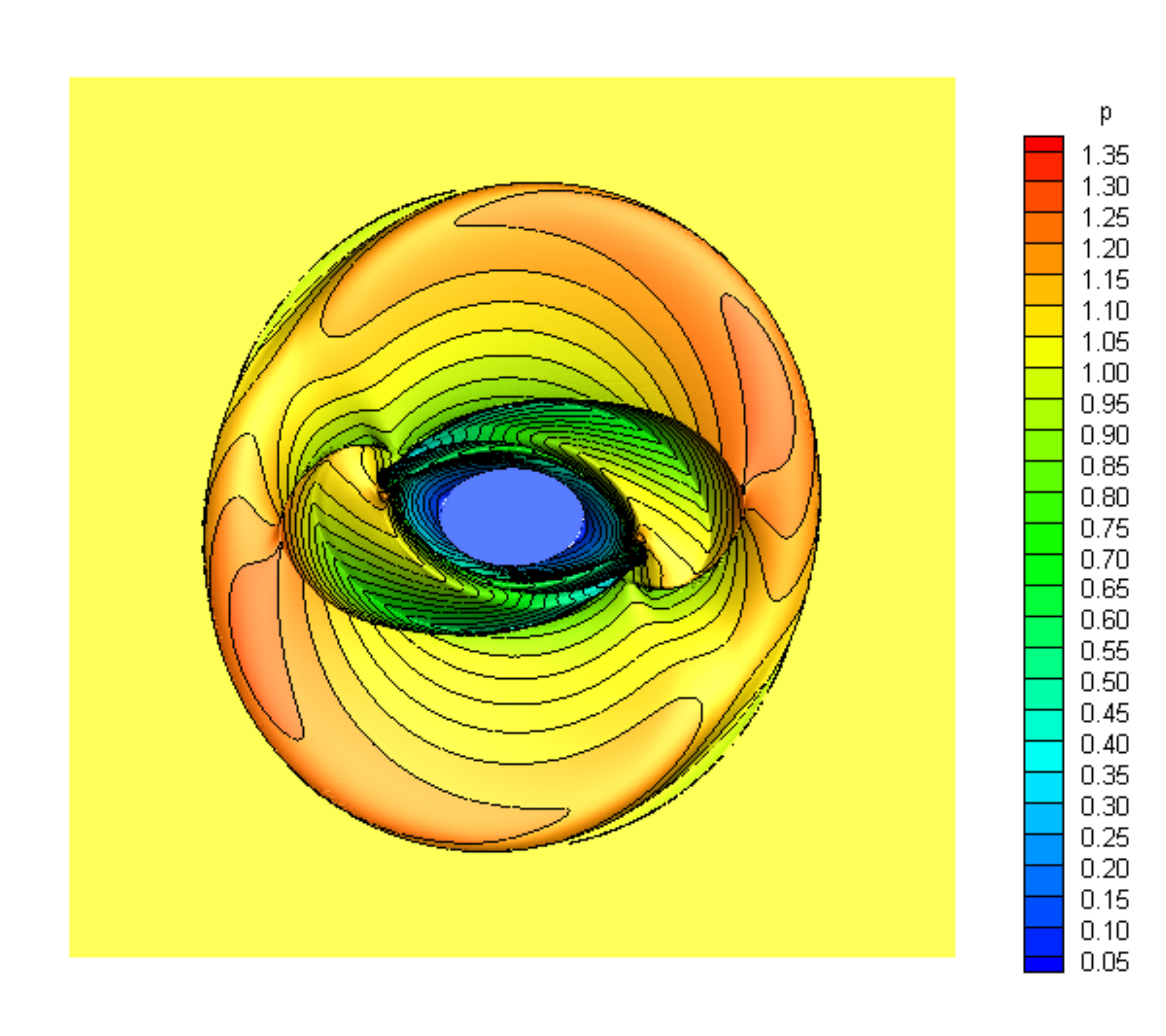}     &  
\includegraphics[width=0.33\textwidth]{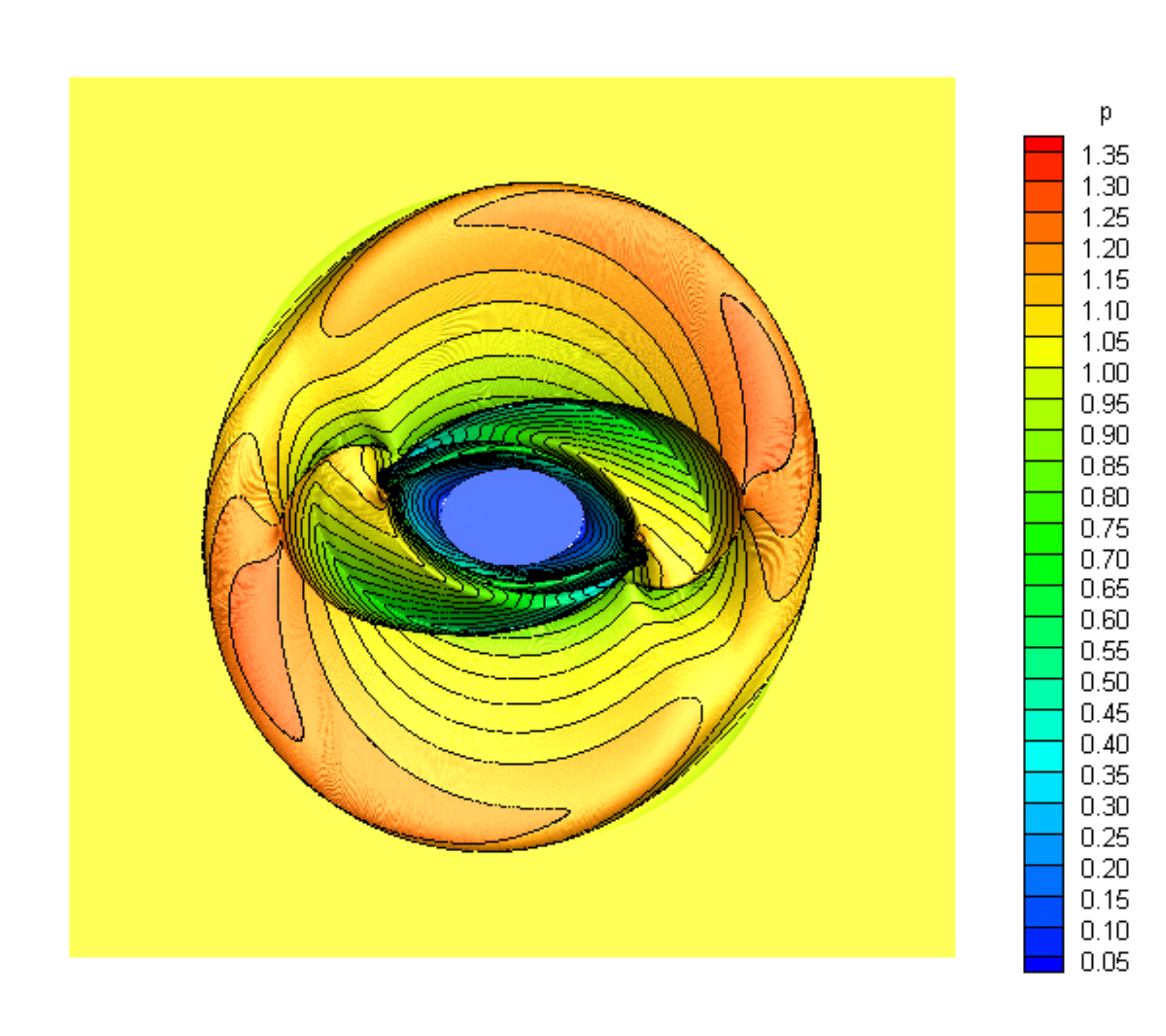}\\
\includegraphics[width=0.33\textwidth]{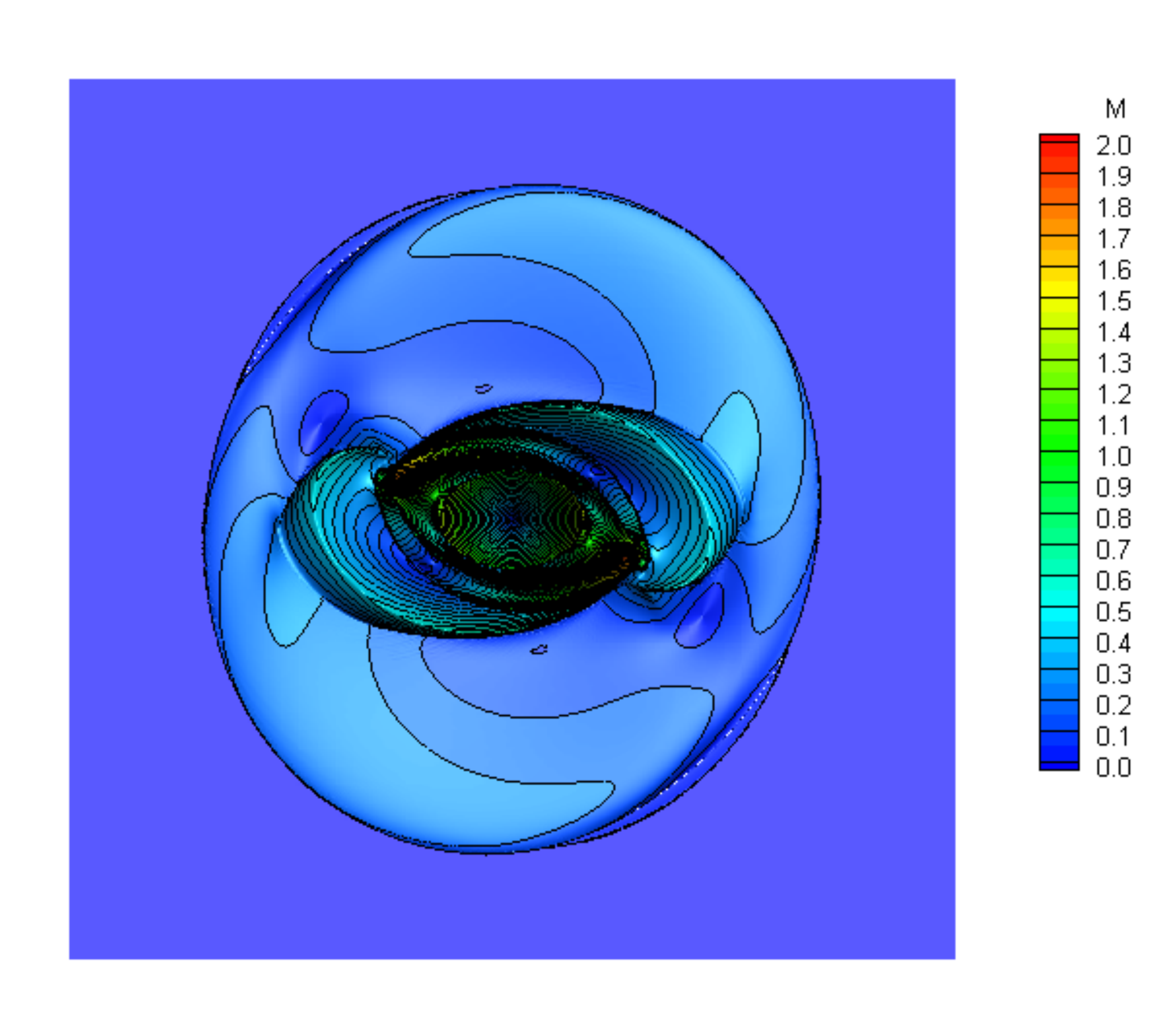}     &  
\includegraphics[width=0.33\textwidth]{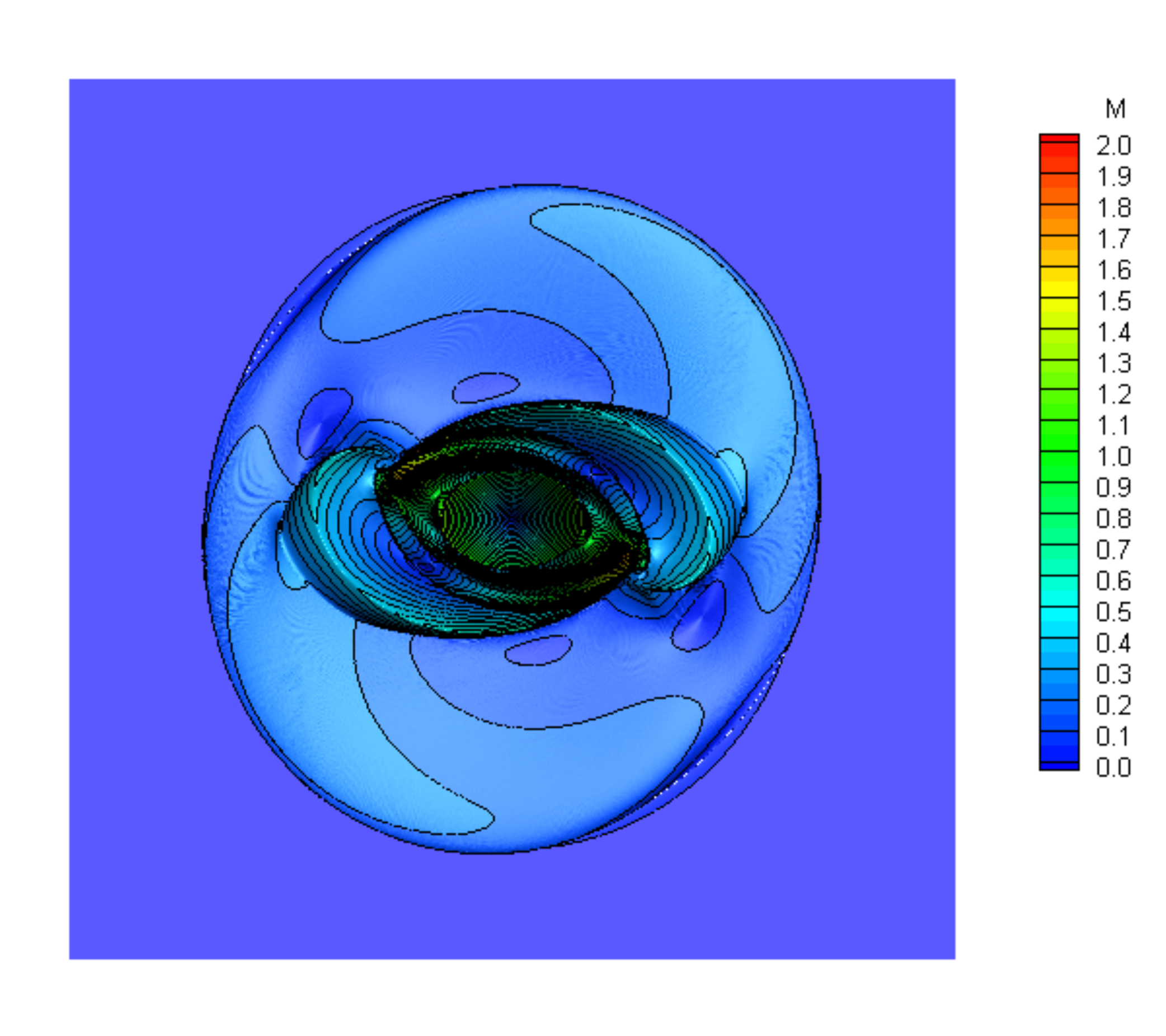}\\
\includegraphics[width=0.33\textwidth]{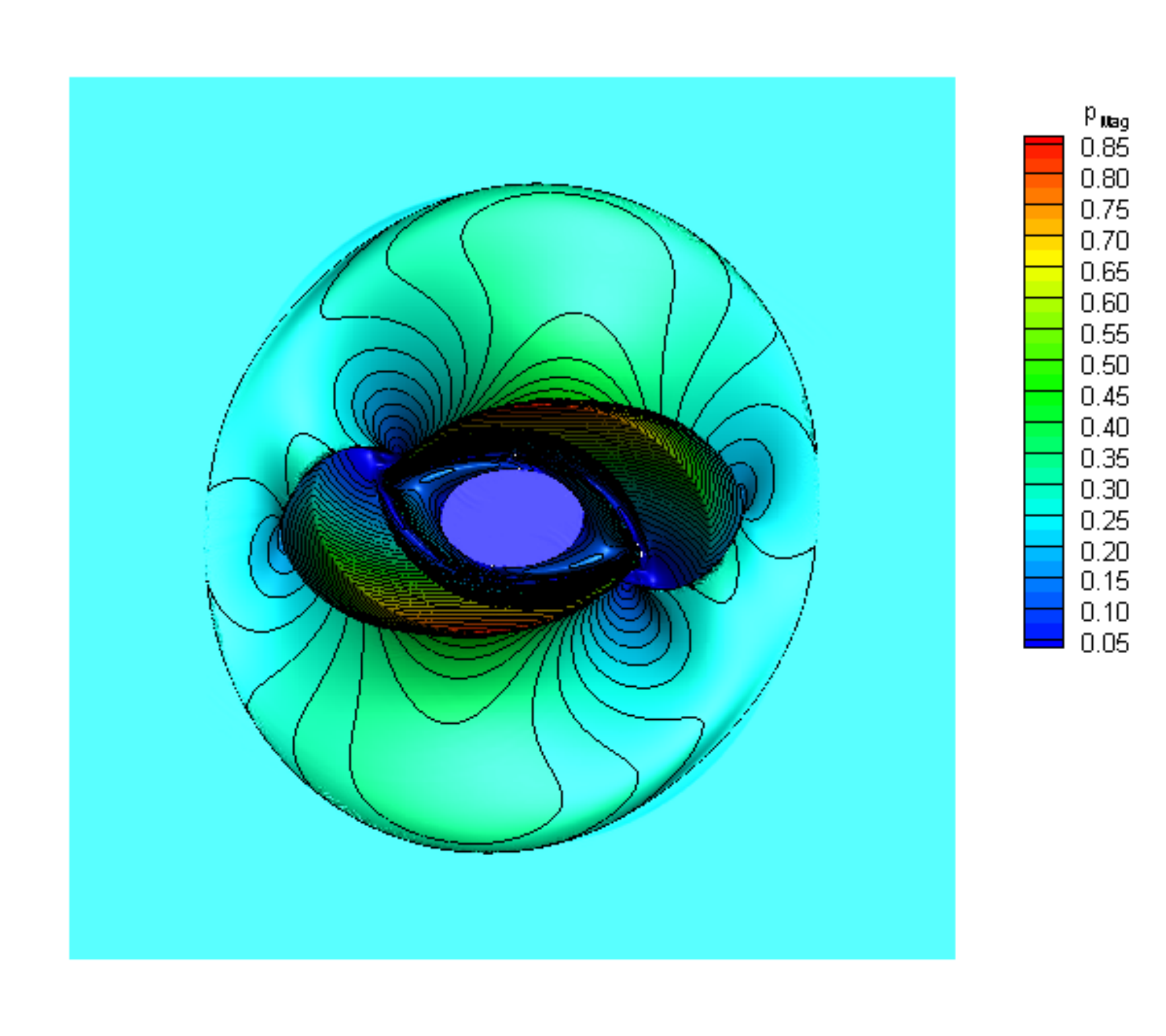}  &  
\includegraphics[width=0.33\textwidth]{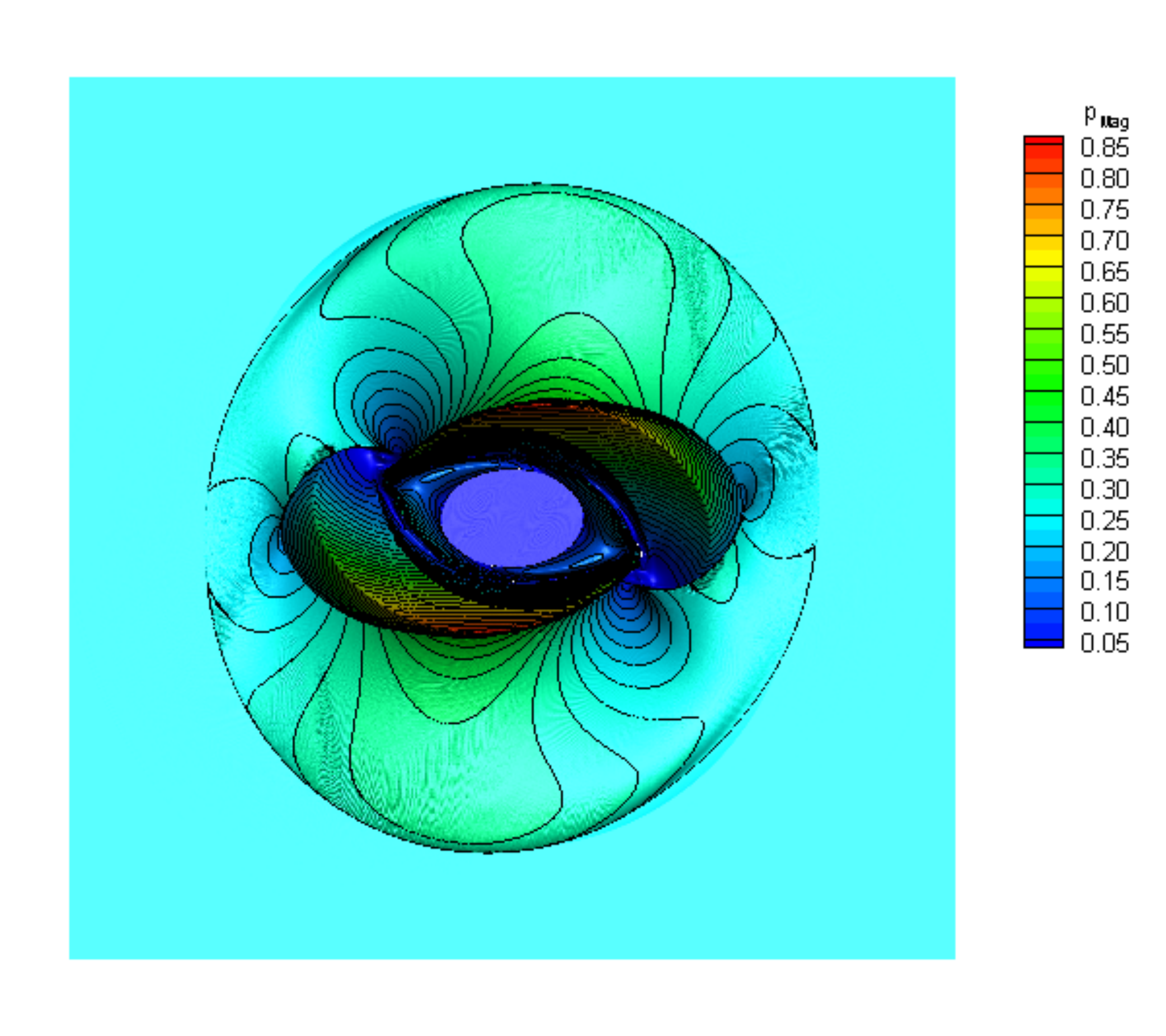}\\
\end{tabular} 
\caption{MHD rotor problem at time $t=0.25$ solved with ADER-DG-$\mathbb{P}_5$.
Left panels: solution obtained
on the AMR grid. Right panels: solution obtained on a fine uniform grid corresponding to the finest AMR grid level.}
\label{fig:MHDrotor}
\end{center}
\end{figure}

\begin{figure}
\begin{center}
\begin{tabular}{lr}
\includegraphics[width=0.45\textwidth]{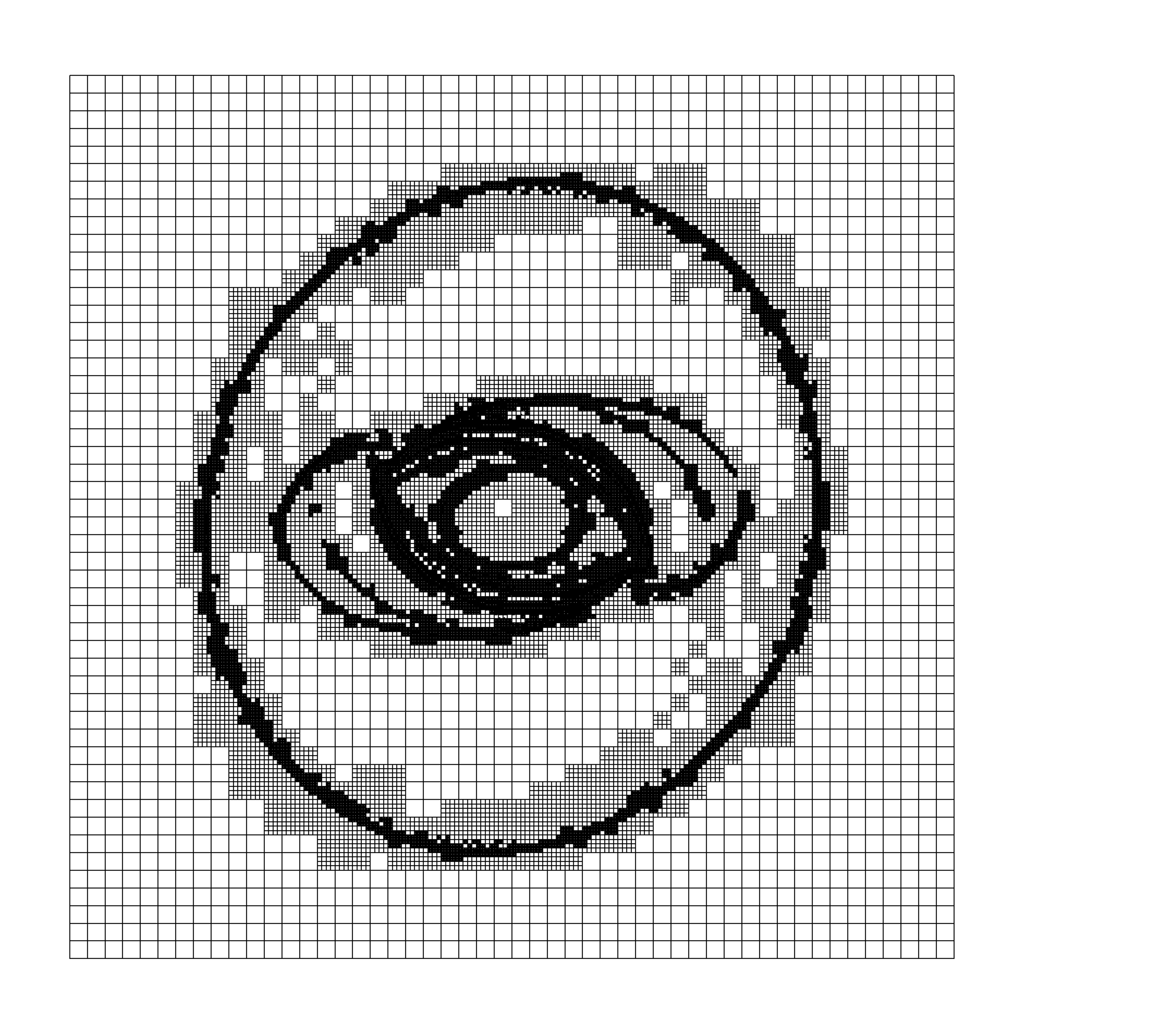}     &  
\includegraphics[width=0.45\textwidth]{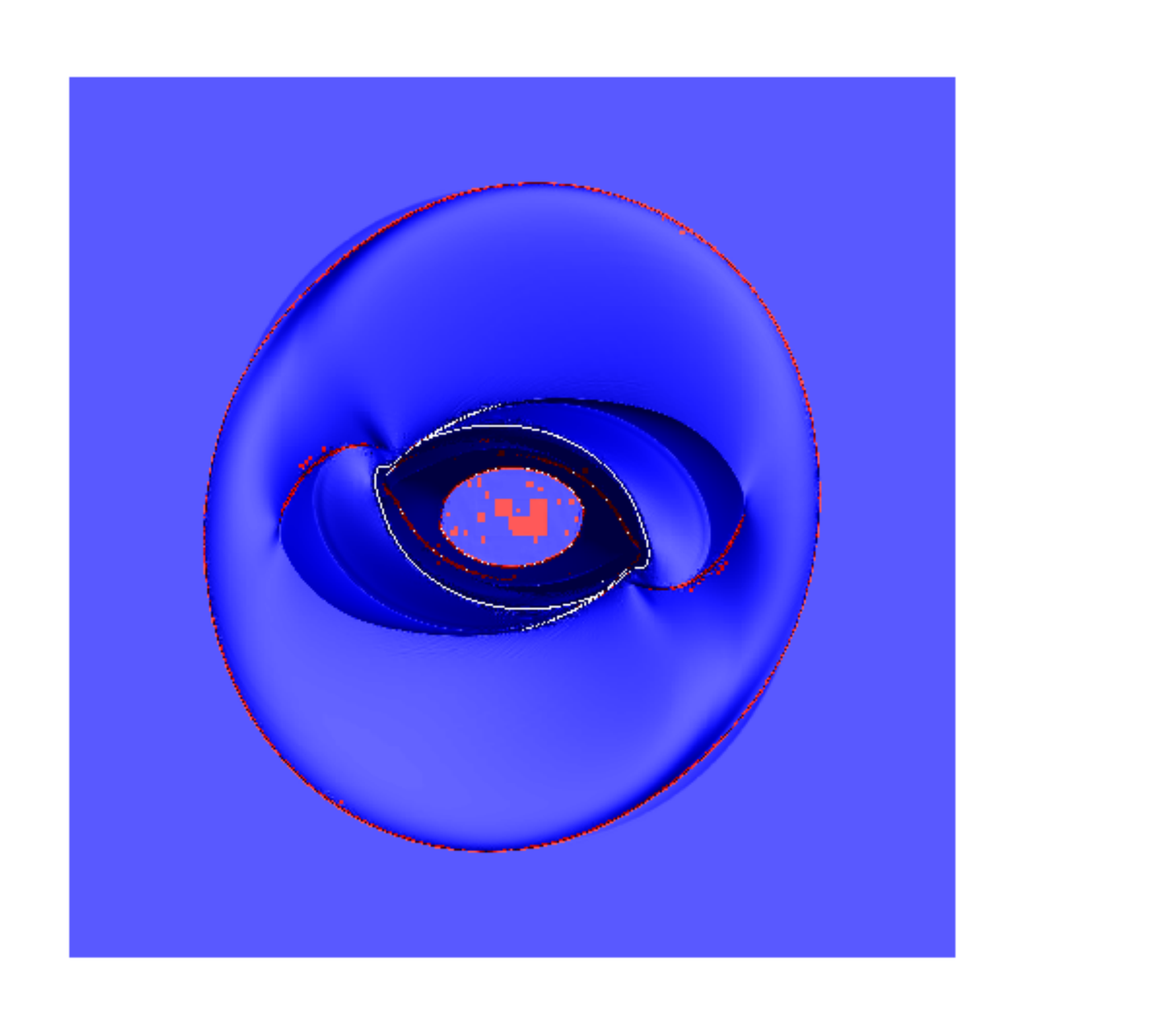}
\end{tabular} 
\caption{MHD rotor problem at time $t=0.25$: AMR grid on the left; troubled cells (red) and unlimited cells (blue) on the right.}
\label{fig:MHDrotor2}
\end{center}
\end{figure}

Figure \ref{fig:MHDrotor} shows the solution for density, pressure, Mach number and magnetic pressure fields at time $t=0.25$.
An excellent agreement between the AMR computation (reported in the left panels) and the uniform grid computation (reported in the right panels) is observed.
Moreover, the numerical results are in very good agreement both with \cite{BalsaraSpicer1999}, and with the results of the ADER-WENO scheme with space-time adaptive mesh refinement presented in \cite{AMR3DCL}. We would like to stress that spurious oscillations
are absent in the density and the magnetic pressure fields, because of the adopted divergence cleaning procedure. In fact, without divergence cleaning, Godunov's schemes would suffer of unphysical oscillations as reported by \cite{BalsaraSpicer1999}. Finally, Fig.~\ref{fig:MHDrotor2} shows the AMR mesh in the left panel and in the right panel the troubled zones in red, for which activation of 
the subcell limiter became necessary. 
\subsubsection{Orszag-Tang vortex system} 
%

\begin{figure}
\begin{center}
\begin{tabular}{lll}
\includegraphics[width=0.32\textwidth]{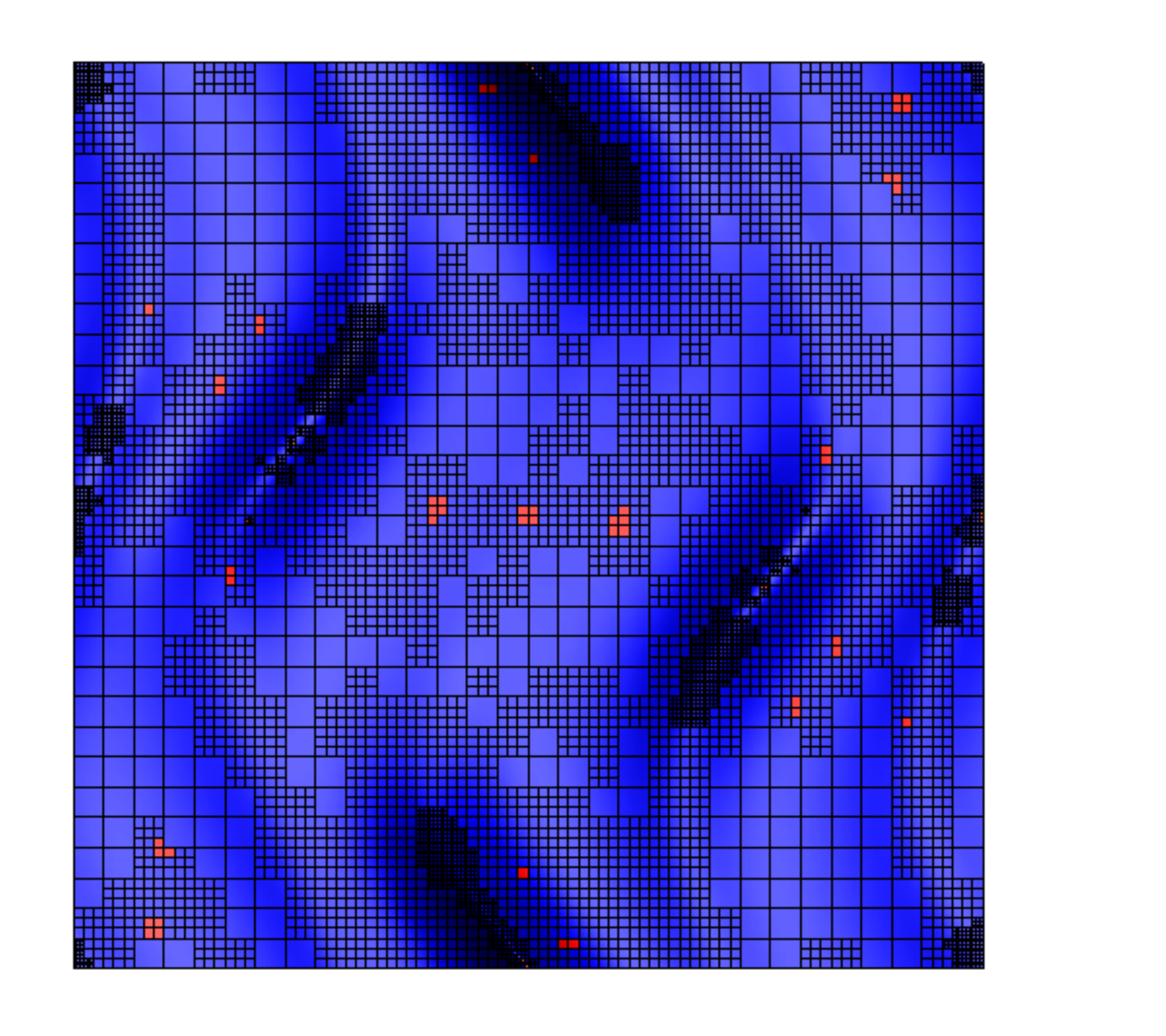}     &  
\includegraphics[width=0.32\textwidth]{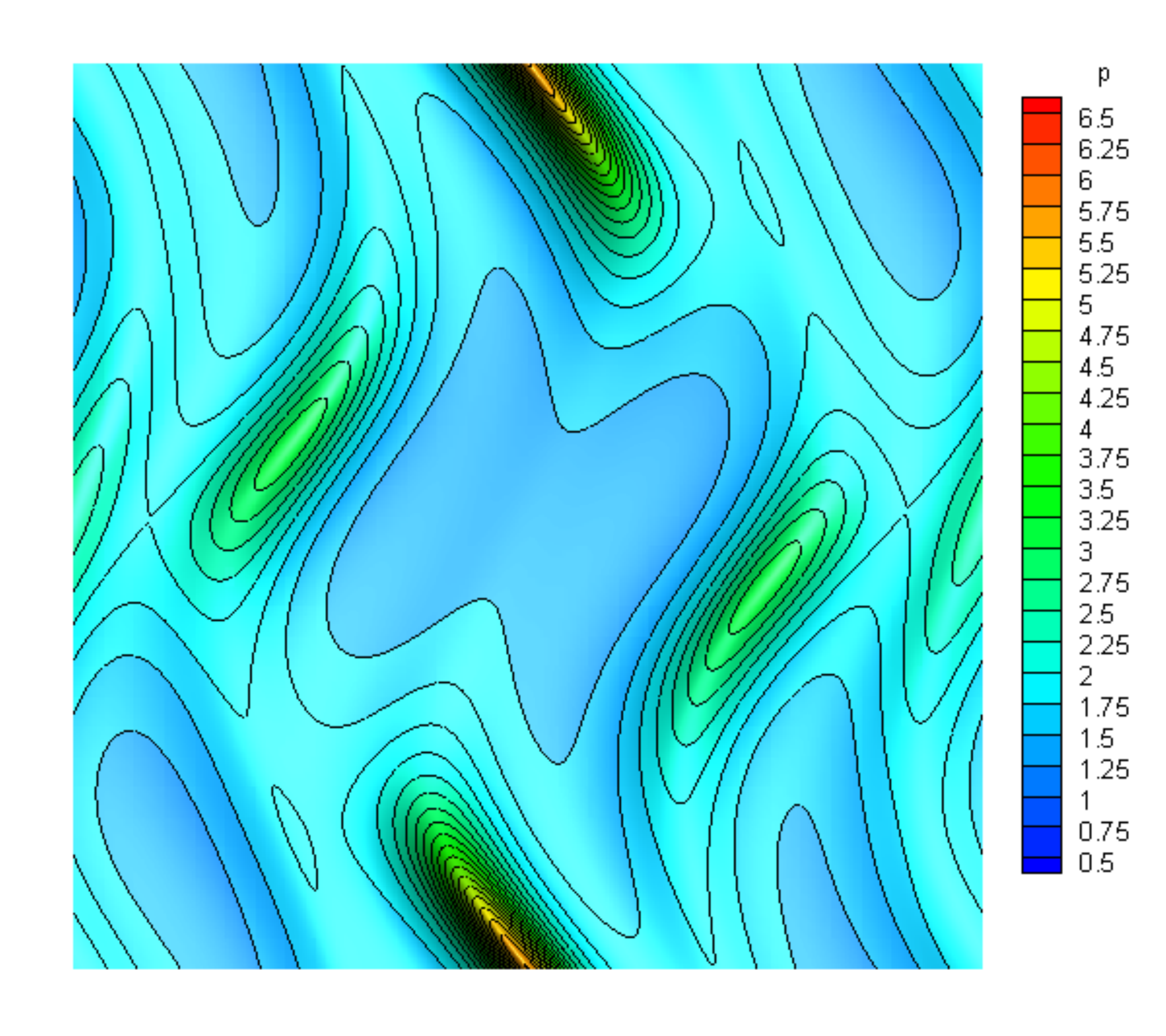}						 &  
\includegraphics[width=0.32\textwidth]{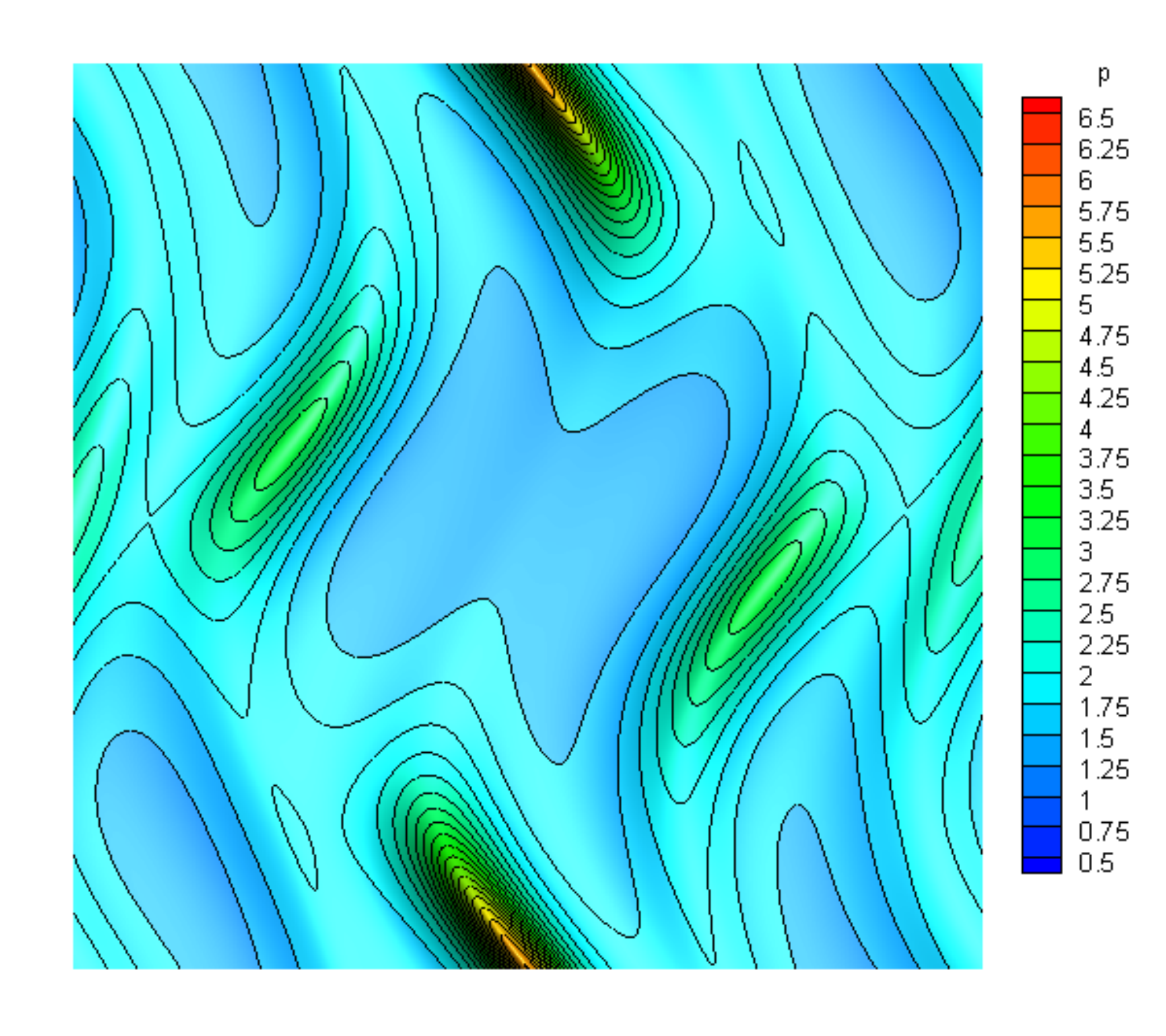}\\
\includegraphics[width=0.32\textwidth]{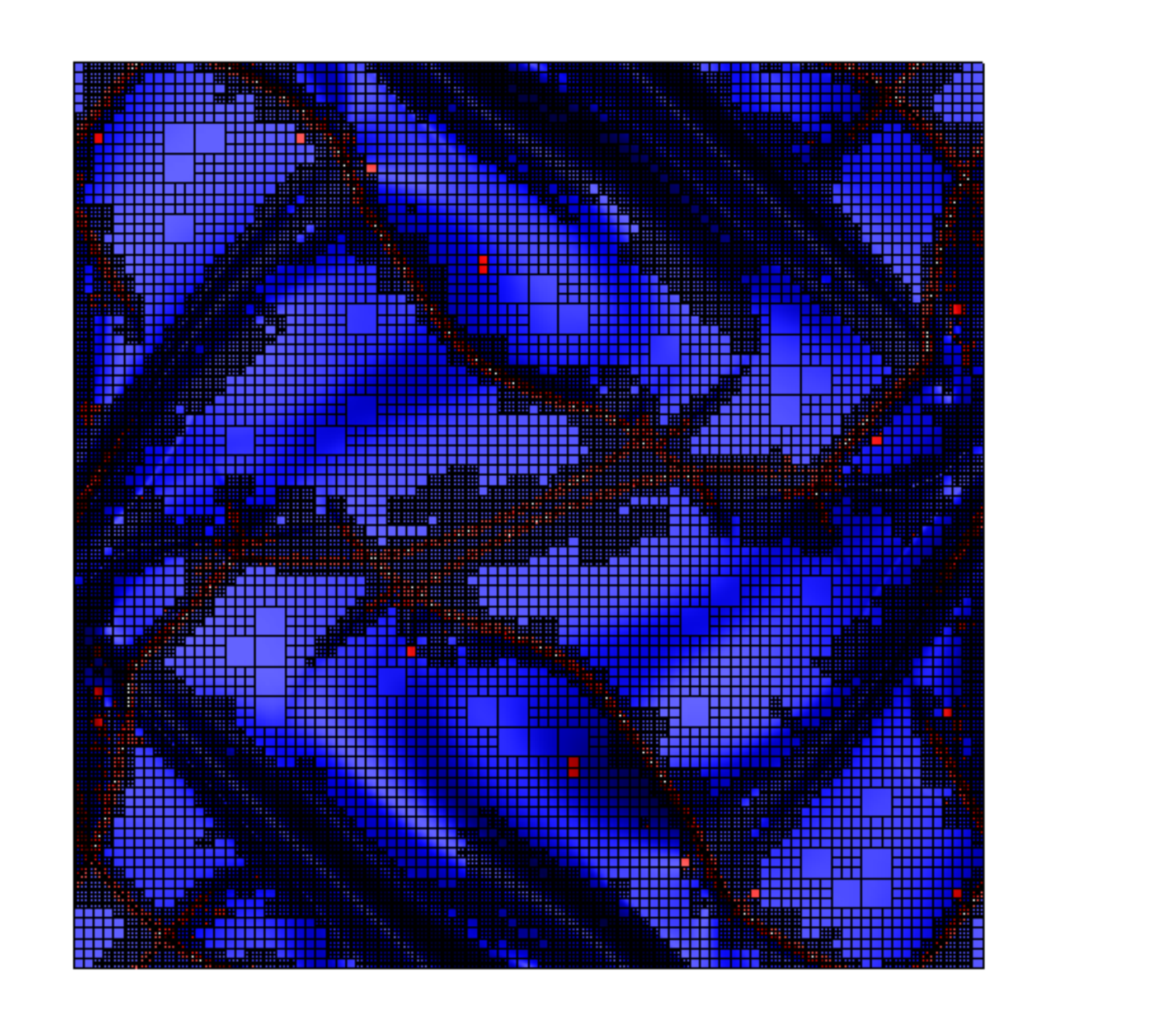}     &  
\includegraphics[width=0.32\textwidth]{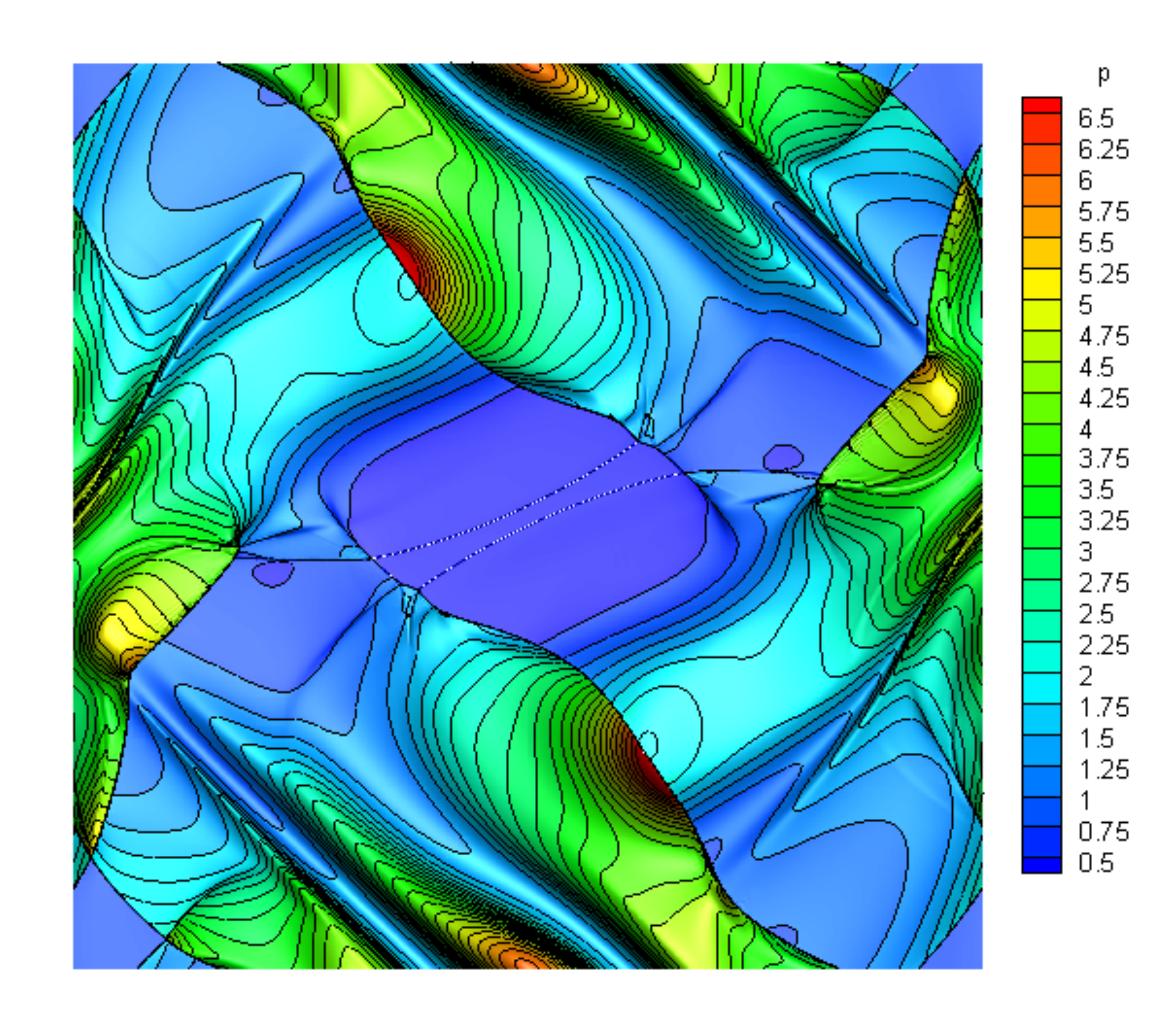}						 &  
\includegraphics[width=0.32\textwidth]{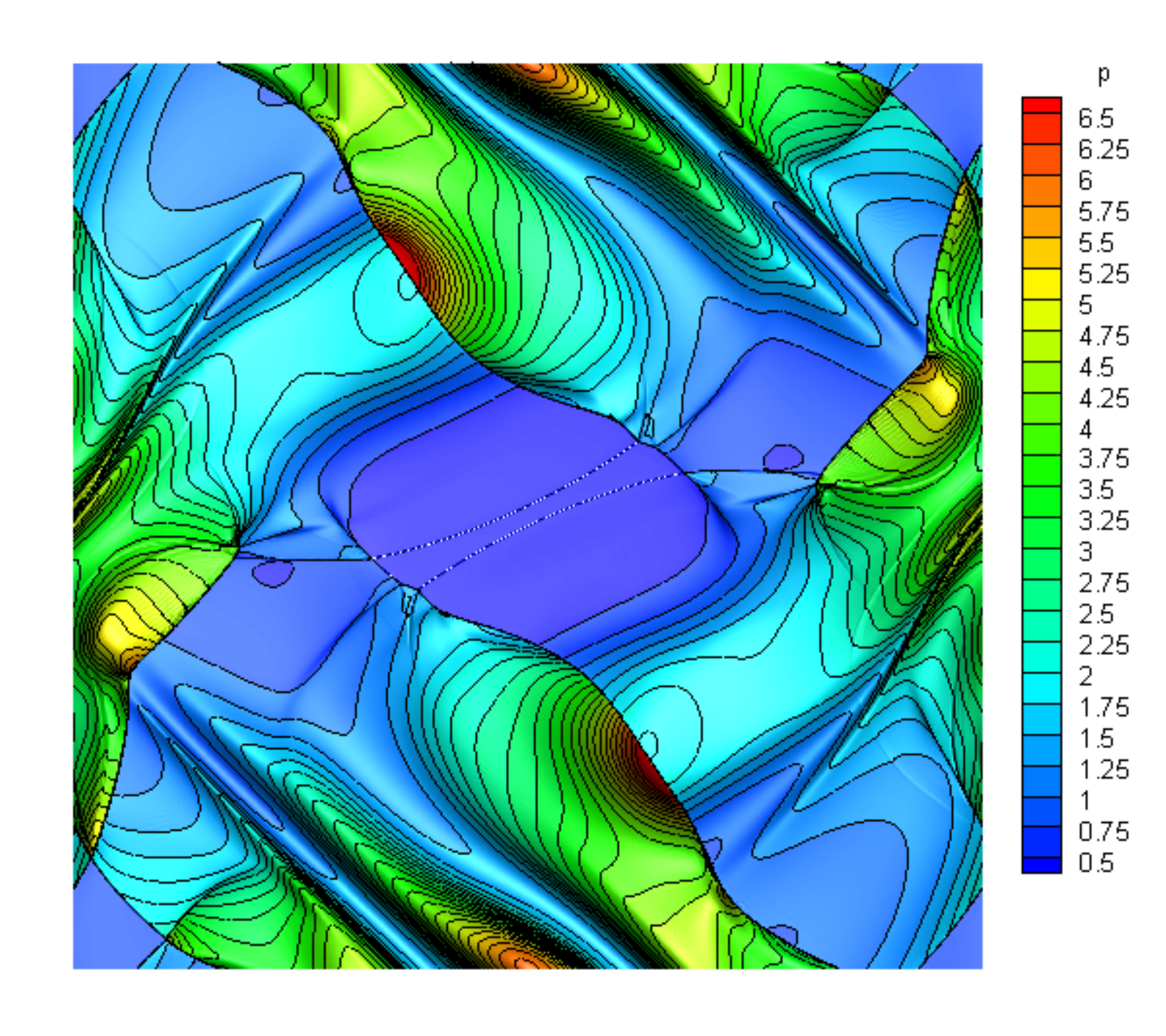}\\
\includegraphics[width=0.32\textwidth]{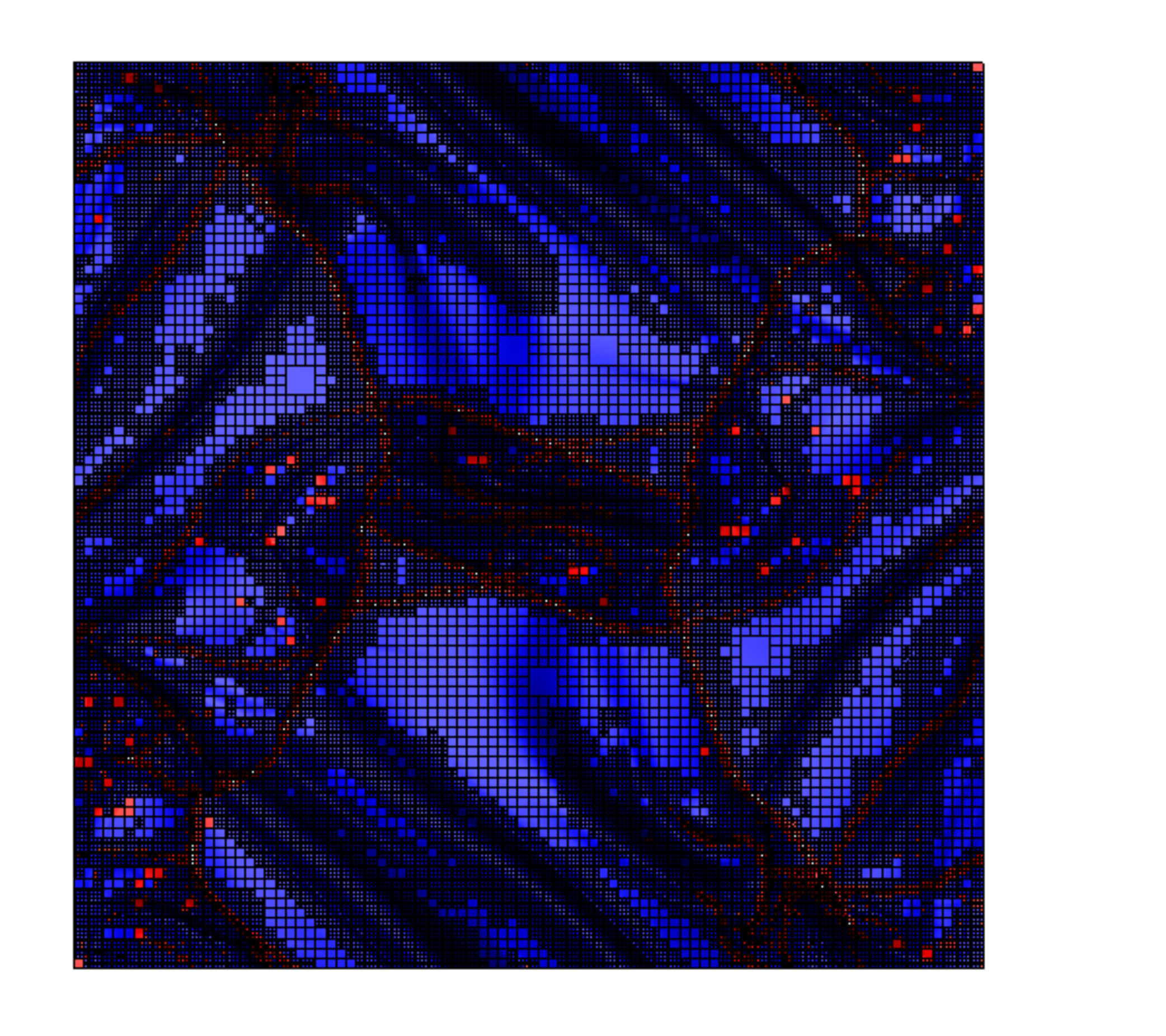}     &  
\includegraphics[width=0.32\textwidth]{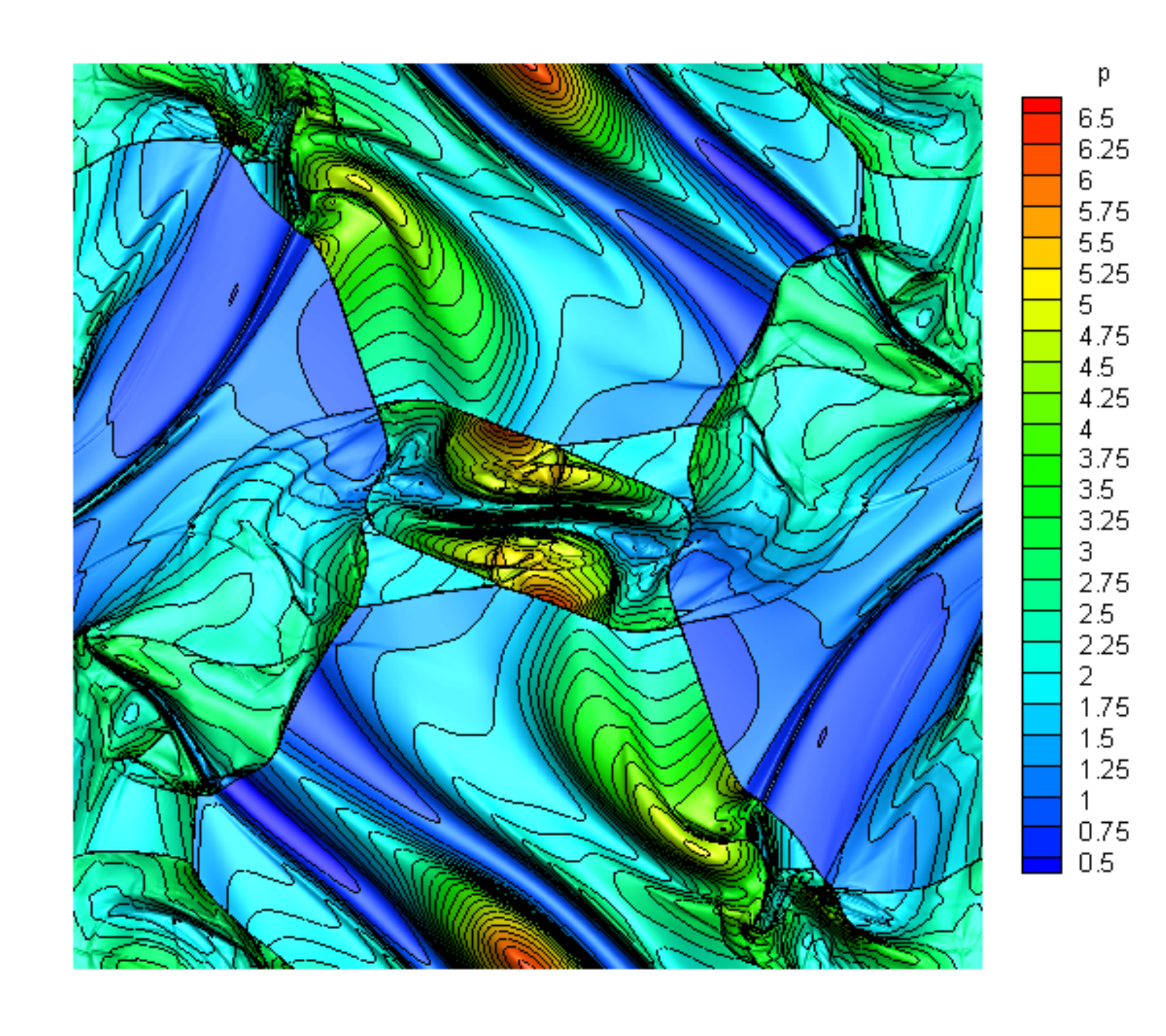}						 &  
\includegraphics[width=0.32\textwidth]{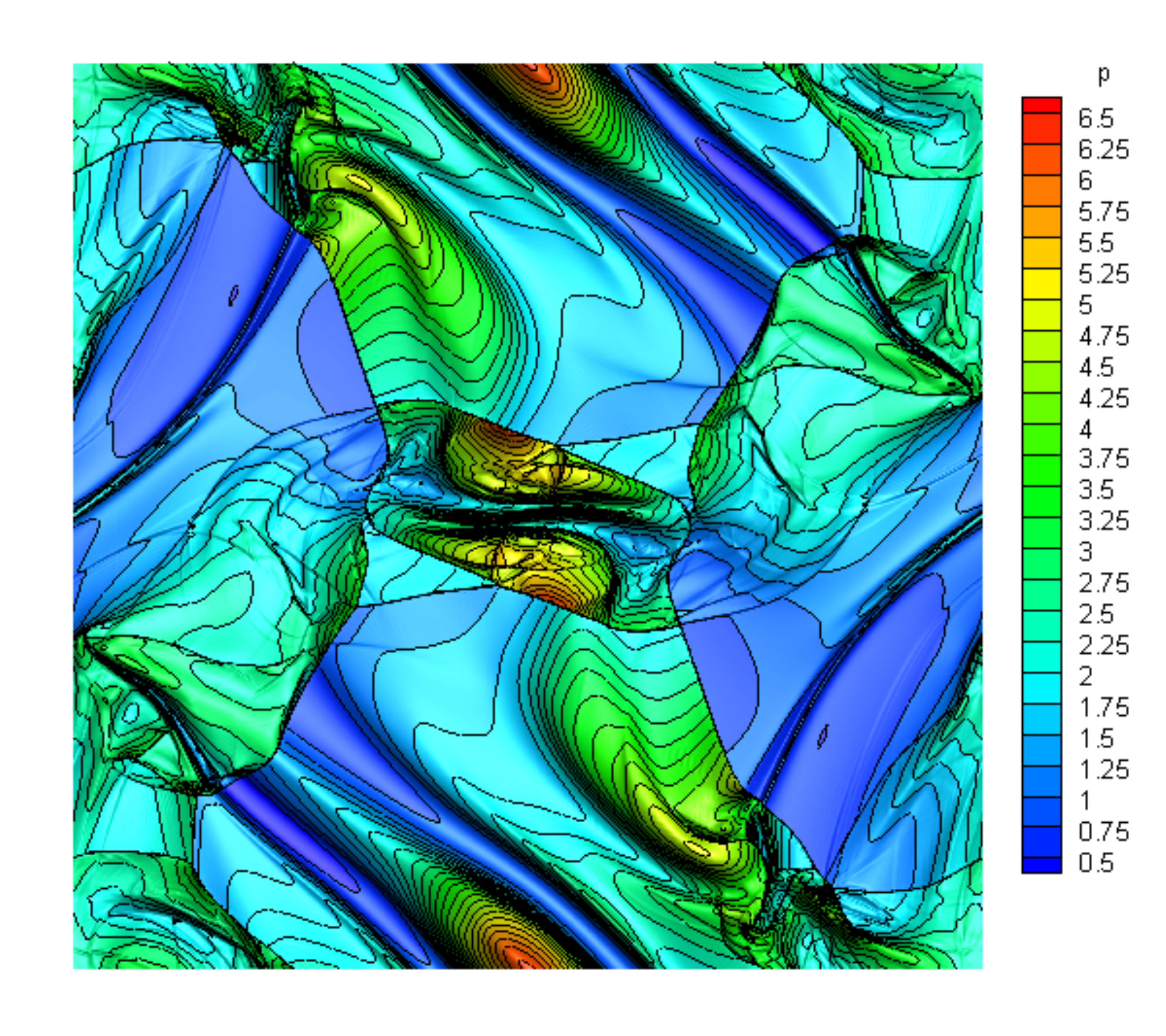}\\
\includegraphics[width=0.32\textwidth]{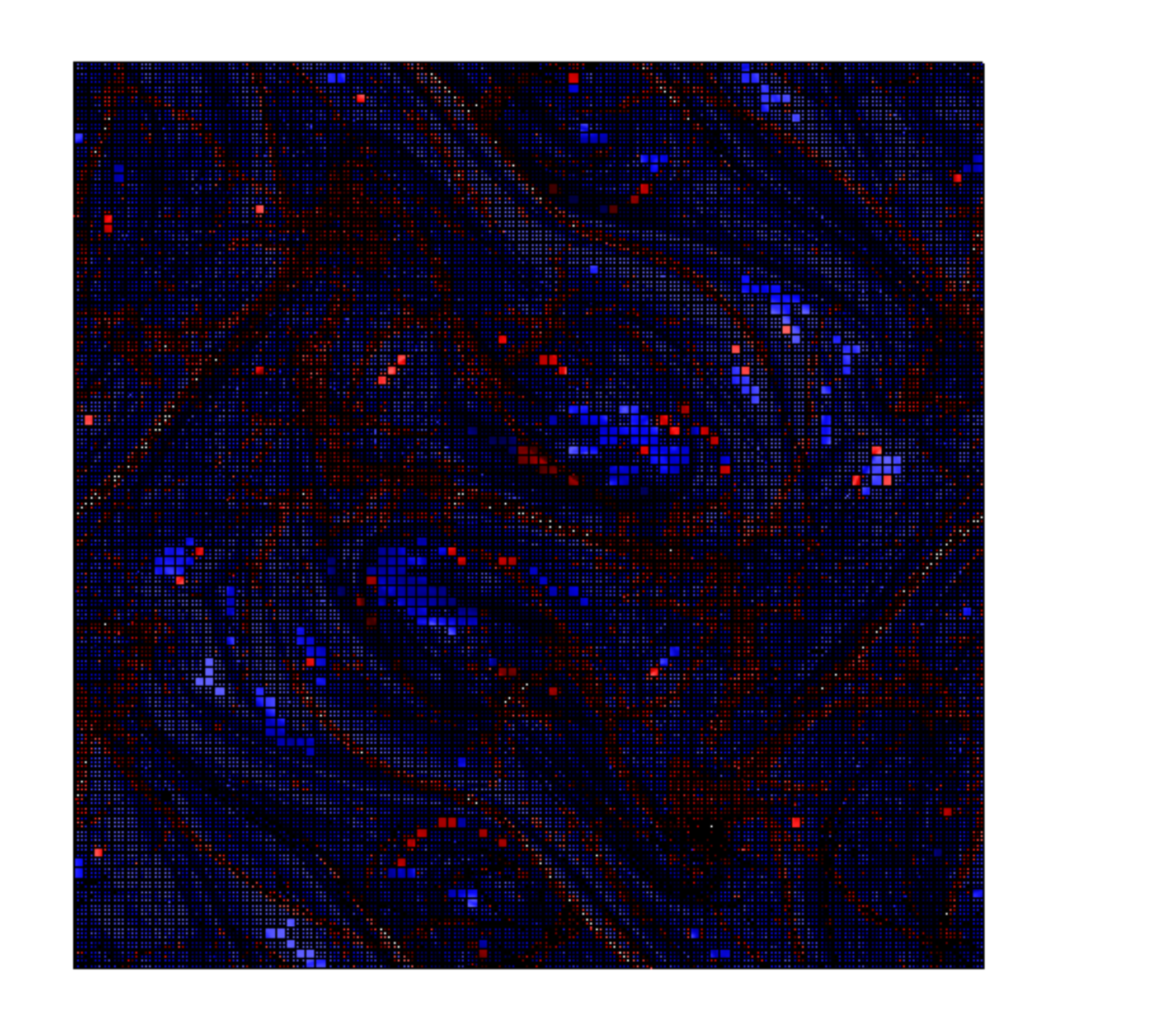}     &  
\includegraphics[width=0.32\textwidth]{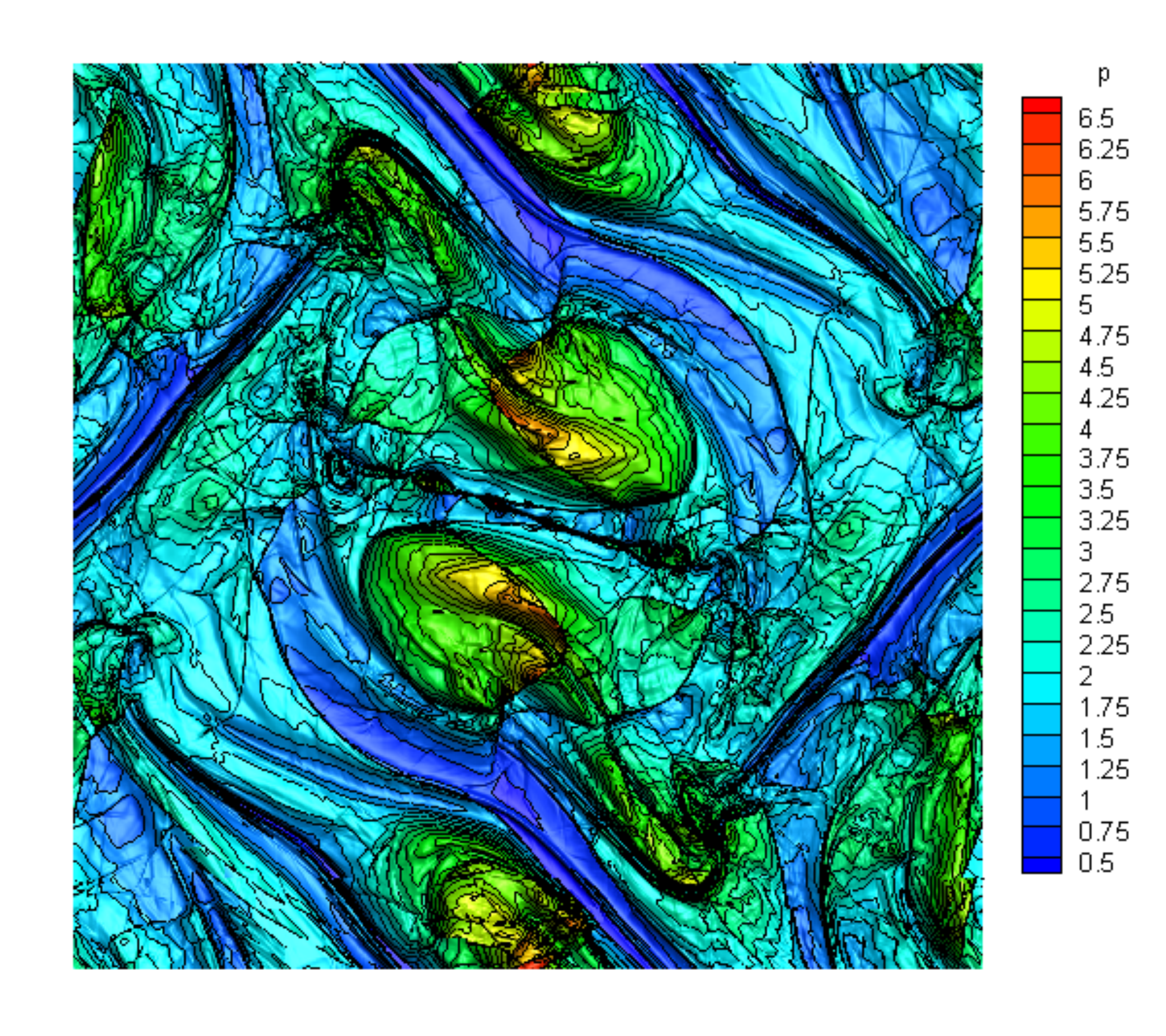}						 &  
\includegraphics[width=0.32\textwidth]{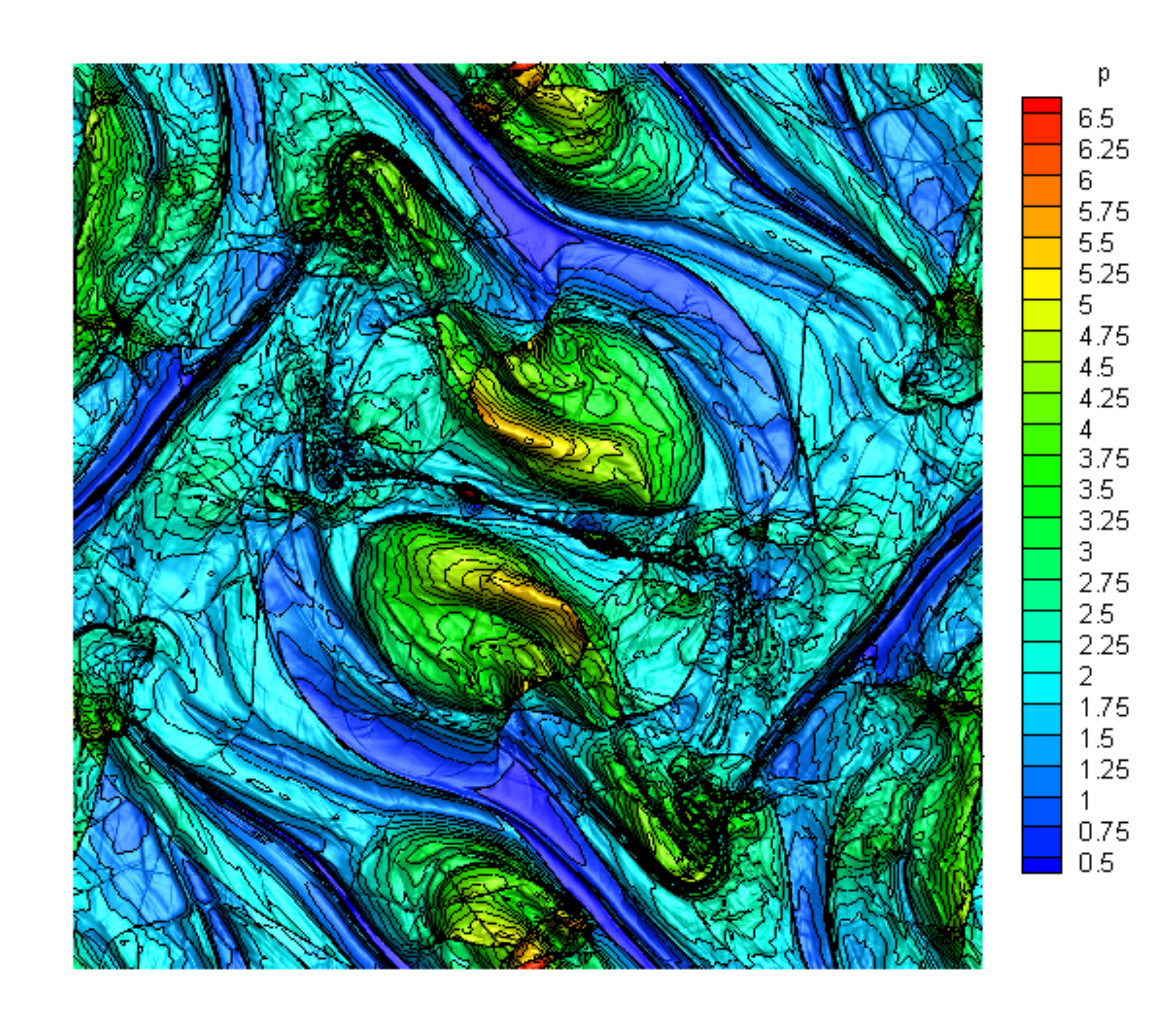}
\end{tabular} 
\caption{Orszag-Tang vortex problem at times $t=0.5$, $t=2.0$, $t=3.0$, $t=5.0$ (from top to bottom) obtained through the ADER-DG-$\mathbb{P}_5$ 
scheme supplemented 
with \aposteriori ADER-WENO3 sub-cell limiter. Left panels: AMR-grid, troubled cells (red) and unlimited cells (blue). Central panels: $\mathbb{P}_5$-solution obtained on the AMR grid. Right panels:  $\mathbb{P}_5$-solution obtained on the uniform grid corresponding to the finest AMR grid level.}
\label{fig:OrszagTang}
\end{center}
\end{figure}

The second test that we have considered concerns the well known Orszag-Tang vortex problem, presented in \cite{OrszagTang}, and later 
investigated by \cite{PiconeDahlburg} and \cite{DahlburgPicone}. 
The adopted parameters refer to the computation performed by \cite{JiangWu}. 
Because of the chosen normalization of the magnetic field, our initial conditions are
\begin{equation}
\left( \rho, u, v, w, p, B_x ,B_y,B_z \right) = \left(  \gamma^2, - \sin\left(y\right), \sin \left(x \right), 0, - \sqrt{4\pi}\sin\left(y\right), \sqrt{4\pi}\sin \left(2x \right), 0\right),
\label{eq:OrszagTang_ic}
\end{equation}
where $\gamma=5/3$. The computational domain $\Omega = [0,2\pi]\times[0,2\pi]$ is discretized with $30\times30$ elements on the coarsest refinement level at $t=0$. Periodic boundary conditions are applied along each edge. 
By using $\mathfrak{r}=3$ and $\ell_\text{max}=2$, the associated maximally refined
uniform grid is formed of $270\times270=72,000$ elements, that correspond to a total resolution of $2,624,400$ spatial degrees of freedom. The resulting solution for density, pressure, Mach number and 
magnetic pressure is plotted at times $t=0.5,2.0,3.0,5.0$ in Fig.\ref{fig:OrszagTang}, both for the AMR and for the uniform grid.  
The AMR results appear to be in very good agreement with the reference solution represented by the calculation over the uniform grid.
Moreover, our computations are in agreement with the the fifth order WENO finite difference results presented by \cite{JiangWu}, with the solution of \cite{Dumbser2008} obtained with an unstructured third order WENO scheme, and also with the ADER-WENO solution computed with space-time adaptive mesh refinement in \cite{AMR3DCL}. 

\section{Conclusion} 
\label{sec:conclusion}

In this paper we have extended the ADER-DG method with \aposteriori ADER-WENO subcell finite volume limiters, which has been recently proposed in \cite{Dumbser2014}, 
to the context of space-time adaptive mesh refinement (AMR) in two and three space dimensions. 
The scheme itself is a modification of the pure discontinuous Galerkin (DG) finite element method that incorporates a novel idea for an \aposteriori limiter. 
In practice, when a cell manifests significant oscillations, which is quite often the case for DG schemes in the presence of discontinuities, a sub-grid composed of $2N+1$  sub-cells is created and the corrupted solution is recomputed. This is done by recovering the solution at the previous time level, by projecting it onto the 
sub-cells to obtain an alternative data representation in terms of sub-cell averages, and, finally, by applying an ADER-WENO evolution in time, so as to completely replace the solution in the cell that was marked as troubled. The interlink with AMR requires a proper communication among the sub-grids at different levels of refinement. In particular, both projection and averaging, 
the two typical AMR operations among refinement levels, 
must be extended to the subcell averages that represent the solution on the sub-grids. 
The new scheme has been validated over a wide sample of test cases for the Euler and for the ideal magnetohydrodynamics equations, both in two and in three spatial dimensions. The nominal order of convergence has been verified up to polynomials of degree $N=8$.
The combination of high order ADER-DG schemes, \aposteriori sub-cell ADER-WENO finite volume limiters within a cell-by-cell AMR framework allows for an unprecedented numerical accuracy. 
In the case of the double Mach reflection problem as well as for the two-dimensional Riemann 
problems that we have considered in this paper, the new numerical approach revealed unexpected details in the dynamics of the fluid. Extending the scheme to equations with stiff source terms 
and true physical dissipative processes is one of the next goals, but even in its present form the method is likely to contribute significantly to those fields of computational fluid dynamics 
where high resolution and low numerical dissipation are needed. Given the unprecedented resolution capabilities shown in various test cases, we are convinced that the method presented in this 
paper belongs to a \textbf{new generation of shock capturing schemes} for computational fluid dynamics. 
Future work may concern the extension of the present space-time adaptive algorithm to high order semi-implicit DG schemes, following the ideas outlined in 
\cite{DumbserCasulli,2STINS,2DSIUSW,BrugnanoCasulli,BrugnanoCasulli2,BrugnanoSestini,CasulliZanolli2010,CasulliZanolli2012,DumbserIbenIoriatti}.

\section*{Acknowledgments}
The research presented in this paper was financed by the European Research Council (ERC) under the
European Union's Seventh Framework Programme (FP7/2007-2013) with the research project \textit{STiMulUs}, 
ERC Grant agreement no. 278267.

The authors are also very grateful for the subsequent financial support of the present research, already granted 
by the European Commission under the H2020-FETHPC-2014 programme with the research project \textit{ExaHyPE}, 
grant agreement no. 671698.  

The authors would like to acknowledge PRACE for awarding access to the SuperMUC 
supercomputer based in Munich, Germany at the Leibniz Rechenzentrum (LRZ). A.H. would also like to thank the 
Spanish Oil Company REPSOL for its support under the project \textit{Development of a numerical multiphase 
flow tool for applications to petroleum industry}.

\bibliographystyle{plain}
\bibliography{./biblio}


\end{document}